\documentclass[10pt]{amsart}
\usepackage{amsmath}
\usepackage{amssymb}
\usepackage{amsfonts}
\usepackage{amsthm}
\usepackage{enumerate}
\usepackage[all]{xy}
\usepackage[latin1]{inputenc}
\usepackage{mathdots}
\usepackage{graphicx}
\usepackage{latexsym}
\usepackage{mathrsfs}
\usepackage{tikz}

\usetikzlibrary{shapes,arrows,positioning,automata}

\definecolor{cof}{RGB}{219,144,71}
\definecolor{pur}{RGB}{186,146,162}
\definecolor{greeo}{RGB}{91,173,69}
\definecolor{greet}{RGB}{52,111,72}

\makeatletter

\newtheoremstyle{definition}
        {5pt}
        {3pt}
        {}
        {0pt}
        {\scshape}
        {.}
        {5pt}
        {\thmname{#1} \thmnumber{#2} \thmnote{[#3]}} 
\newtheoremstyle{theorems}
        {5pt}
        {3pt}
        {\itshape}
        {0pt}
        {\scshape}
        {.}
        {5pt}
        { \thmname{#1} \thmnumber{#2}\thmnote{[#3]}} 

\swapnumbers

\renewcommand\section{\@startsection{section}{1}{\z@}%
        {-3.5ex \@plus -1ex \@minus -.2ex}%
        {2.3ex \@plus .2ex}%
        {\centering\reset@font\scshape}}

\makeatother

\theoremstyle{theorems}

\newtheorem{Theo}{Theorem}[section]
\newtheorem{Prop}[Theo]{Proposition}
\newtheorem{Cor}[Theo]{Corollary}
\newtheorem{Lemma}[Theo]{Lemma}

\theoremstyle{definition}

\newtheorem{Defn}[Theo]{Definition}

\def\mystrut(#1,#2){\vrule height #1pt depth #2pt width 0pt}
\newcommand{\rep}{{\rm rep}}

\newcommand{\Hom}{{\rm Hom}}
\newcommand{\Ext}{{\rm Ext}}
\newcommand{\End}{{\rm End}}

\newcommand{\N}{\mathbb{N}}
\newcommand{\Z}{\mathbb{Z}}
\newcommand{\C}{\mathcal{C}}
\newcommand{\A}{\mathcal{A}}
\newcommand{\D}{\mathcal{D}}
\newcommand{\T}{\mathcal{T}}
\newcommand{\bT}{\,\mathbb{T}}
\newcommand{\B}{\mathcal{B}}

\newcommand{\tauc}{\tau_{\hspace{-0.4pt}_{\mathscr{C}}}\hspace{-0.6pt}}
\newcommand{\mf}{\mathfrak}
\newcommand{\ml}{\hspace{0.5pt}\mathfrak{l}\hspace{0.5pt}}
\newcommand{\mr}{\hspace{0.5pt}\mathfrak{r}\hspace{0.5pt}}
\def\id{\hbox{1\hskip -3pt {\sf I}}}
\def\Ta{\hbox{$\mathit\Theta$}}
\def\Ga{\hbox{$\mathit\Gamma$}}
\def\Da{\hbox{${\mathit\Delta}$}}
\def\Sa{\hbox{${\mathit\Sigma}$}}
\def\Oa{\hbox{${\mathit\Omega}$}}
\def\GaA{\hbox{$\Ga_{\hspace{-2pt}\mathcal{A}}$}}

\begin{document}

\title[Cluster categories]{\sc Cluster categories of type $\mathbb{A}_\infty^\infty$ and \\ triangulations of the infinite strip}

\author[Shiping Liu]{Shiping Liu}

\address{Shiping Liu\\ D\'epartement de math\'ematiques, Universit\'e de Sherbrooke, Sherbrooke, Qu\'ebec, Canada}
\email{shiping.liu@usherbrooke.ca}

\author[Charles Paquette]{Charles Paquette}

\address{Charles Paquette, Department of Mathematics, University of Connecticut, Storrs, CT, 06269-3009, USA}
\email{charles.paquette@usherbrooke.ca}



\thanks{Both authors were supported in part by the Natu\-ral Science and Engineering Research Council of Canada, while the second-named author was also supported in part by the Atlantic Association for Research in the Mathematical Sciences and by the Department of Mathematics at the University of Connecticut.}

\maketitle

\begin{abstract}

We first study the (canonical) orbit category of the bounded derived category of finite dimensional representations of a quiver with no infinite path, and we pay more attention on the case where the quiver is of infinite Dynkin type. In particular, its Auslander-Reiten components are explicitly described. When the quiver is of type $\mathbb{A}_\infty$ or $\mathbb{A}_\infty^\infty$, we show that this orbit category is a cluster category, that is, its cluster-tilting subcategories form a cluster structure as defined in \cite{BIRS}. When the quiver is of type $\mathbb{A}_\infty^\infty,$ we shall give a geometrical description of the cluster structure of the cluster category by using triangulations of the infinite strip in the plane. In particular, we shall show that the cluster-tilting subcategories are precisely given by compact triangulations.

\end{abstract}

\section*{Introduction}

One of the most important developments of the representation theory of quivers is its interaction with the theory of cluster algebras. Cluster algebras were introduced by Fomin and Zelevinsky in connection with dual canonical bases and total positivity in semi-simple Lie groups; see \cite{FoZ0, FoZ}. The two theories are linked together via the notion of cluster categories introduced by Buan, Marsh, Reineke, Reiten and Todorov in \cite{BMRRT}. In its original definition, a cluster category is the orbit category of the bounded derived category of finite dimensional representations of a finite acyclic quiver under an auto-equivalence, which is the composite of the shift functor and the inverse Auslander-Reiten translation.
In such a cluster category, cluster tilting objects correspond to clusters of the cluster algebra defined by the same quiver, and replacing an indecomposable direct summand of a cluster tilting object by another non-isomorphic indecomposable object correspond to mutation of a cluster variable within a cluster. For cluster categories of type $\mathbb{A}_n$, Caldero, Chapoton and Schiffler gave a beautiful geometrical realization in terms of triangulations of the $(n+3)$-gon; see \cite{CCS}. Later on, by replacing cluster tilting objects by cluster tilting subcategories, Buan, Iyama, Reiten and Scott introduced the notion of cluster structure in a $2$-Calabi-Yau triangulated category; see \cite{BIRS}. A canonical example of a $2$-Calabi-Yau category is the orbit category of the bounded derived category of a Hom-finite hereditary abelian category under the composite of the shift functor and the inverse Auslander-Reiten translation, when there is such a translation; see \cite{Kel}. In general, cluster tilting subcategories in a $2$-Calabi-Yau category do not necessarily form a cluster structure. If this is the case, then the $2$-Calabi-Yau category will be called a {\it cluster category}. In \cite{HJo}, Holm and J{\o}rgensen studied a cluster category of infinite Dynkin type $\mathbb{A}_\infty$, which is the derived category of differential graded modules with finite dimensional homology over the polynomial ring viewed as a differential graded algebra, and whose Auslander-Reiten quiver is of shape $\Z\mathbb{A}_\infty$. It was mentioned that this category is equivalent to the canonical orbit category of the bounded derived category of finite dimensional representations of a quiver of type $\mathbb{A}_\infty$ with the zigzag orientation. More importantly, its cluster structure admits a geometrical realization in terms of triangulations of the infinity-gon, a natural generalization of the $\mathbb{A}_n$ case.

\medskip

The aim of this paper is to extend the above-mentioned work of Holm and J{\o}rgensen from a representation theoretic point of view. Indeed, let $Q$ be a locally finite quiver with no infinite path. It is well known that the category $\rep(Q)$ of finite dimensional representations of $Q$ is a Hom-finite hereditary abelian category such that $D^b({\rm rep}(Q))$ has almost split triangles and hence admits a Serre functor; see \cite[(7.11)]{BLP}. Therefore, the canonical orbit category $\mathscr{C}(Q)$ of $D^b(\rep(Q))$, as mentioned above, is a natural candidate for a cluster category of type $Q$. By making use of some results obtained in \cite{BLP}, we shall describe the Auslander-Reiten components of $\mathscr{C}(Q)$; see (\ref{ccpt}) and prove that the projective representations in $\rep(Q)$ generate a cluster-tilting subcategory of $\mathscr{C}(Q)$; see (\ref{Proj-ctsc}). In case $Q$ is of infinite Dynkin type, every indecomposable object of $\mathcal{C}(Q)$ is rigid, and consequently, a cluster-tilting subcategory of $\mathcal{C}(Q)$ is simply a maximal rigid subcategory which is functorially finite; see (\ref{rigidity}). However, we shall only prove that $\mathscr{C}(Q)$ is a cluster category in case $Q$ is of types $\mathbb{A}_\infty$ or $\mathbb{A}_\infty^\infty$; see (\ref{FZ}). Nevertheless, we conjecture that $\mathscr{C}(Q)$ is a cluster category whenever $Q$ is locally finite without infinite paths. In case $Q$ is of type $\mathbb{A}_\infty^\infty$, we shall use an infinite strip $\B_\infty$ with marked points in the plane as a geometric realization of $\mathscr{C}(Q)$. We shall parameterize the indecomposable objects in $\mathscr{C}(Q)$ by arcs in $\B_\infty$ in such a way that two indecomposable objects have no non-trivial extension between them if and only the corresponding arcs do not cross; see (\ref{CrossingPair}). Therefore, maximal rigid subcategories of $\mathscr{C}(Q)$ correspond to triangulations of $\B_\infty$; see (\ref{CrossingPair}), and cluster-tilting subcategories correspond to compact triangulations of $\B_\infty$; see (\ref{functioriallyfinite}). We shall give an easy criterion for a triangulation to be compact; see (\ref{TheoClusterTilting}). In this way, it yields a complete description of the cluster-tilting subcategories. Moreover, it also enables us to count the number of connected components of the quiver of a cluster-tilting subcategory; see (\ref{cpt-cts}). To conclude this introduction, we should mention that triangulations of $\B_\infty$ were first considered in \cite{IgT}, and further studied in \cite{HJo1} as for a geometrical model of a class of cluster categories constructed in a different approach.

\section{Preliminaries}

Throughout this paper, $k$ stands for an algebraically closed field. The standard duality for the category of finite dimensional $k$-vector spaces will be denote by $D$. Throughout this section, $\A$ stands for a Hom-finite Krull-Schmidt additive $k$-category, whose Jacobson radical will be written as ${\rm rad}(\A)$. A {\it radical} morphism in $\A$ is a morphism lying in ${\rm rad}(\A)$. A {\it short cycle} in $\A$ is a sequence $\hspace{-2pt}\xymatrixcolsep{18pt}\xymatrix{X\ar[r] & Y\ar[r] & X}\hspace{-2pt}$ of non-zero radical morphisms between indecomposable objects. Given two collections $\Ta$ and $\Oa$ of objects of $\A$, we write $\Hom_{\A}(\Ta, \Oa)=0$ if $\Hom_{\A}(X, Y)=0$ for all objects $X\in \Ta$ and $Y\in \Oa$. We shall say that $\Ta$ and $\Oa$ are {\it orthogonal} if $\Hom_{\A}(\Ta, \Oa)=0$ and $\Hom_{\A}(\Oa, \Ta)=0.$

\medskip

A morphism $f: M \to N$ in $\A$ is called \emph{right almost split} in $\A$ if it is not a retraction and every non-retraction morphism $g: X \to N$ in $\A$ factors through $f$; and \emph{right minimal} if every factorization $f = fh$ implies that $h$ is an automorphism. If $f$ is right minimal and right almost split in $\A$, then it is called a \emph{sink morphism} for $N$. In dual situations, one says that $f$ is \emph{left almost split, left minimal}, and a {\it source morphism} for $M$. One says that $\A$ has source (respectively, sink) morphisms if there exists a source (respectively, sink) morphism for every indecomposable object of $\A$. An {\it almost split sequence} in $\A$ is a sequence of morphisms
$$\xymatrix{L \ar[r]^f & M\ar[r]^g & N}$$ with $M\ne 0$ such that $f$ is a source morphism and a pseudo-kernel of $g,$ and $g$ is a sink morphism and a pseudo-cokernel of $f$.
Such an almost split sequence is unique for $L$ and for $N$; see \cite{Liu}. Thus, we may write $L=\tau_{_\mathcal{A}}N$ and $N=\tau_{_\mathcal{A}}^-L$, and call $\tau_{_\mathcal{A}}$ the {\it Auslander-Reiten translation} for $\A$. If no confusion is possible, $\tau_{_\mathcal{A}}$ will be simply written as $\tau$.

\medskip

For convenience, we state the following well known result; see, for example, \cite{KSa}.

\medskip

\begin{Lemma}\label{mini-map}

Given any morphism $f: X\to Y$ in $\A$, there exist decompositions

\begin{enumerate}[$(1)$]

\item $f=(f_1, f_2): X=X_1\oplus X_2\to Y$, where $f_1$ is right minimal and $f_2=0;$ and

\vspace{3pt}

\item $f={g_1\choose g_2}: X\to Y_1\oplus Y_2$, where $g_1$ is left minimal and $g_2=0$.

\end{enumerate} \end{Lemma}

\medskip

The {\it Auslander-Reiten quiver} $\GaA$ of $\mathcal{A}$ is a translation quiver endowed with the Auslander-Reiten translation $\tau_{_\mathcal{A}}$, whose underlying quiver is defined as follows; see, for example, \cite{Liu}. The vertex set is a complete set of representatives of the isomorphism classes of indecomposable objects in $\A$. The number of arrows from a vertex $X$ to a vertex $Y$ is the $k$-dimension of
$${\rm irr}(X, Y):={\rm rad}_{\mathcal A}(X, Y)/{\rm rad}^2_{\mathcal A}(X, Y).$$

\medskip

Let $\Ga$ be a connected component of $\Ga_\mathcal{A}$. One defines the {\it path category} $k\Ga$, as well as the mesh category $k(\Ga)$, of $\Ga$ over $k$; see, for example, \cite{Rin}. If $u\in k\Ga$, then its image in $k(\Ga)$ will be denoted by $\bar{u}$. Recall that $\Ga$ is called \emph{standard} if $k(\Ga)$ is equivalent to the full subcategory $\A(\Ga)$ of $\A$ generated by the objects lying in $\Ga$; see \cite{LP, Rin}.

\medskip

\begin{Lemma}\label{para-path}

Let $\Ga$ be a connected component of $\GaA$ of shape $\Z \mathbb{A}_\infty$ or $\Z \mathbb{A}_\infty^\infty$. If $p, q$ are two parallel paths in $\Ga$, then $\bar{p} = \pm \, \bar{q}$ in $k(\Ga)$.

\end{Lemma}

\noindent{\it Proof.} Let $p, q$ be parallel paths from an object $X$ to an object $Y$. Observing that parallel paths in $\Ga$ have the same length, we shall proceed by induction on the length of $p$. The lemma holds trivially if $p$ is of length zero. Suppose that $p$ is of positive length $n$.
If $p,q$ have the same terminal arrow, then the lemma follows immediately from the induction hypothesis. Otherwise, $p =\alpha_1 p_1$ and $q = \alpha_2 q_1$, where $p_1, q_1$ are paths, and $\alpha_1, \alpha_2$ are two distinct arrows in $\Ga$. If $p$ is sectional, then it is the only path in $\Ga$ from $X$ to $Y$, a contradiction. Thus, $p$ is not sectional, and hence, $\Ga$ has a path $u: X\rightsquigarrow \tau Y$. Let $\beta_1, \beta_2$ be the two arrows starting in $\tau Y$ such that $\alpha_1\beta_1 + \alpha_2\beta_2$ is a mesh relation in $k\Ga$. By the induction hypothesis, $\bar{p}_1 = \pm\, \bar{\beta}_1 \bar{u}$ and $\bar{q}_1 = \pm \, \bar{\beta}_2 \bar{u}$. This yields
$$\bar p= \bar{\alpha}_1 \, \bar{p}_1=\pm\, \bar{\alpha}_1 \bar{\beta}_1 \bar{u}= \pm\, \bar{\alpha}_2\bar{\beta}_2 \bar{u}=\pm\, \bar{\alpha}_2 \, \bar{q}_1=\pm\, \bar{q}.$$
The proof of the lemma is completed.

\medskip

Let $\Ga$ be a connected component of $\GaA$ of shape $\Z \mathbb{A}_\infty\vspace{0.5pt}$ or $\Z \mathbb{A}_\infty^\infty$. For each object $X\in \Ga$, we define the {\it forward rectangle} $\mathscr{R}^X$ of $X$ to be the full subquiver of $\Ga$ generated by its successors $Y$ such that, for any path $p: X\rightsquigarrow Y$ and any factorization $p=vu$ with paths $u: X\rightsquigarrow Z$ and $v: Z\rightsquigarrow Y$, either $u$ is sectional, or else, $Z$ has two distinct immediate predecessors. The {\it backward rectangle} $\mathscr{R}_X$ in $\Ga$ of $X$ is defined
in a dual manner.

\medskip

\begin{Prop}\label{rectangle}

Let $\A$ be a Hom-finite Krull-Schmidt additive $k$-category, and let $\Ga$ be a standard component of $\GaA$ of shape $\Z \mathbb{A}_\infty$ or $\Z \mathbb{A}_\infty^\infty\vspace{0.5pt}.$ If $X, Y\in \Ga$, then $\Hom_\A(X, Y)\ne 0$ if and only if $Y\in \mathscr{R}^X$ if and only if $X \in \mathscr{R}_Y.$ In this case, moreover,
$\Hom_\A(X, Y)$ is one-dimensional over $k$.

\end{Prop}

\noindent{\it Proof.} By hypothesis, there exists an equivalence $\varphi: k(\Ga)\to \A(\Ga)$, which acts identically on the objects. For each arrow $\alpha: M \to N \in \Ga$, write $f_\alpha=\varphi(\bar{\alpha})$, which is irreducible in $\A$; see \cite[(1.3)]{LP}. Since $\bar \alpha$ forms a $k$-basis for $\Hom_{k(\mathit\Gamma)}(M,N)$, the set $\{f_\alpha\}$ is a $k$-basis for $\Hom_{\A}(M,N)$. More generally, for a path $p$ in $\Ga$, we write $f_p=\varphi(\bar{p})$.

Let $X, Y\in \Ga$. First of all, it is clear that $Y\in \mathscr{R}^X$ if and only if $X \in \mathscr{R}_Y.$  Hence, we need only to consider forward rectangles. Suppose that $Y \not \in \mathscr{R}^X$. Then either $Y$ is not a successor of $X$ in $\Ga$, or else, $\Ga$ has a path $p: X\rightsquigarrow Y$ which factors through a monomial mesh relation. In the first case, it is clear that $\Hom_{k(\mathit\Gamma)}(X, Y)= 0$. In the second case, $\bar p = 0$, and by Lemma \ref{para-path}, $\bar q =0$ for all paths $q: X\rightsquigarrow Y$. Thus, $\Hom_{k(\mathit\Gamma)}(X, Y)= 0$. As a consequence, $\Hom_\A(X, Y)= 0$.

Suppose now that $Y \in \mathscr{R}^X$. Let $p: X \rightsquigarrow Y$ be a path in $\Ga$. It suffices to prove that $\{f_p\}$ is a $k$-basis for $\Hom_\A(X, Y)$. For this purpose, we proceed by induction on the length of $p$. If $p$ is trivial, then the claim is evident. Otherwise, $p=\alpha q$, where $q: X\rightsquigarrow U$ is a path and $\alpha: U \to Y$ is an arrow. Observe that $U \in \mathscr{R}^X$. By the induction hypothesis, $\{f_q\}$ is a $k$-basis for $\Hom_\A(X, U)\ne 0$. If $p$ is sectional, then $\bar p\ne 0$, and $p$ is the only path in $\Ga$ from $X$ to $Y$. Therefore, $\{\bar p\}$ is a $k$-basis for $\Hom_{k(\mathit\Gamma)}(X, Y)$, and consequently, $\{f_p\}$ is a $k$-basis for $\Hom_\A(X, Y)$. Assume $p$ is not sectional. Then $\Ga$ has a path $w: X\rightsquigarrow \tau Y$. By the induction hypothesis, $\{f_w\}$ is a $k$-basis for $\Hom_{\mathcal A}(X, \tau Y).$ Moreover, by definition, there exists a binomial
mesh relation $\alpha \gamma+\beta\delta$ from $\tau Y$ to $Y$, where $\gamma: \tau Y\to U$, and $\beta: V\to Y$, and $\delta: \tau Y\to V$ are arrows in $\Ga$. Then $\A$ has an almost split sequence
$$\xymatrixcolsep{30pt}\xymatrix{(*) & \tau Y \ar[r]^-{f_\gamma \choose f_\delta} & U\oplus V \ar[r]^-{(f_\alpha, f_\beta)}& Y.}$$

Suppose to the contrary that $f_p=0$. That is, $f_\alpha f_q=0$. In view of the pseudo-exactness of $(*)$, $\A$ has a morphism $u: X\to \tau Y$ such that $f_q=f_\gamma u$ and $f_\delta u=0$. Since $f_q\ne 0$, we see that $u\ne 0$. Again by the induction hypothesis, $\{f_w\}$ is a $k$-basis for $\Hom_{\mathcal A}(X, \tau Y)$ and $\{f_{\delta w}\}$ is a $k$-basis for $\Hom_{\mathcal A}(X, V)$. In particular, $u=\lambda f_w$ for some $\lambda\in k^*$. However, this gives rise to
$\lambda f_{\delta w}= f_\delta (\lambda f_w)=f_\delta u=0$, and hence $f_{\delta w}=0$, a contradiction. Therefore, $f_p\ne 0$. As a consequence, $\bar p\ne 0$. Applying Lemma \ref{para-path}, we deduce easily that $\{\bar{p}\}$ is a $k$-basis for $\Hom_{k(\mathit\Gamma)}(X, Y)$, and therefore, $\{f_p\}$ is a $k$-basis for $\Hom_\A(X, Y)$. The proof of the proposition is completed.

\medskip

If $\Ga$ is a translation quiver of shape $\Z \mathbb{A}_\infty^\infty,\vspace{0.5pt}$ then the forward rectangle of a vertex is the full subquiver generated by its successors. Thus, the following statement is an immediate consequence of Proposition \ref{rectangle}.

\medskip

\begin{Cor}\label{A-d-inf}

Let $\A$ be a Hom-finite Krull-Schmidt additive $k$-category, and let $\Ga$ be a standard component of $\GaA$ of shape $\Z \mathbb{A}_\infty^\infty\vspace{0.5pt}.$ For any objects $X, Y\in \Ga$, we have $\Hom_\A(X, Y)\ne 0$ if and only if $Y$ is a successor of $X$ in $\Ga$.

\end{Cor}

\medskip

Let $\mathcal{D}$ be a full subcategory of $\A$. Recall that $\D$ is {\it covariantly finite} in $\A$ provided that every object $X$ of $\A$ admits a {\it left $\D$-approximation}, that is, a morphism $f: X\to M$ in $\A$ such that every morphism $g: X\to N$ with $N \in \mathcal{D}$ factors through $f$; and {\it contravariantly finite} in $\A$ provided that every object $X$ of $\A$ admits a {\it right $\D$-approximation}, that is, a morphism $f: M\to X$ such that every morphism $g: N\to X$ with $N \in \mathcal{D}$ factors through $f$; and {\it functorially finite} in $\mathcal{A}$ if it is is both covariantly and contravariantly finite in $\A$. Now, we shall say that $\mathcal{D}$ is {\it covariantly bounded} in $\A$ if, for any object $X\in \A$, there exists at most finitely many non-isomorphic indecomposable objects $M\in \mathcal{D}$ such that $\Hom_\A(X, M)\ne 0$; and dually, $\mathcal{D}$ is {\it contravariantly bounded} in $\A$ if, for any object $X\in \A$, there exists at most finitely many non-isomorphic indecomposable objects $M\in \mathcal{D}$ such that $\Hom_\A(M, X)\ne 0$.

\medskip

\begin{Lemma}\label{Fonctorial-finiteness}

Let $\mathcal{D}$ be a full subcategory of $\mathcal{A}.$ If $\mathcal{D}$ is covariantly $($respectively, contravariantly$)$ bounded in $\A$, then it is covariantly $($respectively, contravariantly$)$ finite in $\A$.

\end{Lemma}

\noindent{\it Proof.} We shall prove only one part of the statement. Suppose that $\mathcal{D}$ is contravariantly bounded in $\A$. Let $X$ be an object in $\A$. Let $M_1, \ldots, M_n$ be the non-isomorphic indecomposable objects in $\mathcal{D}$ such that $\Hom_\A(M_i, X)\ne 0$, $i=1, \ldots, n.$ For each $1\le i\le n$, choose a $k$-basis $\{f_{i1}, \ldots, f_{i, d_i}\}$ of $\Hom_\A(M_i, X)$ and set $f_i=(f_{i1}, \ldots, f_{i, d_i}): M_i^{d_i} \to X$, where $ M_i^{d_i}$ denotes the direct sum of $d_i$ copies of $M_i$. Since $\A$ is Krull-Schmidt, it is easy to see that
$$f=(f_1, \cdots, f_n): M_1^{d_1}\oplus\cdots \oplus M_n^{d_n}\to X$$
is a right $\mathcal{D}$-approximation for $X$.  The proof of the lemma is completed.

\medskip

For the rest of this section we assume, in addition, that $\A$ is a triangulated category, whose shift functor is denoted by $[1]$. A {\it Serre functor} for $\A$ is an auto-equivalence $\mathbb{S}$ of $\A$ such that, for any objects $X, Y\in \A$, there exists a natural isomorphism
$$\Hom_{\A}(X, Y)\to D \Hom_{\A}(Y, \mathbb{S}(X)).$$ If such a Serre functor $\mathbb{S}$ exists, then $\mathcal A$ has almost split triangles and its Ausander-Reiten translation is given by the composite of the shift by $-1$ functor followed by $\mathbb{S}$; see \cite{RVDB}. Now, $\mathcal A$ is called 2-{\it Calabi-Yau} if the shift by $2$ functor is a Serre functor; see \cite{BMRRT}. The following easy observation is important for our investigation.

\medskip

\begin{Lemma}\label{2-CY}

Let $\A$ be a 2-Calabi-Yau triangulated $k$-category. For any objects $X, Y$ in $\mathcal{A}$, there exists a natural isomorphism
$$\Hom_\mathcal{A}(X, Y[1])\cong D\Hom_\mathcal{A}(Y, X[1]).$$

\end{Lemma}

\medskip

We say that a full subcategory $\T$ of $\A$ is {\it strictly additive} provided that $\T$ is closed under isomorphisms, finite direct sums, and taking summands. We recall the following definition from \cite{BIRS}, which is our main objective of study.

\medskip

\begin{Defn}

Let $\A$ be a 2-Calabi-Yau triangulated category. A strictly additive category $\T$ of $\A$ is called {\it weakly cluster-tilting} provided that, for any $X \in \A$, $\Hom_{\,\A}(\T, X[1]) = 0$ if and only if $X \in \mathcal{T}$; and \emph{cluster-tilting} provided that $\T$ is weakly cluster-tilting and functorially finite in $\A$.

\end{Defn}

\medskip

Let $\mathcal{T}$ be a strictly additive subcategory of $\A$. In particular, $\T$ is a Krull-Schmidt additive category. By definition, the {\it quiver} of $\T$, denoted by $Q_{_\T}$, is the underlying quiver of its Auslander-Reiten quiver. For each indecomposable object $M$ of $\T,$  we shall denote by $\T\hspace{0.3pt}_M$ the full additive subcategory of $\T$ generated by the indecomposable objects not isomorphic to $M$. Observe that $\T\hspace{0.3pt}_M$ is also strictly additive.

\medskip

\begin{Prop} \label{ct-sc}

Let $\A$ be a $2$-Calabi-Yau triangulated category. If $\T$ is a cluster-tilting subcategory of $\A$, then it has source morphisms and sink morphisms$\,;$ and consequently, its quiver $Q_{_\T}$ is locally finite.

\end{Prop}

\noindent{\it Proof.} Let $\mathcal{T}$ be a cluster-tilting subcategory of $\A$. Suppose that $M$ is an indecomposable object in $\T.$ Then $\mathcal{T}_M$ is functorially finite in $\A$; see \cite[(4.1)]{JoP}.
Let $f: X\to M$ be a right $\T\hspace{0.3pt}_M$-approximation of $M$. By Lemma \ref{mini-map}(1), $f$ restricts to a right minimal morphism $g: Y\to M$, where $Y$ is a direct summand of $X$. Observe that $g$ is a minimal right $\T\hspace{0.3pt}_M$-approximation of $M$.

If ${\rm rad}(\End_{\A}(M))=0$, then $g$ is right almost split, and hence, a sink morphism for $M$ in $\T$. Otherwise, pick a $k$-basis $\{h_1, \ldots, h_m\}$ of ${\rm rad}(\End_{\A}(M))$ and set $h=(h_1, \cdots, h_m): M^m \to M$. We see that every radical endomorphism $M \to M$ factors through $h$. As a consequence, $u=(g,h): Y \oplus M^m \to M$ is right almost split in $\T$. By Lemma \ref{mini-map}(1), $u$ restricts to a right minimal morphism $v: Z\to M$, where $Z$ is a direct summand of $Y \oplus M^m$. Note that $v$ is also right almost split, and hence, a sink morphism for $M$ in $\T$. Dually, we may show that $M$ admits a source morphism in $\T$.
The proof of the proposition is completed.

\medskip

For convenience, we reformulate the notion of a cluster structure without coefficients of a 2-Calabi-Yau triangulated category, originally introduced in \cite{BIRS}.

\medskip

\begin{Defn}\label{Cluster-structure}

Let $\mathcal{A}$ be a 2-Calabi-Yau triangulated category. A non-empty collection $\mathfrak{C}$ of strictly additive subcategories of $\mathcal{A}$ is called a {\it cluster structure} provided that, for any $\T \in \mathfrak{C}$ and any indecomposable object $M\in \T$, the following conditions are verified.

\vspace{-2pt}

\begin{enumerate}[(1)]

\item The quiver of $\T$ contains no oriented cycle of length one or two.

\item There exists a unique (up to isomorphism) indecomposable object $M^* \in \mathcal{A}$ with $M^*\not\cong M$ such that the additive subcategory of $\A$ generated by $\T\hspace{0.3pt}_M$ and $M^*$, written as $\mu_{_M}(\T)$, belongs to $\mathfrak{C}$.

\vspace{1pt}

\item The quiver of $\mu_{_M}(\T)$ is obtained from the quiver of $\mathcal{T}$ by the Fomin-Zelevinsky mutation at $M$ as described in \cite[(1.1)]{FoZ}; see also \cite[Section II]{BIRS}.

\vspace{1pt}

\item There exist two exact triangles in $\mathcal A$ as follows\,: \vspace{-4pt}
$$\xymatrixcolsep{18pt}\xymatrix{M \ar[r]^f & N \ar[r]^g & M^* \ar[r] & M[1]} \mbox{ and }
\xymatrixcolsep{18pt}\xymatrix{M^* \ar[r]^u & L \ar[r]^v& M \ar[r]& M^*[1],}
 \vspace{-2pt}$$
where $f, u$ are minimal left $\T\hspace{0.3pt}_M$-approximations, and $g, v$ are minimal right $\T\hspace{0.3pt}_M$-approximations in $\A$.

\end{enumerate}

\end{Defn}

\medskip

The following notion is our main objective of study.

\medskip

\begin{Defn} \label{Cluster-cat}

Let $\A$ be a 2-Calabi-Yau triangulated category. We shall call $\A$ a {\it cluster category} if its cluster-tilting subcategories form a cluster structure.

\end{Defn}

\section{Representations of quivers}

Throughout this section, let $Q$ stand for a connected locally finite quiver with no infinite path, whose vertex set is written as $Q_0$. By K\"{o}nig's Lemma, the number of paths between any two pre-fixed vertices is finite. Thus, $Q$ is strongly locally finite, and hence, the representation theory of $Q$ obtained in \cite{BLP} applies in this case. For each $x \in Q_0$, let $P_x, I_x,$ and $S_x$ be the indecomposable projective representation, the indecomposable injective representation, and the simple representation at $x$, respectively, which are defined in a canonical way; see, for example, \cite[Section 1]{BLP}. Let $\rep(Q)$ denote the category of finite dimensional $k$-linear representations of $Q$. This is a Hom-finite hereditary abelian category, having almost split sequences; see \cite{BLP, LP}. The Auslander-Reiten quiver $\Ga_{\rep(Q)}$ of $\rep(Q)$ is chosen to contain the representations $P_x$, $I_x$ and $S_x$ with $x \in Q_0$, and its Auslander-Reiten translation is written as $\tau_{\hspace{-1.5pt}_Q}$. The following result is implicitly stated in the proof of Theorem 2.8 in \cite{BLP}.

\medskip

\begin{Lemma}\label{AR-formula}

Let $Q$ be a connected locally finite quiver with no infinite path. If $M, N$ are representations lying in $\Ga_{\rep(Q)}$, then
$$D\Hom_{\hspace{0.4pt}\rep(Q)}(N, \tau_{\hspace{-1pt}_Q}\hspace{-1pt}M)\cong \Ext^1_{\rep(Q)}(M, N)\cong D\Hom_{\hspace{0.4pt}\rep(Q)}(\tau_{\hspace{-1.5pt}_Q}^-N, M).$$

\end{Lemma}

\medskip

Recall that $\Ga_{\rep(Q)}$ has a unique preprojective component $\mathcal{P}$ containing all the $P_x$ with $x\in Q_0$, a unique preinjective component $\mathcal{I}$ containing all the $I_x$ with $x\in Q_0$, and possibly some other regular components which contain none of the $P_x$ and $I_x$ with $x\in Q_0$; see \cite{BLP}.

\medskip

\begin{Theo} \label{structureARquiverGen}

Let $Q$ be a connected locally finite quiver with no infinite path.

\begin{enumerate}[$(1)$]

\item The preprojective component $\mathcal{P}$  of $\Ga_{\rep(Q)}$ is standard of shape $\N\hspace{0.5pt} Q^{\rm op}$.

\vspace{1pt}

\item The preinjective component $\mathcal{I}$  of $\Ga_{\rep(Q)}$
is standard of shape $\N^- \hspace{-1pt}Q^{\rm op}$ with $\Hom_{\hspace{0.4pt}\rep(Q)}(\mathcal{I}, \mathcal{P})= 0$.

\vspace{1pt}

\item If $\mathcal R$ is a regular component of $\Ga_{\rep(Q)},\vspace{1.5pt}$ then it is of shape $\Z \mathbb{A}_\infty$ such that $\Hom_{\hspace{0.4pt}\rep(Q)}(\mathcal{I}, \mathcal{R})= 0$ and $\Hom_{\hspace{0.4pt}\rep(Q)}(\mathcal{R},\mathcal{P})=0.$

\end{enumerate}
\end{Theo}

\noindent{\it Proof.} Let  $f: M \to N$ be a non-zero morphism in $\rep(Q)$, where $M, N\in \Ga_{\rep(Q)}$. Suppose that $M$ is preinjective. Then there exists some $r\ge 0$ for which $\tau^{-r}M=I_x$ for some $x\in Q_0$. If $N$ is not preinjective, then $\tau^{-r}N$ is defined and not injective. On the other hand, by Lemma \ref{AR-formula}, we have a non-zero morphism $g: \tau^{-r}M \to \tau^{-r}N$. Since $\rep(Q)$ is hereditary, $\tau^{-r}N$ is injective, a contradiction. Dually, if  $N$ is  preprojective, then so is $M$. The rest of the theorem has already been established in \cite{LP}. The proof of the theorem is completed.

\medskip

Let $\Ga$ be a regular component of $\Ga_{\rep(Q)}$, and let $X$ be a representation lying in $\Ga$. One says that $X$ is \emph{quasi-simple} if it has only one immediate predecessor in $\Ga$. In general, since $\Ga$ is of shape $\Z \mathbb{A}_\infty,$ it has a unique sectional path
$$\xymatrix{
X=X_n\ar[r]& X_{n-1}\ar[r]& \cdots \ar[r] & X_1
}$$
with $X_1$ being quasi-simple. One defines then the {\it quasi-length} $\ell(X)$ of $X$ to be $n$.

\medskip

Let $Q$ be a quiver of type $\mathbb{A}_\infty^\infty$ with no infinite path. A vertex $a$ lying on a path $p$ is called a {\it middle point} of $p$ if $a$ is neither the starting point nor the ending point. A {\it string} in $Q$ is a finite reduced walk $w$, to which one associates a {\it string representation} $M(w)$; see \cite[Section 5]{BLP}. Let $a_i, b_i$, $i\in \Z$, be the source vertices and the sink vertices, respectively, in $Q$ such that there exist paths $p_i: a_i \rightsquigarrow b_i$ and $q_i: a_i\rightsquigarrow b_{i-1}$ in $Q$, for $i\in \Z$. Let $Q_R$ be the union of the paths $p_i$ with $i\in \Z$ and the trivial paths $\varepsilon_a$ with $a$ being a middle point of $q_j$ for some $j\in \Z.$ Dually, let $Q_L$ be the union of the paths $q_i$ with $i\in \Z$ and the trivial paths $\varepsilon_b$ with $b$ being a middle point of $p_j$ for some $j\in \Z$. These notations will allow us to describe the quasi-simple regular representations as follows; see \cite[(5.15)]{BLP}.

\medskip

\begin{Lemma}\label{A-di-qs}

Let $Q$ be a quiver of type $\mathbb{A}_\infty^\infty$. If $Q$ contains no infinite path, then $\Ga_{\rep(Q)}$ has exactly two regular components $\mathcal{R}_R$ and $\mathcal{R}_L$ such that the quasi-simple representations in $\mathcal{R}_R$ are the string representations $M(p)$ with $p\in Q_R,$ and those in $\mathcal{R}_L$ are the string representations $M(q)$ with $q\in Q_L$.

\end{Lemma}

\medskip

In the infinite Dynkin case, we will have an explicit description of all Auslander-Reiten components of $\rep(Q)$.

\medskip

\begin{Theo}\label{reg-idt}

Let $Q$ be an infinite Dynkin quiver. If $Q$ has no infinite path, then the connected components of $\Ga_{\rep(Q)}$ are all standard and consist of the preprojective component $\mathcal{P}$, the preinjective component $\mathcal{I}$, and $r$ regular components, where

\begin{enumerate}[$(1)$]

\item $r=0$ in case $\,Q$ is of type $\mathbb{A}_\infty;$

\vspace{1pt}

\item $r=1$ in case $\,Q$ is of type $\mathbb{D}_\infty;$

\vspace{1pt}

\item $r=2$ in case $\,Q$ is of type $\mathbb{A}_\infty^\infty;$ and in this case, the two regular components are ortho\-gonal.

\end{enumerate}

\end{Theo}

\noindent {\it Proof.} We need only to prove the second part of Statement (3), since all other parts are known; see \cite{LP} and \cite[(5.16), (5.22)]{BLP}. Assume now that $Q$ is of type $\mathbb{A}_\infty^\infty$ with no infinite path. By Lemma \ref{A-di-qs}, $\Ga_{\rep(Q)}$ has exactly two regular components $\mathcal{R}_R$ and $\mathcal{R}_L$, both are of shape $\Z \mathbb{A}_\infty$. Suppose that $\rep(Q)$ has a non-zero morphism $f: M \to N$ with $M\in \mathcal{R}_R$ and $N\in \mathcal{R}_L$. We may assume that $m=\ell(M)+\ell(N)$ is minimal with respect to this property. Suppose that $\ell(N)>1.$ Then $\rep(Q)$ has a short exact sequence
$$\xymatrix{0\ar[r] & X\ar[r]^u & N \ar[r]^v & Y\ar[r] &0,}$$ where $X, Y\in \mathcal{R}_L$ with $\ell(X)=\ell(N)-1$ and $\ell(Y)=1$. By the minimality of $m$, we have $vf=0$, and hence, $f=uw$ for some non-zero morphism $w: M\to X$, which contradicts the minimality of $m$. Therefore, $\ell(N)=1$, and dually, $\ell(M)=1$. That is, $M, N$ are quasi-simple such that their supports intersect.

Let $a_i, b_i$, $i\in \Z$, be the source vertices and the sink vertices, respectively such that $Q$ has paths $p_i: a_i \rightsquigarrow b_i$ and $q_i: a_i\rightsquigarrow b_{i-1}$ in $Q$, for $i\in \Z$. By Lemma \ref{A-di-qs}, $M=M(p_r)$ with $r\in \Z$ or $M=M(\varepsilon_a)$, where $a$ is a middle point of some $q_s$ with $s\in \Z$. In the first case, $N=M(q_i)$ with $r\le i\le r+1$ or $N=M(\varepsilon_b)$ with $b$ a middle point of $p_r$. Since $\Hom_{\,\rep(Q)}(M, N)\ne 0$, the top of $M$, that is the simple representation $S_{a_r}$, appears as a composition factor of $N$. Therefore, $N=M(q_r)$. Hence, the socle of $N$, that is, the simple representation $S_{b_{r-1}}$, is a composition factor of $M$, which is absurd. In the second case, $M=S_a$ and $N=M(q_s)$. Since $a$ is a middle point of $q_s$, we see that $\Hom_{\hspace{0.4pt}\rep(Q)}(S_a, M(q_s))=0$, a contradiction. This shows that $\Hom_{\hspace{0.4pt}\rep(Q)}(\mathcal{R}_R, \mathcal{R}_L)=0.$ Similarly, one can show that $\Hom_{\hspace{0.4pt}\rep(Q)}(\mathcal{R}_L, \mathcal{R}_R)=0.$ The proof of the theorem is completed.

\medskip

In case $Q$ is of type $\mathbb{A}_\infty$ or $\mathbb{A}_\infty^\infty$, we shall be able to obtain more information on the morphisms between the indecomposable representations.

\medskip

\begin{Lemma}\label{rep-maps}

Let $Q$ be a quiver of type $\mathbb{A}_\infty$ or $\mathbb{A}_\infty^\infty$, containing no infinite path. If $X, Y\in \Ga_{\rep(Q)}$, then $\Hom_{\,\rep(Q)}(X, Y)$ is at most one-dimensional over $k$.

\end{Lemma}

\noindent{\it Proof.} Let $X, Y\in \Ga_{\rep(Q)}$ be such that $\Hom_{\hspace{0.4pt}\rep(Q)}(X, Y)\ne 0$. Assume first that $X$ is preprojective, that is $X=\tau^{-n}P_x$ for some $n\ge 0$ and $x\in Q_0$. If $\tau^nY=0$, then $Y=\tau^{-m}P_y$ for some $0\le m<n$. By Lemma \ref{AR-formula}, $$\Hom_{\hspace{0.4pt}\rep(Q)}(\tau^{m-n}P_x, P_y)\cong \Hom_{\hspace{0.4pt}\rep(Q)}(X, Y) \ne 0.$$
Then $\tau^{m-n}P_x$ is projective; see \cite[(4.3)]{BLP}, which is absurd. Now, applying Lemma \ref{AR-formula} again, we obtain
$$\Hom_{\hspace{0.4pt}\rep(Q)}(X, Y) \cong  \Hom_{\hspace{0.4pt}\rep(Q)}(P_x, \tau^nY),$$
which is one-dimensional over $k$ because $\tau^nY$ is a string representation; see \cite[(5.9)]{BLP}. Dually, the result holds if $Y$ is
preinjective.

Assume now that $X$ is regular. By Theorem \ref{structureARquiverGen}, $Y$ is regular or preinjective. We only need to consider the case where $Y$ is regular. By Theorem \ref{reg-idt}, $X, Y$ belong to a connected component of $\Ga_{\rep(Q)},$ which is standard of shape $\Z\mathbb{A}_\infty$. Thus, the result follows from Proposition \ref{rectangle}. Finally, if $X$ is preinjective, then so is $Y$ by Theorem \ref{structureARquiverGen}, and hence, the statement holds. The proof of the lemma is completed.

\medskip

Let $\Ga$ be a connected component of $\Ga_{\rep(Q)}$ of shape $\mathbb{Z}\mathbb{A}_\infty$, containing a quasi-simple representation $S$. Observe that $\Ga$ has a unique ray starting in $S$, written as $(S\hspace{-3pt}\to)$, and a unique co-ray ending in $S$ written as $(\to \hspace{-2.5pt}S)$. We denote by $\mathcal{W}(S)$ the full subquiver of $\Ga$ generated by the representations $X$ for which there exist paths $M\rightsquigarrow X \rightsquigarrow N$, where $M\in (\to \hspace{-3pt}S)$ and $N\in (S\hspace{-2.5pt}\to)$, and call it the {\it infinite wing} with {\it wing vertex} $S$; compare \cite{Rin}.

\medskip

\begin{Prop} \label{wing}

Let $Q$ be a quiver of type $\mathbb{A}_\infty^\infty$, containing no infinite path. If $X\in \Ga_{\rep(Q)}\vspace{1pt}$ is preprojective, then each regular component $\mathcal{R}$ of $\Ga_{\rep(Q)}$ has a unique quasi-simple $S_X$ such that, for any $Y\in \mathcal{R}$, we have $\Hom_{\,\rep(Q)}(X, Y)\ne 0\vspace{1pt}$ if and only if $\hspace{0.4pt}Y\in \mathcal{W}(S_X);$ and any non-zero morphism $f: X\to Y$ factors through a representation lying on the co-ray $(\to \hspace{-2.5pt}S_X)$ ending in $S_X$.

\end{Prop}

\noindent{\it Proof.} Let $a_i, b_i$, $i\in \Z$, be the source vertices and the sink vertices, respectively such that $Q$ has paths $p_i: a_i \rightsquigarrow b_i$ and $q_i: a_i\rightsquigarrow b_{i-1},$ $i\in \Z$. By Lemma \ref{A-di-qs}, $\Ga_{\rep(Q)}$ has exactly two regular components $\mathcal{R}_R$ and $\mathcal{R}_L$. We shall consider only the case where $\mathcal{R}=\mathcal{R}_R$. Then the quasi-simple representations in $\mathcal{R}$ are the string representations $M(p)$ with $p\in Q_R$, where $Q_R$ denotes the union of the $p_i$ with $i\in \Z$, and the trivial paths $\varepsilon_a$, where $a$ ranges over the set of middle points of the $q_j$ with $j\in \Z.$

\vspace{1pt}

Let $X$ be a preprojective representation in $\Ga_{\rep(Q)}$. Since $\mathcal{R}$ is $\tau_{\hspace{-1.5pt}_Q}$-stable, in view of Lemma \ref{AR-formula}, we may assume that $X = P_x$, for some $x \in Q_0$. Since $x$ appears in exactly one of the paths in $Q_R$, we have a unique quasi-simple representation $S_X\in \mathcal{R}$ such that $\Hom_{\hspace{0.4pt}\rep(Q)}(P_x, S_X)\ne 0$. Observe that each almost split sequence
$\xymatrixcolsep{20pt}\xymatrix{0\ar[r]& U\ar[r] &V\ar[r] &W\ar[r]&0}$ with $W\in \mathcal{R}$ yields an exact sequence
\vspace{2pt}$$(*)\qquad \xymatrixcolsep{20pt}\xymatrix{0\ar[r]& \Hom(P_x, U)\ar[r] & \Hom(P_x, V)\ar[r] & \Hom(P_x, W)\ar[r]&0.}$$

Let $Y \in\mathcal{R}$ be of quasi-length $n$. By Lemma \ref{rep-maps}, $\Hom_{\hspace{0.4pt}\rep(Q)}(P_x, Y)$ is at most one-dimensional over $k$. We claim that ${\rm Hom}_{\rep(Q)}(P_x, Y)\ne 0\vspace{1pt}$ if and only if $Y\in \mathcal{W}(S_X)$.

If $n=1$, then $Y\in \mathcal{W}(S_X)$ if and only if $Y=S_X$, and hence, the claim holds. Assume that $n>1$.
Consider first the case where $Y\not\in \mathcal{W}(S_X)$. Then $\rep(Q)$ has an almost split sequence
$$\xymatrixcolsep{20pt}\xymatrix{0\ar[r]& U\ar[r] &V\ar[r] &W\ar[r]&0,}$$ where $U, W$ belong to $\mathcal{R} \backslash \mathcal{W}(S_X)$ and are of quasi-length $n-1$, and $Y$ is a direct summand of $V$. By the induction hypothesis, $\Hom_{\hspace{0.4pt}\rep(Q)}(P_x, U)=0\vspace{1pt}$
and $\Hom_{\hspace{0.4pt}\rep(Q)}(P_x, W)=0.$ In view of the sequence $(*)$, we have $\Hom_{\hspace{0.4pt}\rep(Q)}(P_x, Y)=0$.

\vspace{1pt}

Consider now the case where $Y\in \mathcal{W}(S_X)$. Then $Y=\tau^iM$ for some $0\le i\le n-1$ and $M\in (S_X\hspace{-3pt}\to).$ We denote by $N$ the immediate predecessor of $M$ in $(S_X\hspace{-3pt}\to)$.
Suppose first that $i=0$, that is, $Y=M$. Since there exists a monomorphism from $S_X$ to $Y$, we see that $\Hom_{\hspace{0.4pt}\rep(Q)}(P_x, Y)\ne 0$.
Suppose that $0<i\le n-1$. Then, there exists an almost split sequence
$$\xymatrixcolsep{25pt}\xymatrix{0\ar[r]&\tau^iN\ar[r]^-{u_1\choose u_2}& Z\oplus Y \ar[r]^-{(v_1, v_2)} & \tau^{i-1}N \ar[r] &0,}$$
where $Z, \tau^iN, \tau^{i-1}N\in \mathcal{W}(S_X)$ are of quasi-length $n-1$ or $n-2$. By the induction hypothesis on $n$, $\Hom(P_x, \tau^iN)$, $\Hom(P_x, Z)$ and $\Hom(P_x, \tau^{i-1}N)$ are one-dimensional over $k$, and by the exactness of the sequence $(*)$, so is $\Hom(P_x, Y)$. This establishes our claim.

Finally, let $f: P_x\to Y$ be a non-zero morphism. By our claim, $Y\in \mathcal{W}(S_X)$. Then there exists a unique representation $Y'\in (\to \hspace{-3pt}S_X)$ which is a sectional predecessor of $Y$ in $\mathcal{R}$. Observe, moreover, that there exists  a monomorphism $g: Y'\to Y$. This yields a monomorphism $\Hom(P_x, g): \Hom_{\,\rep(Q)}(P_x, Y')\to  \Hom_{\,\rep(Q)}(P_x, Y)$. Since $\Hom_{\,\rep(Q)}(P_x, Y')$ and $\Hom_{\,\rep(Q)}(P_x, Y)$ are one-dimensional, $\Hom(P_x, g)$ is an isomorphism. In particular, $f$ factors through $g$. The proof of the proposition is completed.

\section{Derived categories }

Throughout this section, $Q$ stands for a locally finite quiver with no infinite path. We shall study the bounded derived category $D^b(\rep(Q))$ of $\rep(Q)$. It is well known that $D^b(\rep(Q))$ is a Hom-finite Krull-Schmidt triangulated $k$-category ha\-ving almost split triangles; see \cite[(7.11)]{BLP}, whose Auslander-Reiten translation will be written as $\tau_{\hspace{-1pt}_D}.\vspace{0.5pt}$ In particular, $D^b(\rep(Q))$ admits a Serre functor $\mathbb{S}$ such that $\mathbb{S}\circ [-1]=\tau_{\hspace{-1pt}_D}$; see \cite{RVDB}. This fact is expressed in the following result.

\medskip

\begin{Lemma}\label{Serre}

Let $Q$ be a locally finite quiver with no infinite path. If $X, Y$ are indecomposable objects in $D^b(\rep(Q))$, then
$$\Hom_{D^b(\rep(Q))}(Y, \tau_{_D}X)\cong D\Hom_{D^b(\rep(Q)}(X, Y[1]).$$

\end{Lemma}

\medskip

As usual, we shall regard $\rep(Q)$ as a full subcategory of $D^b(\rep(Q))$ by identifying a representation $X$ with the stalk complex $X[0]$ concentrated in degree $0$. The Auslander-Reiten quiver $\Ga_{D^b(\rep(Q))}\vspace{1pt}$ of $D^b(\rep(Q))$ has a \emph{connecting component} $\mathcal{C}$, which is obtained by gluing the preprojective component $\mathcal{P}$ of $\Ga_{\rep(Q)}$ with the shift by $-1$ of the preinjective component $\mathcal I$ in such a way that each arrow $x \to y$ in $Q$ induces a unique arrow $I_x[-1] \to P_y$
in $\mathcal C$; see \cite[Section 7]{BLP}, and compare \cite{Hap}. For convenience, we quote the following result from  \cite[Section 7]{BLP}.

\medskip

\begin{Prop} \label{Der-cpt}

Let $Q$ be an infinite connected quiver, which is locally finite and contains no infinite path.

\vspace{-2pt}

\begin{enumerate}[$(1)$]

\item The connecting component $\mathcal C$ of $\Ga_{D^b(\rep(Q))}\vspace{1pt}$ is standard of shape $\mathbb{Z}Q^{\rm \, op}.$

\item The connected components of $\Ga_{D^b(\rep(Q))}\vspace{1pt}$ are the shifts of $\mathcal C$ and the shifts of the regular components of $\Ga_{\rep(Q)}.$ 

\end{enumerate}

\end{Prop}

\medskip

Specializing to the infinite Dynkin case, we will have a better description of the morphisms between indecomposable objects of $D^b(\rep(Q)).$

\medskip

\begin{Theo} \label{cpt-der-idt}

Let $Q$ be an infinite Dynkin quiver with no infinite path.

\vspace{-2pt}

\begin{enumerate}[$(1)$]

\item All the connected components of $\Ga_{D^b(\rep(Q))}\vspace{1pt}$ are standard.

\vspace{1pt}

\item If $\mathcal{R}, \mathcal{S}$ are distinct regular components of $\Ga_{\rep(Q)}$, then $\mathcal{R}[i], \mathcal{S}[j]$ with $i, j\in \Z$ are orthogonal in $D^b(\rep(Q))$.

\item If $M, N$ are objects lying in different connected components of $\Ga_{D^b(\rep(Q))}$, then $\Hom_{D^b(\rep(Q))}(M, N)=0$ or $\Hom_{D^b(\rep(Q))}(N, M)=0.$

\vspace{1pt}

\item There exists no short cycle in $D^b(\rep(Q)).$
\end{enumerate}
\end{Theo}

\noindent {\it Proof.} Statement (1) follows from Propositions \ref{reg-idt} and \ref{Der-cpt}. For Statement (2), let $\mathcal{R}, \mathcal{S}$ be two distinct regular components of $\Ga_{\rep(Q)}\vspace{1pt}$ with $M\in \mathcal{R}$ and $N\in \mathcal{S}$. Since $\rep(Q)$ is hereditary, $\Hom_{D^b(\rep(Q))}(M, N[j])=0$ for $j\ne 0, 1.\vspace{2pt}$ By Theorem\ref{reg-idt}, $\Hom_{D^b(\rep(Q))}(M, N)=0$, and by Lemma \ref{AR-formula},
$$\Hom_{D^b(\rep(Q))}(M, N[1])\cong \Ext_{\rep(Q)}^1(M, N)\cong D\Hom_{\hspace{0.4pt}\rep(Q)}(N, \tau_{\hspace{-1pt}_Q} M)=0.$$ This shows that $\Hom_{D^b(\rep(Q))}(\mathcal{R}, \mathcal{S}[j])=0\vspace{2pt}$, for all $j\in \Z$. By symmetry, we have $\Hom_{D^b(\rep(Q))}(\mathcal{S}, \mathcal{R}[j])=0,\vspace{1pt}$ for all $j\in \Z$. As a consequence, $\mathcal{R}[i]$ and $\mathcal{S}[j]$ are orthogonal for all $i, j\in \Z$.

\vspace{1pt}

Suppose now that there exist distinct components $\Ga, \Oa$ of $\Ga_{D^b(\rep(Q))}$ such that $\Hom_{D^b(\rep(Q))}(\Ga, \Oa)\ne 0$ and $\Hom_{D^b(\rep(Q))}(\Oa, \Ga)\ne 0.\vspace{1pt}$ Since $\rep(Q)$ is hereditary, making use of Statement (2), we may assume that $\Ga$ is the connecting component $\mathcal{C}$ of $\Ga_{D^b(\rep(Q))}.\vspace{1pt}$
Suppose that $\Oa=\mathcal{C}[i]$ for some $i\ne 0$. Since $\rep(Q)$ is hereditary, $i\in \{-1, 1, 2\}$, and since $\Hom_{\hspace{0.4pt}\rep(Q)}(\mathcal{I}, \mathcal{P})=0,$ we have $i\in \{1, 2\}$. This yields $$\Hom_{D^b(\rep(Q))}(\mathcal{C}, \mathcal{C}[-i])\cong \Hom_{D^b(\rep(Q))}(\mathcal{C}[i], \mathcal{C})=\Hom_{D^b(\rep(Q))}(\Oa, \Ga)\ne 0$$ with $-i\in \{-1, -2\}$, a contradiction. Therefore, $\Oa=\mathcal{R}[j]$ with $\mathcal{R}$ some regular component of $\Ga_{\rep(Q)}$ and $j\in \Z.$ Since $\rep(Q)$ is hereditary, $j\in \{-1, 0, 1\}$. Since $\Hom_{\hspace{0.4pt}\rep(Q)}(\mathcal{I}, \mathcal{R})=0$, we have $j\in \{0, 1\}$;
since $\Hom_{\hspace{0.4pt}\rep(Q)}(\mathcal{R}, \mathcal{P})=0$, we conclude that $j=1$. Now, if $X\in \mathcal{P}$ and $Y\in \mathcal{R}$, we have
$$\Hom_{D^b(\rep(Q))}(X, Y[1])\cong \Ext^1_{\rep(Q)}(X, Y)\cong D\Hom_{\hspace{0.4pt}\rep(Q)}(Y, \tau_{\hspace{-1.5pt}_Q}X)=0.$$
As a consequence, $\Hom_{D^b(\rep(Q))}(\Ga, \Oa)=\Hom_{D^b(\rep(Q))}(\mathcal{C}, \mathcal{R}[1])=0,\vspace{1pt}$ a contradiction. Statement (3) is established. \vspace{-4pt}

Finally, suppose that $D^b(\rep(Q))$ admits a short cycle $\xymatrixcolsep{20pt}\xymatrix{X\ar[r]^f & Y\ar[r]^g & X,}$ where $X, Y \in \Ga_{D^b(\rep(Q))}.\vspace{1pt}$ By Statement (2), $X, Y$ lie in a connected component of $\Ga_{D^b(\rep(Q))}.\vspace{1pt}$ This is absurd, since all the connected components of $\Ga_{D^b(\rep(Q))}$ are standard without oriented cycles. The proof of the theorem is completed.

\section{Cluster Categories}

Throughout this section, let $Q$ be a locally finite quiver with no infinite path. In particular, the bounded derived category $D^b(\rep(Q))$ has almost split triangles, whose Auslander-Reiten translation is written as $\tau_{\hspace{-1pt}_D}$. Setting $F=\tau_{\hspace{-1.5pt}_D}^{-1}\circ [1]$, one obtains
the canonical orbit category $$\mathscr{C}(Q)=D^b(\rep(Q))/F.$$

Recall that the objects of $\mathscr{C}(Q)$ are those of $D^b(\rep(Q));$ and for any objects $X, Y$, the morphisms are given by
$$\Hom_{\mathscr{C}(Q)}(X, Y)={\oplus}_{i\in \Z} \Hom_{D^b(\rep(Q))}(X, F^iY).$$ The composition of morphisms is given by
$$(g_i)_{i\in \mathbb{Z}} \circ (f_i)_{i\in \mathbb{Z}}=(h_i)_{i\in \mathbb{Z}},$$
where $h_i=\sum_{p+q=i} F^p(g_q)\circ f_p$. There exists a canonical projection functor $$\pi: D^b(\rep(Q))\to \mathscr{C}(Q),$$ which acts identically on the objects, and sends a morphism $f: X\to Y$ to $$(f_i)_{i\in \mathbb{Z}}\in \oplus_{i\in \mathbb{Z}} \Hom_{D^b({\rm rep}(Q))}(X, F^iY)
=\Hom_{\mathscr{C}(Q)}(X, Y),$$ where $f_i=f$ if $i=0;$ and otherwise, $f_i=0$. In particular, $\pi$ is faithful.

\medskip

\begin{Theo}\label{ccpt}

Let $Q$ be an infinite connected quiver, which is locally finite and contains no infinite path.

\vspace{-2pt}

\begin{enumerate}[$(1)$]

\item The canonical orbit category $\mathscr{C}(Q)$ is a Hom-finite Krull-Schmidt 2-Calabi-Yau triangulated $k$-category.

\item The canonical projection $\pi: D^b(\rep(Q))\to \mathscr{C}(Q)$ is a faithful triangle-functor, sending Auslander-Reiten triangles to Auslander-Reiten triangles.

\item If $\Ga$ is a connected component of $\Ga_{D^b(\rep(Q))}$, then $\pi(\Ga)$ is a connected component of $\Ga_{\hspace{-1pt}_{\mathscr{C}(Q)}}$
such that $\pi(\Ga)\cong \Ga$ as translation quivers.

\item The connected components of $\Ga_{\hspace{-1pt}_{\mathscr{C}(Q)}}$ are $\pi(\Ga)$, where $\Ga$ is either the connecting component of $\Ga_{D^b(\rep(Q))}$ or a regular component of $\Ga_{\rep(Q)}$.

\end{enumerate} \end{Theo}

\noindent{\it Proof.} Firstly, by a well known result of Keller stated in \cite[Section 9]{Kel}, $\mathscr{C}(Q)$ is a 2-Calabi-Yau triangulated $k$-category such that the canonical projection functor $\pi: D^b(\rep(Q))\to \mathscr{C}(Q)$ is exact and faithful. Moreover, since $\rep(Q)$ is Hom-finite and hereditary, $\mathscr{C}(Q)$ is Hom-finite.

We denote by $\mathscr{D}$ a skeleton of $D^b(\rep(Q))$, containing the objects in $\Ga_{D^b(\rep(Q))}.\vspace{1pt}$ Then $\mathscr{D}$ is a triangulated category such that the inclusion functor $\mathscr{D}\to D^b(\rep(Q))$ is a triangle-equivalence and $\Ga_\mathscr{D}=\Ga_{D^b(\rep(Q))}$. Observe that the Auslander-Reiten translation $\tau_{\hspace{-1.5pt}_D}$ for $D^b(\rep(Q))$ induces an automorphism of $\mathscr{D}$, which is denoted again by $\tau_{\hspace{-1.5pt}_D}$. Setting $F=\tau_{\hspace{-1.5pt}_D}^{-1}\circ [1]$, we obtain a group $G=\{F^n\mid n\in \Z\}$ of automorphisms of $\mathscr{D}$, whose action on $\mathscr{D}$ is free and locally bounded. Now, the image $\mathscr{C}$ of $\mathscr{D}$ under the canonical projection $\pi: D^b(\rep(Q))\to \mathscr{C}(Q)$ is a dense full subcategory of $\mathscr{C}(Q)$. Restricting $\pi: D^b(\rep(Q))\to \mathscr{C}(Q)$, we obtain a triangle functor $\mathscr{D}\to \mathscr{C}$ which, by abuse of notation, is denoted again by $\pi$.
For $X\in \mathscr{D}$ and $n\in \Z$, we define
$$\delta_{n, X}=(\delta_{n,i})_{i\in \mathbb{Z}}\in{\oplus}_{i\in \Z} \,\Hom_\mathscr{D}(F^nX, F^iX)=\Hom_\mathscr{C}(F^nX, X),$$ where $\delta_{n,i}=\id_{F^nX}$ if $i=n$; otherwise, $\delta_{n,i}=0$. It is easy to see that $\delta_{n,X}$ is an isomorphism which is natural in $X$ such that
$$\delta_{n, X} \circ\delta_{m, F^nX}=\delta_{n+m,X}, \mbox{ for } m, n\in \Z.$$
This yields functorial isomorphisms $\delta_n: \pi\circ F^n\to \pi$, $n\in \Z$, such that $\delta=(\delta_n)_{n\in \mathbb{Z}}$ is a $G$-stabilizer for $\pi$. It not hard to verify that the map
$$\pi_{_{X, Y}}: \oplus_{i\in \mathbb{Z}}\,\Hom_{\mathscr{D}}(X, F^iY)\to \Hom_{\,\mathscr{C}}(X, Y): (f_i)_{i\in \mathbb{Z}}\mapsto
{\textstyle\sum}_{i\in \mathbb{Z}}\, \delta_{i, Y}\circ \pi(f_i)$$ is the identity map. Hence, $\pi$ is a $G$-precovering. Since $\rep(Q)$ is Hom-finite and abelian, it is well known that $D^b({\rm rep}(Q))$ is Hom-finie and Krull-Schmidt. In particular, the endomorphism algebra of an indecomposable object in $D^b({\rm rep}(Q))$ is local with a nilpotent radical. By Lemma 2.9 stated in \cite{BL}, the functor $\pi$ satisfies Conditions (2) and (3) stated in \cite[(2.8)]{BL}. In particular, $\mathscr{C}$ is Krull-Schmidt. Moreover, since $\pi$ is dense, it is a Galois $G$-covering; see \cite[(2.8)]{BL}.

In view of Proposition 3.5 stated in \cite{BL}, we see that the exact functor $\pi: \mathscr{D}\to \mathscr{C}$ sends Auslander-Reiten triangles to Auslander-Reiten triangles, and hence Statement (2) holds. For proving Statement (3), observe that $\Ga_{\hspace{-1pt}_{\mathscr{C}(Q)}}=\Ga_\mathscr{C}$. Let $\Ga$ be connected component of $\Ga_{\mathscr{D}}$. By Theorem 4.7 stated in \cite{BL}, $\pi(\Ga)$ is a connected component of $\Ga_\mathscr{C}$ such that $\pi$ restricts to Galois $G_{\it\Gamma}$-covering $\pi_{\it\Gamma}: \Ga\to \pi(\Ga)$, where $G_{\it\Gamma}=\{F^n \mid F^n(\Ga)=\Ga\}$. Since $Q$ is infinite, $F^n(\Ga)\ne \Ga$ for each $n\ne 0.$ That is, $G_{\it\Gamma}$ is trivial, and hence, $\pi_{\it\Gamma}$ is an isomorphism of translation quivers.

Finally, since $\pi$ is dense, $\Ga_\mathscr{C}$ consists of the connected components $\pi(\it\Theta)$ with $\it\Theta$ ranging over the connected components of $\Ga_\mathscr{D}$. Now, if $\it\Theta$ is such a component, then $\it\Theta=F^n(\Ga)$, where $n\in \Z$ and $\Ga$ is the connecting component of $\Ga_{\mathscr{D}}$ or a connected component of $\Ga_{\rep(Q)}$. This yields $\pi(\it\Theta)=\pi(\Ga)$. The proof of the theorem is completed.

\medskip

\noindent{\sc Remark.} (1) As seen in the proof of Theorem \ref{ccpt}, the objects of $D^b({\rm rep}(Q))$ lying in the connecting component of $\Ga_{D^b(\rep(Q))}\vspace{1pt}$ or a regular components of $\Ga_{\rep(Q)}$ form a {\it fundamental domain} for $\mathscr{C}(Q)$, denoted by $\mathscr{F}(Q)$, that is every indecomposable object in $\mathscr{C}(Q)$ is isomorphic to a unique object in $\mathscr{F}(Q)$.

(2) By abuse of language and notation, we shall identify the connecting component $\C$ of $\Ga_{D^b(\rep(Q))}$ with $\pi(\C)$ and call it the {\it connecting component} of $\Ga_{\hspace{-1pt}_{\mathscr{C}(Q)}}$, and identify a regular component $\mathcal{R}$ of $\Ga_{\rep(Q)}$ with $\pi(\mathcal{R})$ and call it a {\it regular component} of $\Ga_{\hspace{-1pt}_{\mathscr{C}(Q)}}$. Note, however, that a standard component of $\Ga_{D^b(\rep(Q))}$ is never standard as a connected component of $\Ga_{\hspace{-1pt}_{\mathscr{C}(Q)}}$.

\medskip

For our purpose, we shall need the following explicit description of the morphisms in $\mathscr{C}(Q)$ between the objects in the fundamental domain.

\medskip

\begin{Lemma} \label{maps} Let $Q$ be a locally finite quiver with no infinite path, and let $X,Y$ be representations lying in $\Ga_{\rep(Q)}$.

\begin{enumerate}[$(1)$]

\item $\Hom_{\,\mathscr{C}(Q)}(X,Y) \cong \Hom_{D^b(\rep(Q))}(X,Y) \oplus D\Hom_{D^b(\rep(Q))}(Y,\tau_{\hspace{-1.5pt}_D}^2X)$.

\vspace{2pt}

\item If $Y$ is preinjective and $X$ is not, then
$$\Hom_{\,\mathscr{C}(Q)}(X,Y[-1]) \cong \Hom_{D^b(\rep(Q))}(X,\tau_{\hspace{-1.5pt}_D}^-Y).$$

\vspace{1pt}

\item If $X$ is preinjective and $Y$ is not, then
$$\Hom_{\,\mathscr{C}(Q)}(X[-1],Y) \cong D\Hom_{D^b(\rep(Q))}(Y,\tau_{\hspace{-1.5pt}_D} X).$$

\end{enumerate}\end{Lemma}

\noindent{\it Proof.} (1) By definition, we have
$$\Hom_{\,\mathscr{C}(Q)}(X,Y)=\oplus_{i\in \mathbb{Z}}\Hom_{D^b(\rep(Q))}(X, F^iY).$$
Since $\rep(Q)$ is hereditary, $\Hom_{D^b(\rep(Q))}(X, F^iY)=0,\vspace{1pt}$ for $i\ne 0, 1$. Since
$\tau_{\hspace{-1.5pt}_D}$ is an auto-equivalence and $\mathbb{S}=\tau_{\hspace{-1.5pt}_D}^2\circ F$ is a Serre functor, we have
$$\Hom_{D^b(\rep(Q))}(X, FY)\cong \Hom_{D^b(\rep(Q))}(\tau_{\hspace{-1.5pt}_D}^2X, \mathbb{S}(Y))\cong D\Hom_{D^b(\rep(Q))}(Y, \tau_{\hspace{-1.5pt}_D}^2X).$$

(2) Firstly, we have $\Hom_{\,\mathscr{C}(Q)}(X,Y[-1]) \cong \Hom_{\,\mathscr{C}(Q)}(X, \tau_{\hspace{-1.5pt}_D}^-Y).\vspace{1pt}$ Suppose that $Y$ is preinjective and $X$ is not. We start with the case where $Y$ is not injective. Then $\tau_{\hspace{-1.5pt}_D}^-Y = \tau_{\hspace{-1.5pt}_Q}^-Y$ is a preinjective representation in $\Ga_{\,\rep(Q)}$. Applying Statement (1), we obtain
$$\Hom_{\,\mathscr{C}(Q)}(X, \tau_{\hspace{-1.5pt}_D}^-Y) \cong \Hom_{D^b(\rep(Q))}(X,\tau_{\hspace{-1.5pt}_D}^-Y) \oplus D\Hom_{D^b(\rep(Q))}(\tau_{\hspace{-1.5pt}_D}^-Y,\tau_{\hspace{-1.5pt}_D}^2X).$$
Observe that $\tau_{\hspace{-1.5pt}_D}^2X$ is a shift by $-1$ of a preinjective representation or a representation which is not preinjective. In both cases, $\Hom_{D^b(\rep(Q))}(\tau_{\hspace{-1.5pt}_D}^-Y,\tau_{\hspace{-1.5pt}_D}^2X)=0$.
Statement (2) holds in this case.

\vspace{1pt}

Assume now that $Y=I_x$ for some $x \in Q_0$. Then $\tau_{\hspace{-1.5pt}_D}^-Y = P_x[1]$. By definition, $F^iP_x[1] = \tau_{\hspace{-1.5pt}_D}^{-i}P_x[i+1]$, for all $i\in \Z$. If $i \ge 1$, then $\Hom_{D^b(\rep(Q))}(X, F^iP_x[1])=0$. If $i\le -1$, then $\tau_{\hspace{-1.5pt}_D}^{-i}P_x[i+1]=N[i]$ for some preinjective representation $N$, and hence, $\Hom_{D^b(\rep(Q))}(X, F^iP_x[1])=0$. This yields $$\Hom_{\,\mathscr{C}(Q)}(X, \tau_{\hspace{-1.5pt}_D}^-Y) = \Hom_{D^b(\rep(Q))}(X, F^0P[1]) = \Hom_{D^b(\rep(Q))}(X,\tau_{\hspace{-1.5pt}_D}^-Y).$$

\vspace{1pt}

(3) Suppose that $X$ is preinjective and $Y$ is not. Then $\tau_{\hspace{-1.5pt}_\mathscr{C}}^{-2}Y=\tau_{\hspace{-1.5pt}_D}^{-2}Y\in \Ga_{\,\rep(Q)}.$
Since $\mathscr{C}(Q)$ is 2-Calabi-Yau and $\tau_{\hspace{-1.5pt}_\mathscr{C}}=[1]$, we have
$$\begin{array}{rcl}
\Hom_{\,\mathscr{C}(Q)}(X[-1],Y) & \cong & D\Hom_{\,\mathscr{C}(Q)}(Y,X[1])\\ \vspace{-10pt}\\
&\cong & D\Hom_{\,\mathscr{C}(Q)}(\tau_{\hspace{-1.5pt}_\mathscr{C}}^{-2}Y, X[-1])\\ \vspace{-10pt}\\
&\cong & D\Hom_{D^b(\rep(Q))}(\tau_{\hspace{-1.5pt}_D}^{-2}Y,\tau_{\hspace{-1.5pt}_D}^{-} X)\\ \vspace{-10pt}\\
&\cong & D\Hom_{D^b(\rep(Q))}(Y,\tau_{\hspace{-1.5pt}_D}X),\end{array}$$
where the third isomorphism follows from Statement (2). The proof of the lemma is completed.

\medskip

\begin{Cor} \label{special-maps} Let $Q$ be a locally finite quiver with no infinite path, and let $X, Y$ be representations in $\Ga_{\rep(Q)}.\vspace{1pt}$ If $X$ is preprojective and $Y$ is regular, then $\Hom_{\mathscr{C}(Q)}(X, Y)\cong \Hom_{\,\rep(Q)}(X, Y)$.

\end{Cor}

\noindent{\it Proof.} Suppose that $X$ is preprojective and $Y$ is regular. Then $\tau_{\hspace{-1.5pt}_D}^2X$ is either a preprojective representation or the shift by $-1$ of a preinjective representation. In both cases, we have
$\Hom_{D^b(\rep(Q))}(Y, \tau_{\hspace{-1.5pt}_D}^2X)=0.\vspace{1pt}$ Thus, by Lemma \ref{maps}(1), $\Hom_{\mathscr{C}(Q)}(X, Y)\cong \Hom_{\,\rep(Q)}(X, Y)$. The proof of the corollary is completed.

\medskip


The following result says that $\mathscr{C}(Q)$ has at least one cluster-tilting subcategory.

\medskip

\begin{Lemma} \label{Proj-ctsc}

Let $Q$ be a locally finite quiver with no infinite path. The strictly additive subcategory $\mathscr{P}$ of $\,\mathscr{C}(Q)$ generated by the representations $P_x$ with $x\in Q_0$, is a cluster-tilting subcategory.

\end{Lemma}

\noindent{\it Proof.} Observe that the Auslander-Reiten translation $\tau_{\hspace{-1.5pt}_\mathscr{C}}$ for $\mathscr{C}(Q)$ coincides with the shift functor. For any $x, y\in Q_0$, making use of Lemma \ref{maps}(2), we obtain
$$\begin{array}{rcl}
\Hom_{\mathscr{C}(Q)}(P_x, P_y[1])
&=&\Hom_{\mathscr{C}(Q)}(P_x, \tau_{\hspace{-1.5pt}_\mathscr{C}}P_y) = \Hom_{\mathscr{C}(Q)}(P_x, I_y[-1])\\
\vspace{-10pt}\\
&\cong &\Hom_{D^b(\rep(Q))}(P_x, \tau_{\hspace{-1.5pt}_D}^-I_y)= \Hom_{D^b(\rep(Q))}(P_x, P_y[1]) \\ \vspace{-10pt}\\
&=& {\rm Ext}_{\,\rep(Q)}^1(P_x, P_y)=0.
\end{array}$$

Let $X$ be an indecomposable object in the fundamental domain of $\mathscr{C}(Q)$, but not in $\mathscr{P}$.
Assume first that $X\in \Ga_{{\rm rep}(Q)}.\vspace{1pt}$ Then, $\tau_{\hspace{-1.5pt}_{\mathscr{C}}}X=\tau_{\hspace{-1.5pt}_D}X=\tau_{\hspace{-1.5pt}_Q}X\in \Ga_{{\rm rep}(Q)}$. Choosing a vertex $x$ in the support of $\tau_{\hspace{-1.5pt}_Q}X$, in view of Lemma \ref{maps}(1), we obtain
$$\begin{array}{rcl}
\Hom_{\mathscr{C}(Q)}(P_x, X[1])&=&\Hom_{\mathscr{C}(Q)}(P_x, \tau_{_{\mathscr{C}}}X)
=\Hom_{\mathscr{C}(Q)}(P_x, \tau_{\hspace{-1.5pt}_Q}X)\\ \vspace{-10pt}\\
&\cong &\Hom_{\,\rep(Q)}(P_x, \tau_{\hspace{-1.5pt}_Q}X) \oplus D\Hom_{D^b(\rep(Q))}(\tau_{\hspace{-1.5pt}_Q}X, \tau_{\hspace{-1.5pt}_D}^2P_x)\\
&\ne & 0.
\end{array}$$

Assume now that $X=Y[-1]$, where $Y$ is a preinjective representation in $\Ga_{{\rm rep}(Q)}$. If $y$ is a vertex in the support of $Y$, then
$$\begin{array}{rcl}
\Hom_{\mathscr{C}(Q)}(P_y, X[1]) &= & \Hom_{\mathscr{C}(Q)}(P_y, Y) \\
&\cong & \Hom_{D^b(\rep(Q))}(P_y, Y) \oplus D\Hom_{D^b(\rep(Q))}(Y,\tau_{\hspace{-1.5pt}_D}^2P_y)\\
&\ne& 0. \end{array}$$

This shows that $\mathscr{P}$ is weakly cluster-tilting. Next, we need to show that $\mathscr{P}$ is functorially finite in $\mathscr{C}(Q)$. For this end, let $Z$ be an indecomposable object in the fundamental domain of $\mathscr{C}(Q)$. We claim that both $\Hom_{\mathscr{C}(Q)}(Z, -)$ and $\Hom_{\mathscr{C}(Q)}(-, Z)$ vanish on $\mathscr{P}$ for all but finitely many objects. Start with the case where $Z\in \Ga_{\rep(Q)}.$
Let $x\in Q_0$ be such that the support of $P_x$ does not intersect the support of $Z\oplus \tau^2_{\hspace{-1.5pt}_Q}Z$.
By Lemma \ref{maps}(1), we have
$$\begin{array}{rcl}
\Hom_{\mathscr{C}(Q)}(P_x, Z)
&\cong& \Hom_{D^b({\rm rep}(Q))}(P_x, Z)\oplus  D\Hom_{D^b({\rm rep}(Q))}(Z, \tau_{\hspace{-1.5pt}_D}^2P_x) \\ \vspace{-10pt}\\
&\cong& \Hom_{{\rm rep}(Q)}(P_x, Z)\oplus  D\Hom_{D^b({\rm rep}(Q))}(Z, (\tau_{\hspace{-1.5pt}_Q}I_x)[-1])\\ \vspace{-10pt}\\
&=&0;
\end{array}$$
and
$$\begin{array}{rcl}
\Hom_{\mathscr{C}(Q)}(Z, P_x)
&\cong& \Hom_{D^b({\rm rep}(Q))}(Z, P_x)\oplus  D\Hom_{D^b({\rm rep}(Q))}(P_x, \tau_{\hspace{-1.5pt}_D}^2Z)\\ \vspace{-10pt}\\
&\cong& \Hom_{\,\rep(Q)}(Z, P_x)\oplus  D\Hom_{D^b({\rm rep}(Q))}(P_x, \tau_{\hspace{-1.5pt}_D}^2Z)\\ \vspace{-10pt}\\
&=& D\Hom_{D^b({\rm rep}(Q))}(P_x, \tau_{\hspace{-1.5pt}_D}^2Z), \\ \vspace{-10pt}\\
\end{array}$$ where the last equation follows from the hypothesis on $x$. If $Z=\tau_{\hspace{-1.5pt}_Q}^{-i}P_y$ for some $y\in Q_0$ and $0\le i\le 1$, then
$$\Hom_{D^b({\rm rep}(Q))}(P_x, \tau_{\hspace{-1.5pt}_D}^2Z)=\Hom_{D^b({\rm rep}(Q))}(P_x, (\tau^{1-i}_{\hspace{-1.5pt}_Q}I_y)[-1])=0.$$ Otherwise, $\tau_{\hspace{-1.5pt}_D}^2Z=\tau_{\hspace{-1.5pt}_Q}^2Z\in \Ga_{\rep(Q)}$, and by the hypothesis on $x$, we obtain
$$\Hom_{D^b({\rm rep}(Q))}(P_x, \tau_{\hspace{-1.5pt}_D}^2Z)\cong \Hom_{{\rm rep}(Q)}(P_x, \tau_{\hspace{-1.5pt}_Q}^2Z)=0.$$
Since $Q$ has no infinite path, our claim holds in this case.

Next, suppose that $Z=N[-1]$, for some preinjective $N\in \Ga_{{\rm rep}(Q)}$. Let $z\in Q_0$ be not in the support of $\tau_{\hspace{-1.5pt}_Q} N \oplus \tau^-_{\hspace{-1.5pt}_Q}N$. In view of Lemma \ref{maps}(3), we obtain
$$\Hom_{\mathscr{C}(Q)}(Z, P_z)\cong D \Hom_{D^b(\rep(Q))}(P_z, \tau_{\hspace{-1.5pt}_D}N)\cong D \Hom_{\,\rep(Q)}(P_z, \tau_{\hspace{-1.5pt}_Q}N)=0; \vspace{1pt}$$
and by Lemma \ref{maps}(2), we have
$$\Hom_{\mathscr{C}(Q)}(P_z, Z)\cong \Hom_{D^b({\rm rep}(Q))}(P_z, \tau_{\hspace{-1.5pt}_D}^-N).$$
If $N=I_b$ for some $b\in Q_0$, then $$\Hom_{D^b({\rm rep}(Q))}(P_z, \tau_{\hspace{-1.5pt}_D}^-N)= \Hom_{D^b({\rm rep}(Q))}(P_z, P_b[1])=0.$$
Otherwise, by the hypothesis on $z$, we obtain
$$\Hom_{D^b({\rm rep}(Q))\vspace{2pt}}(P_z, \tau_{\hspace{-1.5pt}_D}^-N)=\Hom_{\,{\rm rep}(Q)}(P_z, \tau_{\hspace{-1.5pt}_Q}^-N)=0.$$

Since $N\oplus \tau^-N$ is finite dimensional, our claim holds in this case. As a consequence, $\mathscr{P}$ is covariantly and contrvariantly bounded in $\mathscr{C}(Q)$, and by Lemma \ref{Fonctorial-finiteness}, it is functorially finite. The proof of the lemma is completed.

\medskip

For the rest of this section, we shall concentrate on the infinite Dynkin case.

\medskip

\begin{Prop}\label{idt-clust-cpt}

Let $Q$ be an infinite Dynkin quiver with no infinite path.
The connected components of $\Ga_{_{\mathscr{C}(Q)}}$ consist of the connecting component of shape $\Z Q^{\rm op}$ and $r$ regular components of shape $\Z\mathbb{A}_\infty$, where

\begin{enumerate}[$(1)$]

\item $r=0$ if $\,Q$ is of type $\mathbb{A}_\infty;$

\vspace{1pt}

\item $r=1$ if $\,Q$ is of type $\mathbb{D}_\infty;$

\vspace{1pt}

\item $r=2$ if $\,Q$ is of type $\mathbb{A}_\infty^\infty;$ and in this case, the two regular components are ortho\-gonal.

\end{enumerate}\end{Prop}

\noindent{\it Proof.} We need only to show the second part of Statement (3), since the other parts follow from Theorem \ref{ccpt} and some results stated in \cite[(5.16),(5.17),(5.22)]{BLP}. For this purpose, suppose that $Q$ is of type $\mathbb{A}_\infty^\infty$. Let $\mathcal{R}, \mathcal{S}$ be the two distinct regular components of $\Ga_{\rep(Q)}$ with $X \in \mathcal{R}$ and $Y \in \mathcal{S}$. By Theorem \ref{cpt-der-idt}, $\mathcal{R}$ and $\mathcal{S}$ are orthogonal in $D^b(\rep(Q))$. In view of Lemma \ref{maps}(1), we obtain
$$\Hom_{\,\mathscr{C}(Q)}(X,Y) \cong \Hom_{D^b(\rep(Q))}(X,Y) \oplus D\Hom_{D^b(\rep(Q))}(Y,\tau_{\hspace{-1.5pt}_D}^2X)=0.$$
The proof of the proposition is completed.

\bigskip

Recall that an object $X\in \mathscr{C}(Q)$ is called a \emph{brick} if $\End_{\mathscr{C}(Q)}(M)$ is one-dimensional over $k;$ and {\it rigid} if $\Hom_{\mathscr{C}(Q)}(X, X[1])=0$.

\medskip

\begin{Cor}\label{rigid-obj}

Let $Q$ be an infinite Dynkin quiver. If $Q$ has no infinite path, then every indecomposable object in $\mathscr{C}(Q)$ is a rigid brick.

\end{Cor}

\noindent{\it Proof.} Assume that $Q$ has no infinite path. Let $X$ be an indecomposable object of $\mathscr{C}(Q)$. By Theorem \ref{ccpt}(2), $\tau_{\hspace{-1.5pt}_\mathscr{C}}X=\tau_{\hspace{-1.5pt}_D}X$.
In order to show that $X$ is a rigid brick, we may assume $X$ lies in the fundamental domain $\mathscr{F}(Q)$ of $\mathscr{C}(Q)$. Since the Auslander-Reiten translation $\tau_{\hspace{-1.5pt}_\mathscr{C}}$ for $\mathscr{C}(Q)$ coincides with the shift functor, $X$ is a rigid brick if and only if so is $\tau_{\hspace{-1.5pt}_\mathscr{C}}^nX$ for some integer $n$. As a consequence, we may further assume that both $X$ and $\tau_{\hspace{-1.5pt}_D}X$ belong to $\Ga_{\rep(Q)}.\vspace{1pt}$

Let $\Ga$ be the connected component of $\Ga_{D^b(\rep(Q))}\vspace{1pt}$ containing $X$. By Theorem \ref{cpt-der-idt}(1), $\Ga$ is standard with no oriented cycle. Hence, $\Hom_{D^b(\rep(Q))}(X,\tau_{\hspace{-1.5pt}_D}X)=0\vspace{1pt}$, $\Hom_{D^b(\rep(Q))}(X,\tau_{\hspace{-1.5pt}_D}^2X)=0,\vspace{1pt}$ and $\End_{D^b(\rep(Q))}(X)$ is one-dimensional over $k$. In view of Lemma \ref{maps}(1), we deduce that $\End_{\mathscr{C}(Q)}(X)$ is one-dimensional over $k$, and
$$\begin{array}{rcl}
\Hom_{\mathscr{C}(Q)}(X, X[1]) &\cong & \Hom_{\mathscr{C}(Q)}(X, \tau_{\hspace{-1.5pt}_\mathscr{C}}X)=\Hom_{\mathscr{C}(Q)}(X, \tau_{\hspace{-1.5pt}_D}X)\\ \vspace{-10pt} \\
&\cong& \Hom_{D^b(\rep(Q))}(X, \tau_{\hspace{-1.5pt}_D}X)\oplus D\Hom_{D^b(\rep(Q))}(\tau_{\hspace{-1.5pt}_D}X, \tau_{\hspace{-1.5pt}_D}^2X)\\
&=&0.\end{array}$$
The proof of the corollary is completed.

\medskip

More generally, a strictly additive subcategory $\T$ of $\mathscr{C}(Q)$ is called {\it rigid} if $\Hom_{\mathscr{C}(Q)}(X, Y[1])=0$, for all objects $X, Y\in \T;$ and {\it maximal rigid} if it is rigid and maximal with respect to the rigidity property. By definition, a weakly cluster-tilting subcategory of $\mathscr{C}(Q)$ is maximal rigid, and the converse is not true in general.

\medskip

\begin{Lemma}\label{rigidity}

Let $Q$ be an infinite Dynkin quiver with no infinite path. If $\T$ is a strictly additive subcategory of $\mathscr{C}(Q)$, then it is weakly cluster-tilting if and only if it is maximal rigid in $\mathscr{C}(Q)$.

\end{Lemma}

\noindent{\it Proof.} We need only to show the sufficiency. Let $\T$ be a strictly additive subcategory of $\mathscr{C}(Q)$, which is maximal rigid. Let $M\in \mathscr{C}(Q)$ be indecomposable such that $\Hom_{\mathscr{C}(Q)}(\T,M[1])=0$. Since $\mathscr{C}(Q)$ is 2-Calabi-Yau, $\Hom_{\mathscr{C}(Q)}(M, \T[1])=0$. By Corollary \ref{rigid-obj}, $M$ is rigid in $\mathscr{C}(Q)$. Hence, the strictly additive subcategory of $\mathscr{C}(Q)$ generated by $M$ and $\T$ is rigid. Since $\T$ is maximal rigid, $M \in \T$. The proof of the lemma is completed.

\bigskip

The following result is essential for our investigation.

\medskip

\begin{Prop} \label{onedim}

Let $Q$ be a quiver with no infinite path of type $\mathbb{A}_\infty$ or $\mathbb{A}_\infty^\infty\vspace{1pt}$. If $X,Y \in \mathscr{C}(Q)$ are indecomposable, then $\Hom_{\mathscr{C}(Q)}(X,Y)$ is of $k$-dimension at most one.

\end{Prop}

\noindent{\it Proof.} Let $X,Y\in \mathscr{C}(Q)$ be indecomposable and suppose they lie in the fundamental domain $\mathscr{F}(Q)$.
Since $\tau_{\hspace{-1.5pt}_\mathscr{C}}$ is an auto-equivalence of $\mathscr{C}$, we may assume that $\tau_{\hspace{-1.5pt}_\mathscr{C}}^iX, \tau_{\hspace{-1.5pt}_\mathscr{C}}^iY\in \Ga_{\rep(Q)}$ for $i=0,1,2$. Then, $X, Y$ are preprojective or regular representations. Since $\Hom_{\mathscr{C}(Q)}(X,Y) \cong D\Hom_{\mathscr{C}(Q)}(Y, \tau_{\hspace{-1.5pt}_D}^2X)$, we may assume that $X$ is preprojective in case $X$ or $Y$ is not regular. By Lemma \ref{maps}(1), we obtain
$$\begin{array}{rcl}
\Hom_{\mathscr{C}(Q)}(X,Y) &= & \Hom_{D^b(\rep(Q))}(X,Y) \oplus D\Hom_{D^b(\rep(Q))}(Y, \tau_{\hspace{-1.5pt}_Q}^2X)\\ \vspace{-10pt}\\
&=& \Hom_{\hspace{0.4pt}\rep(Q)}(X,Y) \oplus D\Hom_{\hspace{0.4pt}\rep(Q)}(Y, \tau_{\hspace{-1.5pt}_Q}^2X).
\end{array}$$

Assume first that $X$ and $Y$ belong to a connected component $\Ga$ of $\Ga_{D^b(\rep(Q))}.\vspace{1pt}$ By Propositions \ref{Der-cpt} and \ref{cpt-der-idt}(1), $\Ga$ is standard of shape $\Z\mathbb{A}_\infty$ or $\Z\mathbb{A}_\infty^\infty.\vspace{1pt}$
Thus, $\Hom_{D^b(\rep(Q))}(X,Y)=0$ or $\Hom_{D^b(\rep(Q))}(Y, \tau_{\hspace{-1.5pt}_Q}^2X)=0.$ Moreover, in view of Proposition \ref{rectangle}, we see that $\Hom_{\mathscr{C}(Q)}(X,Y)$ is at most one-dimensional over $k$.

Assume now that $X, Y$ belong to two different connected components $\Ga$ of $\Ga_{D^b(\rep(Q))}$. If both $X$ and $Y$ are regular representations in $\Ga_{\rep(Q)}$ then, by Proposition \ref{reg-idt}, $\Hom_{\hspace{0.4pt}\rep(Q)}(X,Y)=0$ and
$\Hom_{\hspace{0.4pt}\rep(Q)}(Y, \tau_{\hspace{-1.5pt}_Q}^2X)=0$. Therefore, $\Hom_{\mathscr{C}(Q)}(X,Y)=0$. Otherwise, by our assumption, $X$ is preprojective and $Y$ is regular. Then $\Hom_{\hspace{0.4pt}\rep(Q)}(Y, \tau_{\hspace{-1.5pt}_Q}^2X)=0$ by Theorem \ref{structureARquiverGen}(2), and $\Hom_{\hspace{0.4pt}\rep(Q)}(X,Y)$ is at most one-dimensional over $k$ by Lemma \ref{rep-maps}. As a consequence, $\Hom_{\mathscr{C}(Q)}(X,Y)$ is at most one-dimensional over $k$. This completes the proof of the proposition.

\medskip

We are ready to obtain our main result of this section.

\medskip

\begin{Theo} \label{FZ}

Let $Q$ be a quiver of type $\mathbb{A}_\infty$ or $\mathbb{A}_\infty^\infty.$ If $Q$ contains no infinite path, then $\mathscr{C}(Q)$ is a cluster category.

\end{Theo}

\noindent{\it Proof.} Assume that $Q$ contains no infinite path. By Theorem \ref{ccpt}(1), $\mathscr{C}(Q)$ is a Hom-finite 2-Calabi-Yau triangulated $k$-category. Since $\mathscr{C}(Q)$ has a cluster-tilting subcategory by Lemma \ref{Proj-ctsc}, we need only to show that the quiver of every cluster-tilting subcategory $\T$ of $\mathscr{C}(Q)$ has no oriented cycle of length one or two; see \cite[(II.1.6)]{BIRS}. For this purpose, it suffices to show that $\T$ has no short cycle.

Suppose that $\mathscr{C}(Q)$ has a cluster-tilting subcategory with short cycles. In particular, $\mathscr{C}(Q)$ has some non-zero radical morphisms $f: X\to Y$ and $g: Y\to X$, where $X, Y\in \Ga_{\hspace{-1pt}_{\mathscr{C}(Q)}}$ form a rigid family contained in the fundamental domain $\mathscr{F}(Q)$. By Corollary \ref{rigid-obj}, $X$ and $Y$ are not isomorphic. For each integer $n$, since the Auslander-Reiten translation $\tau_{_{\mathscr{C}}}$ of $\mathscr{C}(Q)$ coincides with the shift functor, $\mathscr{C}(Q)$ has non-zero radical morphisms $f_n: \tau_{_{\mathscr{C}}}^nX\to \tau_{_{\mathscr{C}}}^nY\vspace{1pt}$ and $g_n: \tau_{_{\mathscr{C}}}^nY\to \tau_{_{\mathscr{C}}}^nX$, where $\tau_{_{\mathscr{C}}}^nX, \tau_{_{\mathscr{C}}}^nY$ form a rigid family contained in the fundamental domain of $\mathscr{C}(Q)$. Therefore, we may assume that $\tau_{\hspace{-1.5pt}_\mathscr{C}}^iX, \tau_{\hspace{-1.5pt}_\mathscr{C}}^iY\in \Ga_{\rep(Q)}$, for $0\le i\le 3$. In particular, $\tau_{\hspace{-1.5pt}_\mathscr{C}}^iY=\tau_{\hspace{-1.5pt}_D}^iY=\tau_{\hspace{-1.5pt}_Q}^iY$, for $0\le i\le 3.$

\vspace{1pt}

Suppose first that $\rep(Q)$ has a non-zero radical morphism $f_0: X\to Y$. By Theorem \ref{cpt-der-idt}(4), ${\rm rad}_{\rep(Q)}(Y, X)=0$, and thus, $\Hom_{\hspace{0.4pt}\rep(Q)}(Y, X)=0.$ In view of Lemma \ref{maps}(1), we obtain a non-zero radical morphism $g_0: X\to \tau_{\hspace{-1.5pt}_Q}^2Y$ in $\rep(Q)$. Moreover,
since $\Hom_{\mathscr{C}(Q)}(X, \tau_{\hspace{-1.5pt}_Q}Y)=\Hom_{\mathscr{C}(Q)}(X, \tau_{\hspace{-1.5pt}_\mathscr{C}}Y)=\Hom_{\mathscr{C}(Q)}(X, Y[1]) = 0,\vspace{1pt}$
we have $\Hom_{D^b(\rep(Q))}(X, \tau_{\hspace{-1.5pt}_Q}Y)=0.$ That is, $\Hom_{\hspace{0.4pt}\rep(Q)}(X, \tau_{\hspace{-1.5pt}_Q}Y)=0.$

\vspace{1pt}

Let $\Ga$ be the connected component of $\Ga_{D^b(\rep(Q))}\vspace{1pt}$ containing $X$. By Propositions \ref{Der-cpt} and \ref{cpt-der-idt}(1), $\Ga$ is standard of shape $\Z\mathbb{A}_\infty$ or $\Z\mathbb{A}_\infty^\infty.\vspace{1pt}$ If $Y\in \Ga$ then, by Proposition \ref{rectangle}, both $\tau_{\hspace{-1.5pt}_D}^2Y$ and $Y$ lie in the forward rectangle $\mathscr{R}^X$ of $X$. Being convex, $\mathscr{R}^X$ also contains $\tau_{\hspace{-1.5pt}_D}Y$. As a consequence, $\Hom_{D^b(\rep(Q))}(X, \tau_{\hspace{-1.5pt}_D}Y)\ne 0,\vspace{1pt}$ a contradiction. Therefore, $Y$ lies in a connected component $\Oa$ of $\Ga_{D^b(\rep(Q))}$ with $\Oa\ne \Ga$.

Observing that $X, Y$ are preprojective or regular representations, we deduce from Theorem \ref{structureARquiverGen} and Proposition \ref{reg-idt} that $X$ is preprojective and $Y$ is regular. Thus, $\Oa$ is a regular component of $\Ga_{\rep(Q)}$. By Proposition \ref{wing}, $\Oa$ has an infinite wing $\mathcal{W}(S)$ determined by a quasi-simple representation $S$ such that, for any $Z\in \Oa$, $\Hom_{\hspace{0.4pt}\rep(Q)}(X,Z) \ne 0$ if and only if $Z \in \mathcal{W}(S)$. In particular, $\tau_{\hspace{-1.5pt}_Q}^2Y, Y\in \mathcal{W}(S)$, and consequently, $\tau_{\hspace{-1.5pt}_Q}Y \in \mathcal{W}(S)$. That is, $\Hom_{\hspace{0.4pt}\rep(Q)}(X, \tau_{\hspace{-1.5pt}_Q}Y)\ne 0$, a contradiction.

Suppose now that $\rep(Q)$ has a non-zero radical morphism $f_1: Y\to \tau_{\hspace{-1.5pt}_Q}^2 X$. Applying Theorem \ref{cpt-der-idt}(1) and (3), we deduce that $\Hom_{D^b(\rep(Q))}(X, \tau_{\hspace{-1.5pt}_Q}^2Y)=0.$ By Lemma \ref{maps}(1), we have a non-zero radical morphism $g_1: Y\to X$. This reduces to the case we have just treated. The proof of the theorem is completed.

\section{Triangulations of the infinite strip}

From now on, we shall study the cluster category of type $\mathbb{A}_\infty^\infty$ from a geometric point of view. Our geometric model will be triangulations of the infinite strip, which was introduced by Igusa and Todorov in \cite{IgT} and further studied by Holm and J{\o}rgensen in \cite{HJo1}.

\medskip

For the rest of this paper, we shall denote by $\mathcal{B}_\infty$ the infinite strip in the plane, consisting of the points $(x, y)$ with $0\le y\le 1.$ The line defined by $y=0$ is called {\it lower boundary line} of $\B_\infty$, while the line defined by $y=1$ is called {\it upper boundary line}. The points $\ml_i=(i, 1)$, $i\in \Z$, are called {\it upper marked points}\hspace{0.3pt}; and the points $\mr_i=(-i, 0)$, $i\in \Z$, are called {\it lower marked points}.
An upper or lower marked point in $\B_\infty$ will be simply called a {\it marked point}. The set of marked points in $\B_\infty$ will be written as $\mf{M}$. Moreover, we denote by $\mf{A}$ the set of all two-element subsets of $\mf{M}$ except the subsets $\{\mr_i, \mr_j\}\vspace{1pt}$ and $\{\ml_i, \ml_j\}$ with $|i-j| \le 1$.

\medskip

By a {\it curve} in $\B_\infty$ we mean a curve in $\B_\infty$ joining two points, called {\it endpoints}. Two curves are said to {\it cross} provided that their intersection contains a point which is not an endpoint of any of the two curves. By a {\it simple curve} in $\mathcal{B}_\infty$ we mean a curve joining two marked points which does not cross itself and intersects the two boundary lines only at the endpoints. Moreover, a simple curve in $\mathcal{B}_\infty$ is called an {\it edge curve} if the set of its endpoints is either $\{\ml_i, \ml_{i+1}\}$ or $\{\mr_i, \mr_{i+1}\}$ for some $i\in \Z;$ and an {\it arc curve} if the set of its endpoints belongs to $\mf{A}$.

\medskip

Let $\mathfrak{p}, \mathfrak{q}$ be two distinct marked points in $\B_\infty$. There exists a unique isotopy class of simple curves in $\mathcal{B}_\infty$ joining $\mathfrak{p}$ and $\mathfrak{q}$, which we denote by $[\mathfrak{p}, \mathfrak{q}]$, or equivalently, by $[\mathfrak{q}, \mathfrak{p}]$. We shall call $\mathfrak{p}, \mathfrak{q}$ the {\it endpoints} of the isotopy class $[\mathfrak{p}, \mathfrak{q}]$. Now, the isotopy class of an edge curve is called an {\it edge} in $\mathcal{B}_\infty$, while the isotopy class of an arc curve is called an {\it arc} in $\mathcal{B}_\infty$. An arc in $\mathcal{B}_\infty$ is called an {\it upper arc} if its endpoints are upper marked points, and a {\it lower arc} if its endpoints are lower marked points, and a {\it connecting arc} if its endpoints do not lie on the same boundary line.

\medskip

Let $u,v$ be arcs in $\mathcal{B}_\infty$. We shall say that $u$ \emph{crosses} $v$, or equivalently, $(u, v)$ is a {\it crossing pair}, provided that every curve lying in the isotopy class $u$ crosses each of the curve lying in the isotopy class $v$. Observe that an arc does not cross itself. For convenience, we state without a proof the following easy observation.

\medskip

\begin{Lemma} \label{CombTech}

The following statement hold true for arcs in $\B_\infty$.

\vspace{-2pt}

\begin{enumerate}[$(1)$]

\item An upper arc does not cross any lower arc.

\vspace{1pt}

\item An upper arc $[\ml_i, \ml_j]$ with $i<j$ crosses a connecting arc $[\ml_r, \mr_s]$ if and only if $i<r<j$.

\item A lower arc $[\mr_i, \mr_j]$ with $i>j$ crosses a connecting arc $[\ml_r, \mr_s]$ if and only if $i>s>j$.

\item Two connecting arcs $[\ml_i, \mr_j]$ and $[\ml_p, \mr_q]$ cross if and only if
  $i>p$ and $j>q$, or else, $i<p$ and $j<q$.

\item Two upper arcs $[\ml_i, \ml_j]$ and $[\ml_p, \ml_q]$ with $i<j$ and $p<q$ cross if and only if $i<p<j<q$ or $p<i<q<j.$

\item Two lower arcs $[\mr_i, \mr_j]$ and $[\mr_p, \mr_q]$ with $i>j$ and $p>q$ cross if and only if $i>p>j>q$ or $p>i>q>j.$

\end{enumerate}\end{Lemma}

\medskip

We denote by ${\rm arc}(\mathcal{B}_\infty)$ the set of arcs in $\B_\infty$, which is equipped with a natural translation $\tau: {\rm arc}(\B_\infty)\to {\rm arc}(\B_\infty)$ as defined below.

\medskip

\begin{Defn}\label{arc-translation}

For each arc $u$ in $\B_\infty$, we define its translate $\tau u$ as follows.

\begin{enumerate}[$(1)$]

\item If $u=[\ml_i, \ml_j]$ with $i<j-1$, then $\tau u =[\ml_{i+1}, \ml_{j+1}]$.

\item If $u=[\mr_i, \mr_j]$ with $i>j+1$, then $\tau u =[\mr_{i+1}, \mr_{j+1}]$.

\item If $u=[\ml_i, \mr_j]$, then $\tau u =[\ml_{i+1}, \mr_{j+1}]$.

\end{enumerate}\end{Defn}

\medskip

The following statement follows easily from the definition of $\tau$ and Lemma \ref{CombTech}.

\medskip

\begin{Cor} \label{tau-cross}

If $u, v$ arcs in $\B_\infty$, then the following statement hold.

\vspace{-1pt}

\begin{enumerate}[$(1)$]

\item The arc $u$ crosses both $\tau u$ and $\tau^-u$.

\item The pair $(u, v)$ is crossing if and only if $(\tau u, \, \tau v)$ is crossing.

\end{enumerate}

\end{Cor}

\medskip

We shall see that the connecting arcs in $\B_\infty$ are very special. Indeed, we shall need a partial order on them.

\medskip

\begin{Lemma}\label{poset}

The set of connecting arcs in $\B_\infty$ is partially ordered in such a way that $[\ml_i, \mr_j]\le [\ml_r, \mr_s]$ if and only if $i\le r$ and $j\ge s.$

\end{Lemma}

\medskip

The following notion is the central objective of study in this section.

\medskip

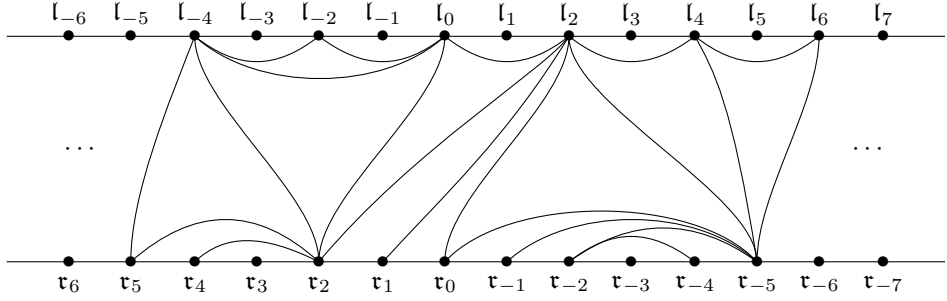
\begin{figure}[h] \label{Figure5.5}
  \centering
  \begin{tikzpicture}[xscale=2.50,yscale=1.5]

    \path (-2.26,0) node{$\cdots$};
    \path (1.93,0) node{$\cdots$};

    \draw (-2.66,1) -- (2.33,1);
    \draw (-2.66,-1) -- (2.33,-1);

   \node at (-2.33,1.2) {$\ml_{-6}$};
    \node at (-2.0,1.2) {$\ml_{-5}$};
    \node at (-1.66,1.2) {$\ml_{-4}$};
    \node at (-1.33,1.2) {$\ml_{-3}$};
    \node at (-1,1.2) {$\ml_{-2}$};
 \node at (-0.66,1.2) {$\ml_{-1}$};
    \node at (-0.33,1.2) {$\ml_{0}$};
    \node at (0,1.2) {$\ml_{1}$};
    \node at (0.33,1.2) {$\ml_2$};
    \node at (0.66,1.2) {$\ml_3$};
 \node at (1,1.2) {$\ml_{4}$};
    \node at (1.33,1.2) {$\ml_{5}$};
    \node at (1.66,1.2) {$\ml_{6}$};
    \node at (2,1.2) {$\ml_{7}$};

    \node at (-2.33,-1.2) {$\mr_{6}$};
    \node at (-2.0,-1.2) {$\mr_{5}$};
    \node at (-1.66,-1.2) {$\mr_{4}$};
    \node at (-1.33,-1.2) {$\mr_3$};
    \node at (-1,-1.2) {$\mr_2$};
 \node at (-0.66,-1.2) {$\mr_{1}$};
    \node at (-0.33,-1.2) {$\mr_{0}$};
    \node at (0,-1.2) {$\mr_{-1}$};
    \node at (0.33,-1.2) {$\mr_{-2}$};
    \node at (0.66,-1.2) {$\mr_{-3}$};
 \node at (1,-1.2) {$\mr_{-4}$};
    \node at (1.33,-1.2) {$\mr_{-5}$};
    \node at (1.66,-1.2) {$\mr_{-6}$};
    \node at (2,-1.2) {$\mr_{-7}$};

    \node at (-2.33,1.00) {$\bullet$};
    \node at  (-2.00,1.00) {$\bullet$};
    \node at  (-1.66,1.0){$\bullet$};
    \node at  (-1.33,1.0){$\bullet$};
    \node at  (-1.00,1.0){$\bullet$};
    \node at  (-0.66,1.0){$\bullet$};
    \node at  (-0.33,1.0){$\bullet$};
    \node at  (0.00,1.0){$\bullet$};
    \node at  (0.33,1.0){$\bullet$};
    \node at  (0.66,1.0){$\bullet$};
    \node at  (1.00,1.0){$\bullet$};
    \node at  (1.33,1.0){$\bullet$};
    \node at  (1.66,1.0){$\bullet$};
    \node at  (2.00,1.0){$\bullet$};

    \node at (-2.33,-1.0){$\bullet$};
    \node at (-2.00,-1.0){$\bullet$};
    \node at (-1.66,-1.0){$\bullet$};
   \node at (-1.33,-1.0){$\bullet$};
   \node at (-1.00,-1.0){$\bullet$};
    \node at (-0.66,-1.0){$\bullet$};
    \node at (-0.33,-1.0){$\bullet$};
    \node at (0.00,-1.0){$\bullet$};
    \node at (0.33,-1.0){$\bullet$};
    \node at(0.66,-1.0){$\bullet$};
    \node at (1.00,-1.0){$\bullet$};
   \node at (1.33,-1.0){$\bullet$};
    \node at(1.66,-1.0){$\bullet$};
    \node at (2.00,-1.0){$\bullet$};

    \draw (-2.0,-1) .. controls (-2.00,-0.4) and (-1.8,0.4) .. (-1.66,1);
    \draw (-2.00,-1) .. controls (-1.66,-0.5) and (-1.33,-0.5) .. (-1.00,-1);
\draw (-1.66, -1) .. controls (-1.5, -0.75) and (-1.3, -0.75) .. (-1.0, -1);

    \draw (-1.66,1) .. controls (-1.66,0.4) and (-1.00,-0.4) .. (-1.00,-1);
    \draw (-1.66,1) .. controls (-1.33,0.5) and (-0.66,0.5) .. (-0.33,1);
    \draw (-1.66,1) .. controls (-1.43,0.7) and (-1.23,0.7) .. (-1.00,1);
    \draw (-1.00,1) .. controls (-0.76,0.7) and (-0.56,0.7) .. (-0.33,1);

    \draw (-1.00,-1) .. controls (-1.00,-0.4) and (-0.33,0.4) .. (-0.33,1);
    \draw (-0.66,-1) .. controls (-0.33,-0.4) and (0,0.1) .. (0.33,1);

    \draw (0.33,1) .. controls (0.2,0.7) and (-0.5,0) .. (-1.00,-1);
    \draw (-0.33,1) .. controls (-0.1,0.7) and (0.1,0.7) .. (0.33,1);

    \draw (-0.33,-1) .. controls (-0.33,-0.4) and (0.33,0.4) .. (0.33,1);
    \draw (-0.33,-1) .. controls (0.00,-0.4) and (1.00,-0.4) .. (1.33,-1);
    \draw (0.00,-1) .. controls (0.33,-0.5) and (1.00,-0.5) .. (1.33,-1);
    \draw (0.33,-1) .. controls (0.66,-0.6) and (1,-0.6) .. (1.33,-1);
    \draw (0.33,-1) .. controls (0.56,-0.7) and (0.76,-0.7) .. (1.00,-1);

    \draw (0.33,1) .. controls (0.33,0.4) and (1.33,-0.4) .. (1.33,-1);
    \draw (1,1) .. controls (1.1,0.7) and (1.2,0.4) .. (1.33,-1);
    \draw (0.33,1) .. controls (0.56,0.7) and (0.76,0.7) .. (1.00,1);
    \draw (1.00,1) .. controls (1.23,0.7) and (1.43,0.7) .. (1.66,1);

    \draw (1.33,-1) .. controls (1.33,-0.4) and (1.66,0.4) .. (1.66,1);

  \end{tikzpicture}

  \caption{A triangulation of $\mathcal{B}_\infty$}

\label{fig:triangulation1} \end{figure}

\begin{Defn}

A maximal set $\bT$ of pairwise non-crossing arcs in $\mathcal{B}_\infty$ is called a \emph{triangulation} of $\mathcal{B}_\infty$; see Figure \ref{Figure5.5}. In this case, $C(\bT)$ denotes the set of connecting arcs of $\bT$.

\end{Defn}

\medskip

It is easy to see that two connecting arcs are comparable with respect to the partial order defined in Lemma \ref{poset} if and only if they do not cross each other. This gives us immediately the following useful observation.

\medskip

\begin{Lemma}

Let $\bT$ be a triangulation of $\B_\infty$. If $\bT$ contains some connecting arcs, then $C(\bT)$ is well ordered.

\end{Lemma}

\medskip

We shall also need the following easy result.

\medskip

\begin{Lemma} \label{nca}

Let $\bT$ be a triangulation of $\B_\infty$, and let $p$ be an integer.

\vspace{-0.5pt}

\begin{enumerate}[$(1)$]

\item If there exist infinitely many $i<p$ such that $[\ml_i, \ml_{j_i}]\in \bT\vspace{1pt}$ for some $j_i\ge p$ or infinitely many $j>-p$ such that $[\ml_{i_j}, \mr_j]\in \bT\vspace{1pt}$ for some $i_j\ge p$, then no $\ml_i$ with $i<p$ is an endpoint of an arc of $C(\bT)$.

\item If there exist infinitely many $i>p$ such that $[\ml_{j_i}, \ml_i]\in \bT\vspace{1pt}$ for some $j_i\le p$ or infinitely many $j<-p$ such that $[\ml_{i_j}, \mr_j]\in \bT\vspace{1pt}$ for some $i_j\le p$, then no $\ml_i$ with $i>p$ is an endpoint of an arc of $C(\bT)$.

\end{enumerate} \end{Lemma}

\noindent {\it Proof.} We shall prove only Statement (1). Consider a connecting arc $v=[\ml_r, \mr_s]$ with $r<p$. If the first situation in Statement (1) occurs, then there exists some integer $i<r$ such that $[\ml_i, \ml_{j_i}]\in \bT$ for some $j_i\ge p$. In this case, $v$ crosses $[\ml_i, \ml_{j_i}]$, and hence, $v\notin \bT$. If the second situation occurs, then there exists some $j>s$ such that $[\ml_{i_j}, \mr_j]\in \bT$ for some $i_j\ge p$. In this case, $v$ crosses $[\ml_{i_j}, \ml_j]$, and hence, $v\notin \bT$. The proof of the lemma is completed.

\medskip

\noindent{\sc Remark.} A similar statement holds for lower marked points.

\medskip

Let $\bT$ be a triangulation of $\B_\infty$. For each  $u\in {\rm arc}(\B_\infty)$, we shall denote by $\mathbb{T}_u$ the set of arcs of $\mathbb{T}$ which cross $u$.

\medskip

\begin{Lemma} \label{inf-crossing}

Let $\bT$ be a triangulation of $\B_\infty$ containing some connecting arcs, and let $u$ be an arc in $\B_\infty$. If $\bT_u$ is infinite, then some marked point in $\B_\infty$ is an endpoint of infinitely many arcs in $\bT_u$.

\end{Lemma}

\noindent{\it Proof.} Assume that every marked point in $\B_\infty$ is an endpoint of at most finitely many arcs in $\bT_u$. We shall prove that $\bT_u$ is finite. Consider first the case where $u$ is an upper arc, say $u=[\ml_r, \ml_s]$ with $r<s-1$. Then every arc in $\bT_u$ has as an endpoint some marked point $\ml_i$ with $r<i<s$. Thus, it follows from the assumption that $\bT_u$ is finite. The case where $u$ is a lower arc can be treated in a similar manner.

Consider finally the case where $u$ is a connecting arc, say $u=[\ml_r, \mr_s]$. We claim that $\bT_u$ contains at most finitely connecting arcs. Indeed, assume that $\bT_u$ contains some connecting arc $v_0=[\ml_p, \mr_q]$. Then $p>r$ and $q>s$, or $p<r$ and $q<s$. Suppose that the first situation occurs. If $v=[\ml_i, \mr_j]\in \bT_u$, then $v$ does not cross $v_0$, and hence, either $i> p$ and $q\ge j>s$ or $p\ge i>r$ and $j> q$. Thus, our claim follows from the assumption on the marked points. Similarly, our claim holds if the second situation occurs.

Suppose now that $\bT_u$ contains infinitely many upper arcs $u_i=[\ml_{r_i}, \ml_{s_i}]$ with $r_i<r<s_i$, for $i=1, 2 \cdots.$ By the hypothesis on the marked points, we may assume that $r_{i+1}<r_i,$ for all $i\ge 1$. Since the $u_i$ do not cross each other, we obtain $s_i<s_{i+1}$ for all $i\ge 1$. As a consequence, no upper marked point in $\B_\infty$ is an endpoint of some connecting arc of $\bT$, a contradiction to the hypothesis stated in the lemma. Similarly, $\bT_u$ contains at most finitely many lower arcs. That is, $\bT_u$ is finite. The proof of the lemma is completed.

\medskip

We say that an upper marked point $\ml_i$ in $\B_\infty$ is {\it covered} by an upper arc $[\ml_r, \ml_s]$ if $r < i < s$; and a lower marked point $\mr_j$ is {\it covered} by a lower arc $[\mr_p, \mr_q]$ if $p> j > q$.

\medskip

\begin{Lemma} \label{WC-NCA}

Let $\bT$ be a triangulation of $\B_\infty$. If $C(\bT)$ is empty, then one of the following situations occurs.

\vspace{-1pt}

\begin{enumerate}[$(1)$]

\item Every upper marked point in $\B_\infty$ is covered by infinitely many upper arcs of $\bT$ and is an endpoint of at most finitely many upper arcs of $\bT$.

\vspace{1pt}

\item Every lower marked point in $\B_\infty$ is covered by infinitely many lower arcs of $\bT$ and is an endpoint of at most finitely many lower arcs of $\bT$.

\end{enumerate}

\end{Lemma}

\noindent{\it Proof.} Let $C(\bT)=\emptyset$. Assume first that each upper marked point $\ml_t$ with $t\in \Z$ is covered by an upper arc $u_t=[\ml_{i_t}, \ml_{j_t}]$ of $\bT$. We shall show that Statement (1) holds. Suppose on the contrary that some upper marked point, say $\ml_0$, is covered by only finitely many upper arcs of $\bT$. Let $r <0$ be minimal such that $\ml_r$ is an endpoint of some arc $v_0$ covering $\ml_0$ of $\bT$. Then $v_0=[\ml_r, \ml_s]$ with $r< 0 < s.$ By our assumption, $\ml_r$ is covered by an upper arc $[\ml_p, \ml_q]\in \bT$ with $p<r<q$. Since $[\ml_p, \ml_q]$ does not cross $[\ml_r, \ml_s]$, we obtain $s\le q$. That is, $[\ml_p, \ml_q]$ covers $\ml_0$, contrary to the minimality of $r$. This establishes the first part of Statement (1). Next, given any $t\in \Z$, each of the $[\ml_i, \ml_t]$ with $i<i_t$ and the $[\ml_t, \ml_j]$ with $j>j_t$ crosses the upper arc $u_t\in \bT$, and hence, it does not belong to $\bT$. Thus, $\ml_t$ is an endpoint of at most finitely many upper arcs of $\bT$. This establishes Statement (1).

Assume now that some upper marked point, say $\ml_0$, is not covered by any upper arc of $\bT$. For any $j\in \Z$, the connecting arc $[\ml_0, \mr_j]$ does not belong to $\bT$, and hence, it crosses some arc $v\in \bT$. By the hypothesis on $\ml_0$, the arc $v$ is a lower arc which covers the lower marked point $\mr_j$. That is, every lower marked point is covered by a lower arc of $\bT$. Using a similar argument as above, we are able to show that Statement (2) holds. The proof of the lemma is completed.

\medskip

Let $\bT$ be a triangulation of $\B_\infty$.
An upper marked point $\mf{l}_p$ is called {\it left $\bT$-bounded} if $[\ml_i, \ml_p]\notin \bT$ for all but finitely many $i<p$ and $[\ml_p, \mr_j]\notin \bT$ for all but finitely many $j>-p\,;$ and {\it left $\bT$-unbounded} if $[\ml_i, \ml_p], [\ml_p, \mr_j]\in \bT$ for infinitely many $i<p$ and infinitely many $j>-p$.

\medskip

By symmetry, $\ml_p$ is called {\it right $\bT$-bounded} if $[\ml_p, \ml_i]\not\in \bT$ for all but finitely many $i>p$ and $[\ml_p, \mr_j]\notin \bT$ for all but finitely many $j<-p\,;$ and {\it right $\bT$-unbounded} if $[\ml_p, \ml_i], [\ml_p, \mr_j]\in\bT$ for infinitely many $i>p$ and infinitely many $j<-p$.

\medskip

In a similar manner, we shall define a lower marked point to be {\it left $\bT$-bounded}, {\it left $\bT$-unbounded},
{\it right $\bT$-bounded}, and {\it right $\bT$-unbounded}.

\medskip

\begin{Lemma} \label{min-arc}

Let $\bT$ be a triangulation of $\B_\infty$ with $[\ml_p, \mr_q]\in C(\bT)$.

\vspace{0pt}

%
%
%
%
%
%
%
%

\begin{enumerate}[$(1)$]

\item If $\ml_p$ is left $($respectively, right$)$ $\bT$-bounded, then some of the $\ml_i$ with $i<p$ $($respectively, $i>p)$ is an endpoint of some arc of $C(\bT)$.

\vspace{2pt}

\item If $\mr_q$ is left $($respectively, right$)$ $\bT$-bounded, then some of the $\mr_j$ with $j> q$ $($respectively, $j<q)$ is an endpoint of some arc of $C(\bT)$.

\end{enumerate}\end{Lemma}

\noindent {\it Proof.} We shall prove only the first part of Statement (1). Assume that none of the $\ml_i$ with $i<p$ is an endpoint of an arc of $C(\bT)$, moreover, $[\ml_p, \mr_j]\in \bT$ for at most finitely many $j> -p$. Thus, we may assume with no loss of generality that $q$ is maximal such that $[\ml_p, \mr_q]\in \bT$. We shall need to prove that
$[\ml_i, \ml_p]\in \bT$ for infinitely many integers $i<p$.

Indeed, in view of the assumption on $\ml_p$, we see that $[\ml_p, \mr_q]$ is a minimal element in $C(\bT)$. Since $[\ml_{p-1}, \mr_q]\notin \bT$ by the assumption, $\bT_{[\ml_{p-1}, \mr_q]}$ contains some arc $v$. Crossing $[\ml_{p-1}, \mr_q]$ but not $[\ml_p, \mr_q]$, the arc $v$ is neither a lower arc nor a connecting arc greater than $[\ml_p, \mr_q].$ Then, $v$  is not a connecting arc by the minimality of $[\ml_p, \mr_q]$. As a consequence, $v=[\ml_i, \ml_p]$ for some $i<p-1$. This proves that $\bT_{[\ml_{p-1}, \mr_q]}$ contains only upper arcs of the form $[\ml_i, \ml_p]$ with $i<p-1.$

\vspace{2pt}

Thus, it suffices to show that $\bT_{[\ml_{p-1}, \mr_q]}$ is infinite. If this is not the case, then there exists a minimal $ r <p-1$ for which $[\ml_r, \ml_p]\in \bT_{[\ml_{p-1}, \mr_q]}$. Not belonging to $\bT$, the connecting arc $[\ml_r, \mr_q]$ crosses some arc $w$ of $\bT$. As argued above, $w$ is neither a lower arc nor a connecting arc. Thus $w=[\ml_s, \ml_t]$ with $s<r <t$. If $t\le p-1$, then $w$ crosses $[\ml_r, \ml_p]$, a contradiction. Thus $p-1<t$, and hence, $w\in \bT_{[\ml_{p-1}, \mr_q]}\vspace{1pt}.$ As we have shown, $t=p$, a contradiction to the minimality of $r$. The proof of the lemma is completed.

\medskip

Let $\Sa$ be a set of arcs in $\B_\infty$. We shall denote by $\tau \Sa$ the set of arcs of the form $\tau u$ with $u\in \Sa;$ and by $\tau^-\Sa$ the set of arcs of the form $\tau^-v$ with $v\in \Sa.$

\medskip

\begin{Defn}\label{compact-set}

A set $\Da$ of arcs in $\B_\infty$ is called {\it compact} if it admits a finite subset $\Sa$ such that every arc in $\Da$ crosses some arc of $\tau \Sa$ and some arc of $\tau^-\Sa$.

\end{Defn}

\medskip

By definition, the empty subset of ${\rm arc}(\B_\infty)$ is compact. Moreover, as shown below, any finite set of arcs in $\B_\infty$ is compact.

\medskip

\begin{Lemma}\label{Cofinite-compact}

Let $\Da$ be a set of arcs in $\B_\infty$, and let $\Oa$ be a co-finite subset of $\Da$. If $\Oa$ is compact, then $\Da$ is compact.

\end{Lemma}

\noindent{\it Proof.} Assume that $\Oa$ is compact. Let $\Sa$ be a finite subset of $\Oa$ satisfying the condition stated in Definition \ref{compact-set}.
Let $\Ta$ be the union of $\Sa$ and the complement of $\Oa$ in $\Da$, which is finite by the hypothesis. It is evident that $\tau \Sa\subseteq \tau \Ta$ and $\tau^-\Sa\subseteq \tau^-\Ta$. Let $u$ be an arc in $\B_\infty$. If $u\in \Oa$, then it crosses some arc of $\tau \Sa$ and some arc of $\tau^-\Sa$. Otherwise, $u\in \Ta$, and by Corollary \ref{tau-cross}, $u$ crosses $\tau u$ and $\tau^-u$. The proof of the lemma is completed.

\medskip

The following notion is important for us to determine the cluster-tilting subcategories in the next section.

\medskip

\begin{Defn}

A triangulation $\bT$ of $\B_\infty$ is called {\it compact} provided that $\bT_u$ is compact, for every arc $u$ in $\B_\infty$.

\end{Defn}

\medskip

\begin{Lemma} \label{Comapct-2}

Let $\bT$ be a compact triangulation of $\B_\infty$, and let $p$ be an integer.

\vspace{-1pt}

\begin{enumerate}[$(1)$]

\item If $[\ml_i, \ml_p]\in \bT$ for infinitely many integers $i<p$, then $\ml_p$ is left $\bT$-unbounded.


\vspace{1pt}

\item If $[\ml_p, \ml_i]\in \bT$ for infinitely many integers $i>p$, then $\ml_p$ is right $\bT$-unbounded.


\vspace{1pt}

\item If $[\mr_j, \mr_p]\in \bT$ for infinitely many integers $j>p$, then $\mr_p$ is left $\bT$-unbounded.

\vspace{1pt}

\item If $[\mr_p, \mr_j]\in \bT$ for infinitely many integers $j<p$, then $\mr_p$ is right $\bT$-unbounded.

\end{enumerate}

\end{Lemma}

\noindent{\it Proof.} We shall prove only Statement (1). Assume that $[\ml_i, \ml_p]\in \bT$ for infinitely many $i<p$. We shall need to show that $[\ml_p, \mr_j]\in \bT$ for infinitely many integers $j>-p$. Suppose on the contrary that this is not the case. Then, there exists an integer $q$ such that $[\ml_p, \mr_j]\notin \bT$ for every $j>q$. Consider the connecting arc $u=[\ml_{p-1}, \mr_q]$. By the assumption, $[\ml_i, \ml_p]\in \bT_{u}$ for infinitely many $i<p-1$. Since $\bT_{u}$ is compact, there exists a finite subset $\Sa$ of $\bT_{u}$ satisfying the condition stated in Definition \ref{compact-set}. Let $t<p-1$ be minimal such that $u_0=[\ml_t, \ml_p]\in \Sa$. Observe that $\bT_{u}$ contains some arc $w=[\ml_r, \ml_p]$ with $r<t$.

We claim that $w$ does not cross $\tau^-v$ for any arc $v\in \Sa$. Indeed, this is trivially the case if $v$ is a lower arc. Suppose that $v$ is a connecting arc in $\Sa$. Then $v=[\ml_m, \mr_n]$ with $m>p-1$ and $n>q$, or else, $m<p-1$ and $n<q$. Since $v$ does not cross any of the infinitely many arcs $[\ml_i, \ml_p]$ in $\bT_{u}$ with $i<p-1$, we obtain $m>p-1$ and $n>q$. By the assumption on $q$, we obtain $m>p$, and hence, $w=[\ml_r, \ml_p]$ does not cross $\tau^-v=[\ml_{m-1}, \mr_{n-1}]$.
Suppose now that $v$ is an upper arc in $\Sa$. Then $v=[\ml_m, \ml_n]$ with $m<p-1<n$. Since $v$ does not cross any of the infinitely many arcs $[\ml_i, \ml_p]\in \bT_{u}$ with $i<p$, we have $n=p$, that is, $v=[\ml_m, \ml_p]$ with $m<p-1$. By the minimality of $t$, we obtain $t\le m$, and hence, $w=[\ml_r, \ml_p]$ does not cross $\tau^-v=[\ml_{m-1}, \ml_{p-1}]$. This establishes our claim, which is contrary to the property stated in Definition \ref{compact-set}. The proof of the lemma is completed.

\medskip

The following result exhibits a characteristic property of a compact triangulation. 

\medskip

\begin{Prop} \label{DI-chain}

Let $\bT$ be a triangulation of $\B_\infty$. If $\bT$ is compact, then $C(\bT)$ is a double-infinite chain
as follows$\,:$
$$\cdots < u_{-2} < u_{-1} < u_0 < u_1 < u_2 < \cdots
\vspace{2pt}$$

\end{Prop}

\noindent{\it Proof.} Let $\bT$ be compact. First, we claim that $C(\bT)$ is non-empty. Assume that this is not the case. By Lemma \ref{WC-NCA}, we may assume that every upper marked point is covered by infinitely many upper arcs of $\bT$ and is an endpoint of at most finitely many upper arcs of $\bT$. Consider the connecting arc $u_0=[\ml_0, \mr_0].$ Since every upper arc covering $\ml_0$ crosses $u_0$, the set $U(\bT_{u_0})$ of upper arcs of $\bT_{u_0}$ is infinite. Being compact, $\bT_{u_0}$ has a finite subset $\Sa$ satisfying the condition stated in Definition \ref{compact-set}. Let $r_0, s_0$ be the least and the greatest, respectively, such that $\ml_{r_0}, \ml_{s_0}$ are endpoints of some arcs of $\Sa$. Since each upper marked point is an endpoint of at most finitely many arcs of $U(\bT_{u_0})$, the set $\bT_{u_0}$ contains an upper arc $u_1=[\ml_{r_1}, \ml_{s_1}]$ with $r_1 < r_0-1$ and $s_1> s_0+1$. Given any arc $v\in \Sa$, by the assumption, either $v=[\ml_r, \ml_s]$ with $r_0\le r <s \le s_0$ or $v$ is a lower arc. In either case, $u_1$ does not cross $\tau^-v$ or $\tau v$, a contradiction. This establishes our claim.

Next, suppose on the contrary that $C(\bT)$ contains a minimal element $[\ml_p, \mr_q]$. Since the arcs of $\bT$ do not cross each other, we deduce from the minimality of $[\ml_p, \mr_q]$ that $[\ml_i, \mr_j]\not\in \bT$ for all $i, j$ with $i<p$. By Lemmas \ref{min-arc} and \ref{Comapct-2}, $[\ml_p, \mr_j]\in \bT$ for some $j>q$, contrary to the minimality of $[\ml_p, \mr_q]$. Similarly, one can show that $C(\bT)$ contains no maximal element. Being well ordered and interval-finite, $C(\bT)$ is a double infinite chain. The proof of the lemma is completed.

\medskip

Let $\bT$ be a triangulation of $\B_\infty$. A marked point $\mf{p}$ in $\B_\infty$ is called a {\it left $\bT$-fountain base} if $\mathfrak{p}$ is left $\bT$-unbounded but right $\bT$-bounded. In this case, if $\mf{p} = \ml_p$, then the set of arcs of the form $[\ml_i, \ml_p]$ with $i<p$ and $[\ml_p, \mr_{j}]$ with $j>-p$ is called a {\it left fountain} of $\bT$ at $\mf{p}$; and if $\mf{p} = \mr_q$, then the set of arcs of the form $[\mr_i, \mr_q]$ with $i > q$ and $[\ml_j, \mr_{q}]$ with $j<-q$ is called a {\it left fountain} of $\bT$ at $\mf{p}.$

\medskip

Similarly, $\mf{p}$ is said to be a {\it right $\bT$-fountain base} if $\mathfrak{p}$ is right $\bT$-unbounded but left $\bT$-bounded. In this case, if $\mf{p} = \ml_p$, then the set of arcs of the form $[\ml_p, \ml_i]$ with $p < i$ and $[\ml_p, \mr_{j}]$ with $j<-p$ is called a {\it right fountain} of $\bT$ at $\mf{p}$; and if $\mf{p} = \mr_q$, then the set of arcs of the form $[\mr_q, \mr_i]$ with $q > i$ and $[\ml_j, \mr_{q}]$ with $j>-q$ is called a {\it right fountain} of $\bT$ at $\mf{p}.$

\medskip

Furthermore, $\mf{p}$ is called a {\it full $\bT$-fountain base} if $\mf{p}$ is left and right $\bT$-unbounded; and in this case, the set of arcs of $\bT$ which have $\mf{p}$ as an endpoint is called a {\it full fountain} of $\bT$ at $\mf{p}.$

\medskip

For brevity, a marked point is called a $\bT$-{\it fountain base} if it is a left, right or full $\bT$-fountain base; and a left, right or full fountain of $\bT$ will be simply called a {\it fountain}. If $\mf{p}$ is a $\bT$-fountain base, then the fountain at $\mathfrak{p}$ will be denoted by $\mathbb{F}_{\hspace{-1.5pt}_{\bT}}(\mathfrak{p})$.

\medskip

\begin{Lemma}\label{Full-fountain}

Let $\bT$ be a triangulation  of $\B_\infty$.

\vspace{-2pt}

\begin{enumerate}[$(1)$]

\item If $\mf{p}$ is a full $\bT$-fountain base in $\B_\infty$, then it is the unique $\bT$-fountain base and is an endpoint of all connecting arcs of $\bT$.

\item If $\mf{p}$, $\mf{q}$ are two distinct $\bT$-fountain bases in $\B_\infty$, then they are the only $\bT$-fountain bases with one being a left $\bT$-fountain base and the other one being a right $\bT$-fountain base.

\end{enumerate}

\end{Lemma}

\noindent{\it Proof.} Assume that some $\ml_{p}$ is left $\bT$-unbounded. We claim that $\ml_p$ is the only left $\bT$-unbounded marked point, and none of the $\ml_i$ with $i<p$ is an endpoint of any connecting arc in $\B_\infty$. Indeed, the second part of this claim follows from Lemma \ref{nca}(1). As a consequence, none of the $\ml_i$ with $i<p$ and the $\mr_j$ with $j\in \Z$ is left $\bT$-unbounded. Now, choose arbitrarily a connecting arc $[\ml_p, \mr_q]$ of $\bT$. Since the arcs of $\bT$ do not cross each other, $[\ml_i, \ml_j]\not\in \bT$ for all $i, j$ with $i<p<j.$ In particular, no $\ml_j$ with $j>p$ is left $\bT$-unbounded. This establishes the claim.

Suppose now that $\mf{p}$ is a full $\bT$-fountain base. We shall consider only the case where $\mf{p}$ is an upper marked point, say $\mf{p}=\ml_p$. It follows from the above claim and its dual version that $\mf{p}$ is the only $\bT$-fountain base. Moreover,
none of the $\ml_i$ with $i<p$ or $i>p$ is an end-point of some connecting arc of $\bT$. As a consequence, $\mf{p}$ is an endpoint of all connecting arcs of $\bT$. This establishes Statement (1). Similarly, Statement (2) follows immediately from the above claim and its dual version. The proof of the lemma is completed.

\medskip

Next we shall find some sufficient conditions for a triangulation to be compact.

\medskip

\begin{Lemma} \label{LemmaCrossFountain}

Let $\bT$ be a triangulation of $\B_\infty$, and let $v$ be an arc in $\B_\infty$. If $v$ crosses infinitely many arcs of a full fountain of $\bT$, then $\bT_v$ is compact.

\end{Lemma}

\noindent{\it Proof.} Assume that $v$ crosses infinitely many arcs of a full fountain $\mathbb{F}_{\hspace{-1.5pt}_{\bT}}(\mathfrak{p})$ of $\bT$. We shall consider only the case where $\mathfrak{p}$ is an upper marked point, say $\mathfrak{p} = \ml_p$ for some $p\in \Z$. In particular, $v$ is not be a lower arc.

Assume first that $v$ is an upper arc. It is easy to see that $v=[\ml_r, \ml_s]$ with $r<p<s$. Let $i_0 < r$ be maximal such that $v_1=[\ml_{i_0}, \ml_p]\in \bT$, and let $j_0>s$ be minimal such that $v_2=[\ml_p, \ml_{j_0}]\in \bT$.
Suppose that $u$ is an arc lying in $\bT_v$ but not in $\mathbb{F}_{\hspace{-1.5pt}_{\bT}}(\ml_{p})$. Then $u$ is not a lower arc, and by Lemma \ref{Full-fountain}(1), it is an upper arc. Since $u$ does not cross $v_1$ or $v_2$, we see that $u=[\ml_i, \ml_j]$ with $i_0\le i<r<j<p$ or $p<i<s<j\le j_0$. Therefore, $\bT_v\cap \mathbb{F}_{\hspace{-1.5pt}_{\bT}}(\mf{p})$ is co-finite in $\bT_v$.

We claim that $\bT_v\cap \mathbb{F}_{\hspace{-1.5pt}_{\bT}}(\mf{p})$ is compact. Indeed, $v_1, v_2 \in \bT_v\cap \mathbb{F}_{\hspace{-1.5pt}_{\bT}}(\ml_p)$ with $\tau v_1=[\ml_{i_0+1}, \ml_{p+1}]$ and $\tau^-v_2=[\ml_{p-1}, \ml_{j_0-1}]$. Let $u\in \bT_v\cap \mathbb{F}_{\hspace{-1.5pt}_{\bT}}(\ml_p)$. If $u$ is an upper arc, then $u=[\ml_m, \ml_p]$ with $m<r$ or $u=[\ml_p, \ml_n]$ with $n>s$. By the maximality of $i_0$ and the minimality of $j_0$, we obtain $u=[\ml_m, \ml_p]$ with $m\le i_0$ or $u=[\ml_p, \ml_n]$ with $j_0\le n$.
In the first situation, since $m<i_0+1<r+1\le p<p+1$ and $m<r\le p-1<p<s\le j_0-1$, we see that $u$ crosses both $\tau v_1$ and $\tau^-v_2$. In the second situation, since $i_0+1\le r<p<p+1\le s<n$ and $p<s\le j_0-1<n$, we see that $u$ crosses both $\tau v_1$ and $\tau^-v_2$. This establishes our claim, and hence, $\bT_v$ is compact by Lemma \ref{Cofinite-compact}.

Consider next the case where $v$ is a connecting arc, say $v = [\ml_r, \mr_s]$. Then $r\ne p$. We shall consider only the case where $r < p$. Let $i_0 < r$ be maximal such that $v_1=[\ml_{i_0}, \ml_{p}] \in \bT$ and let $j_0 > s$ be minimal such that $v_2=[\ml_{p}, \mr_{j_0}] \in \bT$. Moreover, since $\ml_p$ is right $\bT$-unbounded, there exists a maximal integer $t_0< s$ such that $w_0=[\ml_p, \mr_{t_0}]\in \bT$. Let $u$ be an arc lying in $\bT_v$ but not in $\mathbb{F}_{\hspace{-1.5pt}_{\bT}}(\mathfrak{p})$. By Lemma \ref{Full-fountain}(1), $u$ is not a connecting arc. Since $u$ does not cross any of the arcs $v_1, v_2, w_0$,  in case $u$ is a lower arc, we have $u=[\mr_i, \mr_j]\in \bT$ with $j_0 \ge i > s > j \ge t_0$; and otherwise, $u=[\ml_i, \ml_j]\in \bT$ with $i_0\le i<r<j< p$. Therefore, $\bT_v\cap \mathbb{F}_{\hspace{-1.5pt}_{\bT}}(\mf{p})$ is co-finite in $\bT_v$.

By Lemma \ref{Cofinite-compact}, it remains to show that $\bT_v\cap \mathbb{F}_{\hspace{-1.5pt}_{\bT}}(\mf{p})$ is compact. Indeed, $v_1, v_2\in \bT_v\cap \mathbb{F}_{\hspace{-1.5pt}_{\bT}}(\mf{p})$ with $\tau v_1=[\ml_{i_0+1}, \ml_{p+1}]$ and $\tau^-v_2=[\ml_{p-1}, \mr_{j_0-1}]$. Let $u$ be an arbitrary arc in $\bT_v\cap \mathbb{F}_{\hspace{-1.5pt}_{\bT}}(\ml_{p}).$ Then $u$ is not a lower arc.
If $u$ is a connecting arc, then $u=[\ml_p, \mr_j]$ with $j>s$. By the minimality of $j_0$, we have $j\ge j_0$.
In this case, $u$ crosses both $\tau v_1=[\ml_{i_0+1}, \ml_{p+1}]$ and $\tau^-v_2=[\ml_{p-1}, \mr_{j_0-1}]$.
If $u$ is an upper arc, then it follows from the maximality of $i_0$ that $u=[\ml_i, \ml_p]$ with $i\le i_0$. It is then easy to see that $u$ crosses $\tau v_1$ and $\tau^-v_2$. The proof of the lemma is completed.

\medskip

Given a marked point $\mf{p}$ in $\B_\infty$, we shall denote by $\mathbb{E}_{\hspace{-0.5pt}_{\bT}}(\mf{p})$ the set of arcs of $\bT$ which have $\mf{p}$ as an endpoint.

\medskip

\begin{Lemma} \label{LemmaCrossFountain2}

Let $\bT$ be a triangulation of $\B_\infty$ with $\mathfrak{p}$ a left or right $\bT$-fountain base, and let $v$ be an arc in $\B_\infty$. If $v$ crosses infinitely many arcs of $\mathbb{F}_{\hspace{-1.5pt}_{\bT}}(\mathfrak{p})$, then $\bT_v \cap \mathbb{F}_{\hspace{-1.5pt}_{\bT}}(\mathfrak{p})$ is compact and co-finite in $\mathbb{E}_{\hspace{-0.5pt}_{\bT}}(\mf{p})$.

\end{Lemma}

\noindent{\it Proof.} We shall consider only the case where $\mathfrak{p}=\ml_p$ for some $p\in \Z$, which is a left $\bT$-fountain base. Assume that $v$ crosses infinitely many arcs of $\mathbb{F}_{\hspace{-1.5pt}_{\bT}}(\mathfrak{p})$.
Since every lower arc crosses at most finitely many arcs of $\mathbb{E}_{\hspace{-0.5pt}_{\bT}}(\ml_p)$, the arc $v$ is not a lower arc. Moreover, since $\bT$ is right $\bT$-bounded, $v$ has as an endpoint some marked point $\ml_r$ with $r<p$. That is, $v=[\ml_r, \ml_s]$ with $p<s$ or $v=[\ml_r, \mr_s]$.

Let $w$ be an arc of $\mathbb{F}_{\hspace{-1.5pt}_{\bT}}(\mathfrak{p})$, which does not cross $v$. If $v=[\ml_r, \ml_s]$ with $r<p<s$, then $w=[\ml_j, \ml_p]$ with $r\le j < p-1$. If $v=[\ml_r, \mr_s]$, then $w=[\ml_j, \ml_p]$ with $r \le j<p-1$ or $w=[\ml_i, \mr_t]$ with $r \le i < p-1$ and $-p < t \le s$. Therefore, all but finitely many arcs of $\mathbb{F}_{\hspace{-1.5pt}_{\bT}}(\mathfrak{p})$ cross $v$.

Now, let $m < r$ be maximal such that $v_1=[\ml_m, \ml_p] \in \bT$. It is clear that $v_1\in \bT_v\cap \mathbb{F}_{\hspace{-1.5pt}_{\bT}}(\ml_p)$. Let $u\in \bT_v\cap \mathbb{F}_{\hspace{-1.5pt}_{\bT}}(\ml_p).$ If $u$ is a connecting arc, then $u=[\ml_p, \mr_t]$ for some $t>-p$, which clearly crosses $\tau v_1=[\ml_{m+1}, \ml_{p+1}].$ Otherwise, $u=[\ml_t, \ml_p]$ for some $t<p-1$. By the maximality of $m$, we obtain $t\le m$, and thus, $u$ crosses $\tau v_1=[\ml_{m+1}, \ml_{p+1}].$

Next, in case $v=[\ml_r, \ml_s]$, let $n>-p$ be minimal such that $[\ml_p, \mr_n]\in \bT$; and in case $v=[\ml_r, \mr_s]$, let $n>{\rm max}\{-p, s\}$ be minimal such that $[\ml_p, \mr_n]\in \bT$. In either case, set $v_2=[\ml_p, \mr_n],$ which clearly belongs to $\bT_v \cap \mathbb{F}_{\hspace{-1.5pt}_{\bT}}(\ml_p)$. Let $u\in \bT_v\cap \mathbb{F}_{\hspace{-1.5pt}_{\bT}}(\ml_p).$ If $u$ is an upper arc, then $u=[\ml_t, \ml_p]$ with $t<r$, which crosses $\tau^-v_2=[\ml_{p-1}, \mr_{n-1}].$ Otherwise, $u=[\ml_p, \mr_t]$, where $t>-p$, and $t>s$ in case $v=[\ml_r, \mr_t]$. By the minimality of $n,$ we obtain $t\ge n$, and hence, $u$ crosses $\tau^-v_2=[\ml_{p-1}, \mr_{n-1}].$ Setting $\Sa=\{v_1, v_2\}$, we see that $\bT_v \cap \mathbb{F}_{\hspace{-1.5pt}_{\bT}}(\mathfrak{p})$ is compact. The proof of the lemma is completed.

\medskip

Let $\bT$ be a triangulation of $\B_\infty$. We shall say that a marked point in $\B_\infty$ is {\it $\bT$-bounded} if it is both left and right $\bT$-bounded, or equivalently, $\mf{p}$ is an endpoint of at most finitely many arcs of $\bT$.

\medskip

\begin{Lemma} \label{LemmaFountain}

Let $\bT$ be a triangulation of $\B_\infty$. If every marked point in $\B_\infty$ is either $\bT$-bounded or an endpoint of infinitely many arcs of $C(\bT)$, then every marked point is either $\bT$-bounded or a $\bT$-fountain base.

\end{Lemma}

\noindent{\it Proof.} Assume that every marked point is either $\bT$-bounded or an endpoint of infinitely many arcs of $C(\bT)$. Let $\mathfrak{p}$ be a marked point, which is an endpoint of infinitely many arcs of $C(\bT)$. We need to show that $\mf{p}$ is a $\bT$-fountain base. First of all, we claim that $C(\bT)$ is a double infinite chain. Indeed, suppose that $\C(\bT)$ has a least element $[\ml_m, \mr_n]$. Since the arcs of $\bT$ do not cross each other, $\,[\ml_i, \mr_j]\notin \bT$ for all $i, j$ with $i<m$ or $j>n$. By Lemma \ref{min-arc}, $[\ml_i, \ml_m]\in \bT$ for infinitely many $i<m;$ and $[\mr_j, \mr_n]\in \bT$ for infinitely many $j>n$. By the hypothesis stated in the lemma, $\C(\bT)$ has connecting arcs $[\ml_m, \mr_{n_0}]$ with $n_0\ne n$ and $[\ml_{m_0}, \mr_n]$ with $m_0\ne m$. By the minimality of $[\ml_m, \mr_n]$, we obtain $m>m_0$ and $n_0<n$. Then, $[\ml_m, \mr_{n_0}]$ crosses $[\ml_{m_0}, \mr_n]$, which is absurd. Similarly, we can show that $\C(\bT)$ has no greatest element. This establishes our claim.

For the rest of the proof, we shall consider only the case where $\mathfrak{p} = \ml_p$ such that $[\ml_p, \mr_j] \in \bT$ for infinitely many $j>-p$. We claim that $\ml_p$ is left $\bT$-unbounded, that is, $[\ml_i, \ml_p] \in \bT$ for infinitely many $i < p$. Indeed, assume that this is not the case. Define $s=p-1$ if $[\ml_j, \ml_p]\not\in \bT$ for every $j < p-1;$ and otherwise, let $s<p-1$ be minimal such that $[\ml_s, \ml_p] \in \bT$. By Lemma \ref{nca}(1), $\ml_s$ is not an endpoint of any arc of $C(\bT)$. By the hypothesis stated in the lemma, $\ml_s$ is an endpoint of at most finitely many arcs of $\bT$. Set $t=s-1$ if $[\ml_i, \ml_s]\notin \bT$ for every $i<s-1$; and otherwise, let $t<s-1$ be minimal such that $[\ml_t, \ml_s] \in \bT$. Consider the upper arc $v=[\ml_t, \ml_p]$ which, by Lemma \ref{nca}(1), does not cross any arc of $C(\bT)$. Suppose that $v$ crosses some upper arc $u$ in $\bT$. Since $u$ does not cross any connecting arc of $\bT$ having $\ml_p$ as an endpoint, $u = [\ml_{t_1}, \ml_{s_1}]$ with $t_1 < t < s_1 < p$. If $s < s_1 < p$, then $s<p-1$, and by the definition of $s$, the arc $[\ml_s, \ml_p]$ belongs to $\bT$ and crosses $u$, a contradiction. If $t < s_1 < s$, then $t<s-1$, and by the definition of $t$, the arc $[\ml_t, \ml_s]$ belongs to $\bT$ and crosses $u$, a contradiction. Thus, $s_1 = s$, which is a contradiction to the definition of $t$. Therefore, $v$ does not cross any arc of $\bT$, and hence, $v\in\bT$, a contradiction to the minimality of $s$. This establishes our claim. In particular, $\ml_p$ is a left $\bT$-fountain whenever it is right $\bT$-bounded.

Assume now that $\ml_p$ is not right $\bT$-bounded. We shall consider only the case where $[\ml_p, \mr_i]\notin \bT$ for all but finitely many $i < -p$. In particular, there exists a minimal integer $q$ such that $[\ml_p, \mr_q]\in \bT$. Since $\ml_p$ is not right $\bT$-bounded, $[\ml_p, \ml_j] \in \bT$ for infinitely many $j > p$. By Lemma \ref{nca}(2), no $\ml_j$ with $j>p$ is an endpoint of an arc of $C(\bT)$. Thus, $[\ml_p, \mr_q]$ is a maximal element in $\C(\bT)$, a contradiction. This shows that $[\ml_p, \mr_i] \in \bT$ for infinitely many $i < -p$. Using a similar argument as above, we see that $\ml_p$ is right $\bT$-unbounded. That is, $\ml_p$ is a full $\bT$-fountain. The proof of the lemma is completed.

\medskip

We are ready to obtain our main result of this section, which gives an easy criterion for a triangulation of $\B_\infty$ to be compact.

\medskip

\begin{Theo} \label{TheoClusterTilting}

If $\bT$ is a triangulation of $\B_\infty$, then it is compact if and only if it contains some connecting arcs such that every marked point in $\B_\infty$ is either $\bT$-bounded or an endpoint of infinitely many connecting arcs of $\bT$.

\end{Theo}

\noindent{\it Proof.} By Lemmas \ref{Comapct-2} and \ref{DI-chain}, we shall need only to prove the sufficiency. For this purpose, assume that $\bT$ is a triangulation of $\B_\infty$ such that $\C(\bT)$ is non-empty and every marked point in $\B_\infty$ is either $\bT$-bounded or an endpoint of infinitely many arcs of $C(\bT)$.

Fix an arc $v$ in $\mathcal{B}_\infty$. We shall need to prove that $\bT_v$ is compact. By Lemma \ref{Cofinite-compact}, we may assume that $\bT_v$ is infinite. Then, by Lemma \ref{inf-crossing}, some marked point is an endpoint of infinitely many arcs of $\bT_v$; and by Lemma \ref{LemmaFountain}, such a marked point is a $\bT$-fountain base. Denoting by $t$ the number of such $\bT$-fountain bases, we have $t\le 2$ by Lemma \ref{Full-fountain}. Let $\mf{p}_i$, with $i\in \{1, t\}$, be the $\bT$-fountain bases such that $\bT_v\cap \mathbb{F}_{\hspace{-1.5pt}_{\bT}}(\mf{p_i})$ is infinite. By Lemma \ref{LemmaCrossFountain}, we may assume that each $\mf{p}_i$ with $i\in \{1, t\}$ is a left or right $\bT$-fountain base. By Lemma \ref{LemmaCrossFountain2}, each $\bT_v\cap \mathbb{F}_{\hspace{-1.5pt}_{\bT}}(\mf{p}_i)$ with $i\in \{1, t\}$ is compact and co-finite in $\mathbb{E}_{\hspace{-0.5pt}_{\bT}}(\mf{p}_i)$. It is then easy to see that $\cup_{1\le i\le t}\, \bT_v\cap \mathbb{F}_{\hspace{-1.5pt}_{\bT}}(\mf{p}_i)\vspace{1pt}$ is compact. We claim that $\cup_{1\le i\le t}\, \bT_v\cap \mathbb{F}_{\hspace{-1.5pt}_{\bT}}(\mf{p}_i)\vspace{1pt}$ is co-finite in $\bT_v$. Indeed, given any marked point $\mf{q}$ in $\B_\infty$, we set
$$\Oa(\mf{q})=\left\{\begin{array}{ll}
\hspace{-2pt} \mathbb{E}_{\hspace{-0.5pt}_{\bT}}(\mf{q}) \backslash  \left(\bT_v\cap \mathbb{F}_{\hspace{-1.5pt}_{\bT}}(\mf{q})\right),
 \hspace{-2pt} & \mbox{ if }  \; \mf{q} \in \{\mf{p}_1, \mf{p}_t\},\\ \vspace{-10pt} \\
\hspace{-2pt} \bT_v\cap \mathbb{E}_{\hspace{-0.5pt}_{\bT}}(\mf{q}), \hspace{-2pt} & \mbox{ if }  \; \mf{q} \not\in \{\mf{p}_1, \mf{p}_t\},
\end{array} \right.$$
which is always finite by Lemma \ref{LemmaCrossFountain2} and the assumption.

Suppose first that $v$ is an upper arc, say $v=[\ml_r, \ml_s]$ with $r<s-1$. Let $u$ be an arc lying in $\bT_v$ but not in $\mathbb{F}_{\hspace{-1.5pt}_{\bT}}(\mf{p}_1)\cup \mathbb{F}_{\hspace{-1.5pt}_{\bT}}(\mf{p}_t)\vspace{1pt}$. Since $u$ crosses $[\ml_r, \ml_s]$, there exists some
$i$ with $r<i<s$ such that $u\in \mathbb{E}_{\hspace{-0.5pt}_{\bT}}(\ml_i),$ and by definition, $u\in \Oa(\ml_i)$. That is, $u\in \cup_{r<i<s}\, \Oa(\ml_i)$. Thus, our claim holds in this case. Similarly, we can establish the claim in case $v$ is a lower arc.

Next, suppose $v$ is a connecting arc, say $v=[\ml_r, \mr_s]$. We shall consider only the case where $\mf{p}_1$ is an upper marked point and a left $\bT$-fountain base. It is easy to see that $\mf{p}_1=\ml_{p_1}$ for some $p_1>r$. Then, $\mathbb{F}_{\hspace{-1.5pt}_{\bT}}(\mf{p}_1)$ contains some connecting arc $w=[\ml_{p_1}, \mr_q]$ with $q>s$.
Let $u$ be an arc lying in $\bT_v$ but not in $\cup_{1\le i\le t}\, \mathbb{F}_{\hspace{-1.5pt}_{\bT}}(\mf{p}_i).$ If $u$ is an upper arc then, since it does not cross $w$, we obtain $u=[\ml_i, \ml_j]$ with $i<r <j\le p_1$; and by definition, $u\in \Oa(\ml_j)$ for some $r<j\le p_1$. If $u$ is a connecting arc, since it does not cross $w$, we deduce from Lemma \ref{nca}(1) that $u=[\ml_i, \mr_j]$ with $i\ge p_1$ and $q\ge j >s$; and by definition, $u\in \Oa(\mr_j)$ for some $q\ge j>s$. If $u$ is a lower arc, since $u$ does not cross $w$, we obtain $u=[\mr_j, \mr_i]$ with $q\ge j >s$; and by definition, $u\in \Oa(\mr_j)$ for some $q\ge j>s$. In any case, $u$ belongs to the finite union of the $\Oa(\ml_i)$ with $r<i\le p_1$ and the $\Oa(\mr_j)$ with $q\ge j>s$. That is, our claim holds in this case. By Lemma \ref{Cofinite-compact}, $\bT_v$ is compact. The proof of the theorem is completed.

\medskip

Here is an example of a compact triangulation of $\B_\infty$ having two fountains.

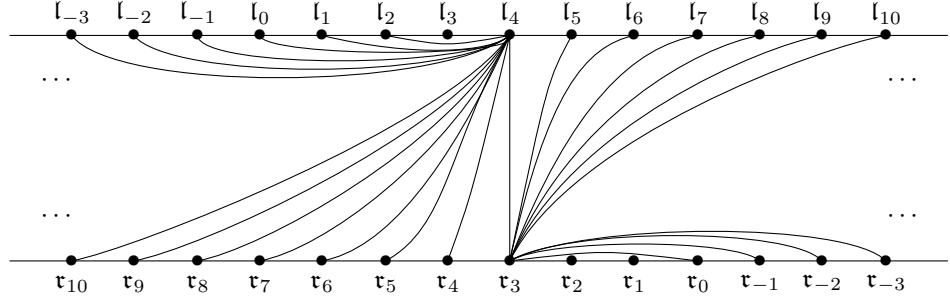
\begin{figure}[h]
  \centering
  \begin{tikzpicture}[xscale=2.50,yscale=1.5]

    \path (-2.4,0.6) node{$\cdots$};
    \path (2.1,0.6) node{$\cdots$};

   \path (-2.4,-0.6) node{$\cdots$};
    \path (2.1,-0.6) node{$\cdots$};

    \draw (-2.66,1) -- (2.33,1);
    \draw (-2.66,-1) -- (2.33,-1);

    \node at (-2.33,1.2) {$\ml_{-3}$};
    \node at (-2.0,1.2) {$\ml_{-2}$};
    \node at (-1.66,1.2) {$\ml_{-1}$};
    \node at (-1.33,1.2) {$\ml_0$};
    \node at (-1,1.2) {$\ml_1$};
 \node at (-0.66,1.2) {$\ml_{2}$};
    \node at (-0.33,1.2) {$\ml_{3}$};
    \node at (0,1.2) {$\ml_{4}$};
    \node at (0.33,1.2) {$\ml_5$};
    \node at (0.66,1.2) {$\ml_6$};
 \node at (1,1.2) {$\ml_{7}$};
    \node at (1.33,1.2) {$\ml_{8}$};
    \node at (1.66,1.2) {$\ml_{9}$};
    \node at (2,1.2) {$\ml_{10}$};

    \node at (-2.33,-1.2) {$\mr_{10}$};
    \node at (-2.0,-1.2) {$\mr_{9}$};
    \node at (-1.66,-1.2) {$\mr_{8}$};
    \node at (-1.33,-1.2) {$\mr_7$};
    \node at (-1,-1.2) {$\mr_6$};
 \node at (-0.66,-1.2) {$\mr_{5}$};
    \node at (-0.33,-1.2) {$\mr_{4}$};
    \node at (0,-1.2) {$\mr_{3}$};
    \node at (0.33,-1.2) {$\mr_2$};
    \node at (0.66,-1.2) {$\mr_1$};
 \node at (1,-1.2) {$\mr_{0}$};
    \node at (1.33,-1.2) {$\mr_{-1}$};
    \node at (1.66,-1.2) {$\mr_{-2}$};
    \node at (2,-1.2) {$\mr_{-3}$};

    \node at (-2.33,1.00) {$\bullet$};
    \node at  (-2.00,1.00) {$\bullet$};
    \node at  (-1.66,1.0){$\bullet$};
    \node at  (-1.33,1.0){$\bullet$};
    \node at  (-1.00,1.0){$\bullet$};
    \node at  (-0.66,1.0){$\bullet$};
    \node at  (-0.33,1.0){$\bullet$};
    \node at  (0.00,1.0){$\bullet$};
    \node at  (0.33,1.0){$\bullet$};
    \node at  (0.66,1.0){$\bullet$};
    \node at  (1.00,1.0){$\bullet$};
    \node at  (1.33,1.0){$\bullet$};
    \node at  (1.66,1.0){$\bullet$};
    \node at  (2.00,1.0){$\bullet$};

    \node at (-2.33,-1.0){$\bullet$};
    \node at (-2.00,-1.0){$\bullet$};
    \node at (-1.66,-1.0){$\bullet$};
   \node at (-1.33,-1.0){$\bullet$};
   \node at (-1.00,-1.0){$\bullet$};
    \node at (-0.66,-1.0){$\bullet$};
    \node at (-0.33,-1.0){$\bullet$};
    \node at (0.00,-1.0){$\bullet$};
    \node at (0.33,-1.0){$\bullet$};
    \node at(0.66,-1.0){$\bullet$};
    \node at (1.00,-1.0){$\bullet$};
   \node at (1.33,-1.0){$\bullet$};
    \node at(1.66,-1.0){$\bullet$};
    \node at (2.00,-1.0){$\bullet$};

\draw (-2.33, -1) .. controls (-2, -0.9) and (-0.33, 0.2) .. (0, 1);
\draw (-2, -1) .. controls (-1.66, -0.9) and (-0.33, 0.1) .. (0, 1);
\draw (-1.66, -1) .. controls (-1.33, -0.9) and (-0.33, 0) .. (0, 1);
\draw (-1.33, -1) .. controls (-1, -0.9) and (-0.33, -0.1) ..  (0, 1);
\draw (-1, -1) .. controls (-0.66, -0.9) and (-0.33, -0.2) ..  (0, 1);
\draw (-0.66, -1) .. controls (-0.33, -0.7) and (-0.33, 0) .. (0, 1);
\draw (-0.33, -1) .. controls (-0.2, -0.4) .. (0, 1);
\draw (0, -1) -- (0, 1);

\draw (0, -1) .. controls (0.2, 0.6) .. (0.33, 1);
\draw (0, -1) .. controls (0.2, 0.5) and (0.4, 0.9) .. (0.66, 1);
\draw (0, -1) .. controls (0.33, 0.4) and (0.66, 0.9) .. (1, 1);
\draw (0, -1) .. controls (0.33, 0.3) and (1, 0.9) .. (1.33, 1);
\draw (0, -1) .. controls (0.33, 0.2) and (1.33, 0.9) .. (1.66, 1);
\draw (0, -1) .. controls (0.33, 0.1) and (1.66, 0.9) .. (2, 1);

\draw (-0.66, 1) .. controls (-0.33, 0.9) .. (0,1);
\draw (-1, 1) .. controls (-0.5, 0.8) and (-0.1, 0.85).. (0,1);
\draw (-1.33, 1) .. controls (-1.30, 0.8) and (0.1, 0.8) .. (0,1);
\draw (-1.66, 1) .. controls (-1.65, 0.7) and (-0.2, 0.7) .. (0,1);
\draw (-2, 1) .. controls (-1.8, 0.6) and (-0.3, 0.6) .. (0,1);
\draw (-2.33, 1) .. controls (-2.1, 0.5) and (-0.4, 0.5) .. (0,1);

\draw (0, -1) .. controls (0.5, -0.9) .. (1,-1);
\draw (0, -1) .. controls (0.2, -0.8) and (1.2, -0.8) .. (1.33,-1);
\draw (0, -1) .. controls (0.3, -0.7) and (1.5, -0.7) .. (1.66,-1);
\draw (0, -1) .. controls (0.4, -0.65) and (1.8, -0.65) .. (2.0,-1);

  \end{tikzpicture}
\label{fig:triangulation3}
\end{figure}

\section{Geometric Realization of cluster categories of type $\mathbb{A}_\infty^\infty$}

The objective of this section is to study the cluster structure of a cluster category of type $\mathbb{A}_\infty^\infty$ in terms of triangulations of the infinite strip $\B_\infty$, as studied in Section 5. We start this section with some algebraic considerations. Throughout, $Q$ stands for a canonical quiver with no infinite path of type $\mathbb{A}_\infty^\infty$, that is, its vertices are the integers and its arrows are of form $n\to (n+1)$ or $(n+1)\to n$. Let $a_i, b_i$, $i\in \Z$, be the sources and the sinks in $Q$ respectively such that $b_{i-1}<a_i<b_i$. We denote by $p_i: a_i \rightsquigarrow b_i$ and $q_i: a_i\rightsquigarrow b_{i-1},$ $i\in \Z$, the maximal paths in $Q$. It will be convenient to picture $Q$ as follows:

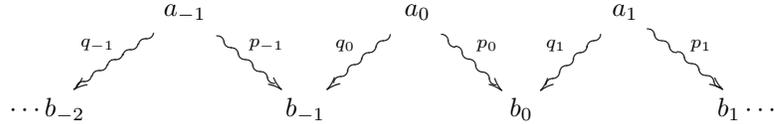
\vspace{-5pt}

$$\xymatrix{
&&a_{-1}\ar@{~>}[dl]_{q_{-1}}\ar@{~>}[dr]^{p_{-1}}&&a_0\ar@{~>}[dr]^{p_0}\ar@{~>}[dl]_{q_0} &&a_1\ar@{~>}[dl]_{q_1}\ar@{~>}[dr]^{p_1}\\
&\cdots b_{-2}&&b_{-1}&&b_0&&b_1\cdots }
\vspace{3pt}$$

Let $\mathscr{S}$ be a set of paths having pairwise distinct starting points. For $p, q\in \mathscr{S}$, define an order $\preceq$ on $\mathscr{S}$ so that $p\preceq q$ if and only if $s(p)\le s(q)$. In case $p\prec q$ and there exists no $u$ in $\mathscr{S}$ such that $p\prec u\prec q,$ we write $p=\sigma_{_{\hspace{-1.5pt}\mathscr{S}}}(q)$ and $q=\sigma_{_{\hspace{-1.5pt}\mathscr{S}}}^-(p).$ This yields an injective map $\sigma_{_{\hspace{-1.5pt}\mathscr{S}}}: \mathscr{S}\to \mathscr{S},$ called {\it source translation} for $\mathscr{S}$.

\medskip

Let $Q_R$ stand for the union of the maximal paths $p_i$ with $i\in \Z$ and the trivial paths $\varepsilon_a$ with $a$ being a middle point of some $q_j$ with $j\in \Z.$ Dually, $Q_L$ stands for union of the maximal paths $q_i$ with $i\in \Z$ and the trivial paths $\varepsilon_b$ with $b$ being a middle point of some $p_j$ with $j\in \Z$. Recall from Lemma \ref{A-di-qs} that the Auslander-Reiten quiver $\Ga_{\rep(Q)}$ of $\rep(Q)$ has two regular components $\mathcal{R}_R$ and $\mathcal{R}_L$ such that the quasi-simple objects in $\mathcal{R\hspace{-1pt}}_R$ are the string representations $M(p)$ with $p \in Q_R$, while those in $\mathcal{R\hspace{-1pt}}_L$ are the string representations $M(q)$ with $q \in Q_L$. For convenience, we quote the following result stated in \cite[(5.13),(5.14)]{BLP}.

\medskip

\begin{Lemma}\label{AR-translation}

Let $\sigma_{\hspace{-1.5pt}_R}$ and $\sigma_{\hspace{-1.5pt}_L}$ be the source translations for $Q_R$ and $Q_L$, respectively.

\vspace{-1pt}

\begin{enumerate}[$(1)$]

\item If $p\in Q_R$, then $\tau_{\hspace{-1pt}_Q}M(p) = M(\sigma_{\hspace{-1.5pt}_R}(p)).$

\vspace{1pt}

\item If $q\in Q_L$, then $\tau_{\hspace{-1pt}_Q} M(q) = M(\sigma_{\hspace{-1.5pt}_L}^-(q))$.

\end{enumerate}

\end{Lemma}

\medskip

We shall choose the vertex set of the Auslander-Reiten quiver $\Ga_{\hspace{-1pt}_{\mathscr{C}(Q)}}$ of $\mathscr{C}(Q)$ to be the fundamental domain $\mathscr{F}(Q)$ of $\mathscr{C}(Q)$, which consists of the regular representations in $\Ga_{\rep(Q)}$ and the objects in the connecting component of $\Ga_{D^b(\rep(Q))}.\vspace{1pt}$ By Proposition \ref{idt-clust-cpt}, $\Ga_{\hspace{-1pt}_{\mathscr{C}(Q)}}$ has exactly three connected components, namely, the connecting component $\C$ of shape $\Z\mathbb{A}_\infty^\infty$, and the two orthogonal regular components $\mathcal{R}_R$ and $\mathcal{R}_L$ of shape $\Z \mathbb{A}_\infty$. As stated in Proposition \ref{idt-clust-cpt}(3), the two regular components $\mathcal{R}_R$ and $\mathcal{R}_L$ are orthogonal. We shall describe morphisms from an object in $\mathcal{C}$ to an object in $\mathcal{R}_R$ or $\mathcal{R}_L$. For this purpose, we shall need some notation for $\mathcal{C}$. First, observe that $\mathcal{C}$ has a section; see \cite[2.1]{Liu1}, which is generated by the projective representations in $\Ga_{\rep(Q)}$ as follows:

\vspace{-10pt}

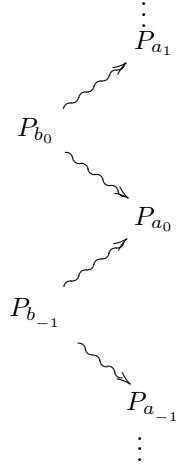
$$
\hspace{-15pt}\xymatrixcolsep{12pt}\xymatrixrowsep{12pt}\xymatrix{
&&&&&&&                \;\vdots &&&&&&&\\
&&&&&                      &&&&&&&
}\vspace{-35pt}$$
$$\xymatrixcolsep{18pt}\xymatrixrowsep{18pt}\xymatrix{
&&&&&                                 & P_{a_1} &&&&&&&\\
&&&&&  P_{b_0}\ar@{~>}[dr]\ar@{~>}[ur] & &&&&&&&&&&&&&&&\\
&&&&&                                  & P_{a_0} &&&&&&&&&&&&&&\\
&&&&& P_{b_{_{-1}}}\ar@{~>}[dr]\ar@{~>}[ur] &         &&&&&&&&&&&&&&\\
&&&&&                                  & P_{a_{_{-1}}} &&&&&&&&&&&&&&
} \vspace{-13pt}$$
$$
\xymatrixcolsep{18pt}\xymatrixrowsep{18pt}\xymatrix{
&&&&&                 & &\vdots&&&&&&&&&&&&& \\
&&&&&&                                 &  &&&&&&&&&&&&&&
}$$

\vspace{-20pt}

We shall denote by $R_0$ the unique double infinite sectional path in $\C$ containing the path $P_{b_0} \rightsquigarrow P_{a_0},$ which corresponds to the path $p_0$ in $Q;$ and by $L_0$ the unique double infinite sectional path containing the path $P_{b_{-1}} \rightsquigarrow P_{a_0},$ which corresponds to the path $q_0$ in $Q$. For each $i \in \Z$, put $R_i=\tau_{\hspace{-1pt}_\mathscr{C}}^iR_0\vspace{0.5pt}$ and $L_i=\tau_{\hspace{-1pt}_\mathscr{C}}^iL_0$. Observe that each object in $\mathcal{C}$ lies in a unique $R_i$ with $i\in \Z$ and in a unique $L_j$ with $j\in \Z$.

\medskip

\begin{Prop} \label{LemmaRays} Let $M$ be an object in $\C$. For each integer $i$, we have

\vspace{-2pt}

\begin{enumerate}[$(1)$]

\item $M \in R_i$ if and only if $\Hom_{\mathscr{C}(Q)}(M, \tau_{\hspace{-1pt}_\mathscr{C}}^iM(p_0)) \ne 0;$ and in this case, for any $Y\in \mathcal{R}_R$, $\Hom_{\,\mathscr{C}(Q)}(M,Y) \ne 0$ if and only if $Y \in \mathcal{W}(\tau_{\hspace{-1pt}_\mathscr{C}}^iM(p_0));$

\item $M \in L_i$ if and only if $\Hom_{\mathscr{C}(Q)}(M, \tau_{\hspace{-1pt}_\mathscr{C}}^iM(q_0)) \ne 0;$ and in this case, for any $Y\in \mathcal{R}_L$, $\Hom_{\,\mathscr{C}(Q)}(M,Y) \ne 0$ if and only if $Y \in \mathcal{W}(\tau_{\hspace{-1pt}_\mathscr{C}}^iM(q_0))$.

\end{enumerate}

\end{Prop}

\noindent{\it Proof.} We shall only prove (1) since the proof of (2) is similar. For each $i\in \Z$, write the path $q_i$ as follows:
$$\xymatrixcolsep{18pt}\xymatrix{b_{i-1}=a_{i,\ell_i} & a_{i,\ell_i-1}\ar[l] & \cdots \ar[l] & a_{i,1}\ar[l] & \ar[l] a_{i,0}=a_i,
}$$
where $\ell_i$ is the length of $q_i$. Write $\sigma=\sigma_{\hspace{-1.5pt}_R},$ the source translation for $Q_R$. By definition, $\sigma^j(p_i)=\varepsilon_{a_{i, j}}$ for $0<j<\ell_i$, and $\sigma^{\ell_i}(p_i)=p_{i-1}$. Setting $t_0 = 0$ and $t_i =t_{i+1}+\ell_{i+1}=\ell_0+\cdots +\ell_{i+1}$ for $i \le -1$, we get a sequence of non-negative integers:
$$0=t_0< t_{-1} < \cdots < t_i < t_{i-1}< \cdots $$
 For each $j\ge 0,$ there exists a unique $i\le 0$ such that $t_i\le j < t_{i-1}$. In this case, it is easy to see that
$$\sigma^j(p_0)=\left\{
  \begin{array}{ll}
    p_i, & \hbox{if $j = t_i\hspace{0.4pt};$} \\
    \varepsilon_{a_{i,  j - t_i}}, & \hbox{if $t_i < j < t_{i-1}$.}
  \end{array}
\right.$$

Consider the section in $\mathcal{C}$ generated by the representations $P_a$ with $a\in Q_0$. Recall that $R_0$ is the double infinite sectional path in $\mathcal{C}$ containing $P_{b_0}\rightsquigarrow P_{a_0}$. For each $x\in Q_0$ with $x\le 0$, there exists in $\mathcal{C}$ a unique sectional path $u_x: P_x\rightsquigarrow M_x$ with $M_x\in R_0$, whose length is denoted by $d_x$. On the other hand, there exists a unique path $\rho_x\in Q_R$ with $\rho_x\preceq p_0$ such that $x\in \rho_x$. We claim that $\rho_x=\sigma^{d_x}(p_0)$.
Indeed, considering the rectangle in $\mathcal{C}$ with vertices $P_{b_i}, M_{b_i}, P_{a_i}$ and $M_{a_i}$ as follows:
$$\xymatrixcolsep{18pt}\xymatrixrowsep{18pt}
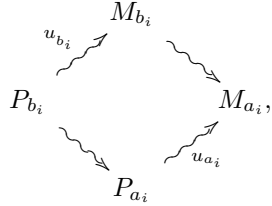
\xymatrix{
&&                                            & M_{b_i}\ar@{~>}[dr]           &        \\
&&P_{\hspace{0.4pt}b_i}\ar@{~>}[ur]^{u_{_{b_i}}}\ar@{~>}[dr] &                                 &M_{a_i},\\
&&                                            & P_{a_i}\ar@{~>}[ur]_{u_{a_i}} &        }$$
we see that $d_{b_i}=d_{a_i}$, for any $i\ge 0$. Now, for each $x\in Q_0$ with $x\le 0$, there exists a unique integer $r\le 0$ such that $p_{r-1}\prec \rho_x\preceq p_r$. If $\rho_x=p_r$, then $x\in p_r$ and hence, $P_x$ lies on the sectional path $P_{b_r}\rightsquigarrow P_{a_r}$. Considering the rectangle in $\mathcal{C}$ with vertices $P_{b_r}, M_{b_r}, P_x$ and $M_x$, we see that $d_x=d_{b_r}$. Now, the path $u_{b_r}: P_{b_r}\rightsquigarrow M_{b_r}$ is the composite of the path
$P_{b_r}\rightsquigarrow P_{a_{r+1}}$ corresponding to $q_{r+1}$ and the path $u_{a_{r+1}}: P_{a_{r+1}}\rightsquigarrow M_{a_{r+1}}$. Hence, $$d_x=d_{b_r}=d_{a_{r+1}}+\ell_{r+1}=d_{b_{r+1}}+\ell_{r+1}=\ell_0+\cdots +\ell_{r+1}=t_r.$$
This yields $\rho_x=p_r=\sigma^{t_r}(p_0)=\sigma^{d_x}(p_0)$. Suppose now that $p_{r-1}\prec \rho_x\prec p_r.$ Then,
$x$ lies in the path $q_r: a_r\rightsquigarrow b_{r-1},$ that is, $x=a_{r,l}$ for some $0<l<\ell_r$. In this case,
$\rho_x=\varepsilon_{a_{r, l}}=\sigma^{t_r+l}(p_0)$. On the other hand, the path
$u_x: P_x\rightsquigarrow M_x$ is the composite of the path
$$\xymatrix{P_x=P_{a_{r, l}} \ar[r]& \cdots \ar[r] & P_{a_{r,0}}=P_{a_r},}$$ corresponding to $a_{r, l}\leftarrow \cdots \leftarrow a_{r,0}=a_r$, and the path $u_{a_r}: P_{a_r}\rightsquigarrow M_{a_r}$. This yields $d_x=l+d_{a_r}=l+d_{b_r}=l+t_r$. Since $t_r<t_r+l<t_{r-1}$, we obtain $\sigma^{d_x}(p_0)=\sigma^{l+t_r}(p_0)=\varepsilon_{a_{r, l}}=\rho_x$. This establishes our claim.

Let $M\in R_0$, which we assume to be a successor of $P_{b_0}$ in $R_0$. Then, $M=\tau_{\hspace{-2pt}_\mathscr{C}}^{-s}\hspace{-1.5pt}P_x$ for a unique pair $(s, x)$, where $s$ is a non-negative integer and $x\in Q_0$ with $x\le 0.$ As shown above, $x\in \rho_x=\sigma^{d_x}(p_0)$, where $d_x$ is the length of the sectional path $u_x: P_x\rightsquigarrow M_x$ in $\mathcal{C}$, where $M_x\in R_0$. Since $R_0$ has a subpath $M_x\rightsquigarrow M$ of length $d_x$, we see that $M=\tau_{\hspace{-2pt}_\mathscr{C}}^{-d_x}\hspace{-1.5pt}P_x$. Since $\tau_{\hspace{-2pt}_\mathscr{C}}$ is an auto-equivalence of $\mathscr{C}(Q)$, we obtain
$$\begin{array}{rcl}
\Hom_{\mathscr{C}(Q)}(M, M(p_0)) & \cong & \Hom_{\mathscr{C}(Q)}(P_x, \tau_{\hspace{-2pt}_\mathscr{C}}^{d_x}\hspace{-2pt}M(p_0))\\ \vspace{-10pt}\\
& \cong & \Hom_{\mathscr{C}(Q)}(P_x, \tau_{\hspace{-2pt}_Q}^{d_x}\hspace{-2pt}M(p_0))\\ \vspace{-10pt}\\
&=& \Hom_{\mathscr{C}(Q)}(P_x, M(\sigma^{d_x}\hspace{-2pt}(p_0)))\\ \vspace{-10pt}\\
&=&\Hom_{\mathscr{C}(Q)}(P_x, M(\rho_x))\\
&\cong &\Hom_{\,\rep(Q)}(P_x, M(\rho_x))\\
&\ne & 0,\end{array}$$
where the first equality follows from Lemma \ref{AR-translation} and the last isomorphism follows from Proposition \ref{special-maps}. In case $M$ is a predecessor of $P_{b_0}$ in $R_0$, we may show by a dual argument that $\Hom_{\mathscr{C}(Q)}(M, M(p_0)) \ne 0.$ That is, $\Hom_{\mathscr{C}(Q)}(M, M(p_0)) \ne 0$ for all $M\in R_0$. As a consequence, if $M\in R_i$ with $i\in \Z$, say $M=\tau_{\hspace{-2pt}_\mathscr{C}}^iU$ with $U\in R_0$, then $\Hom_{\mathscr{C}(Q)}(M, \tau_{\hspace{-2pt}_\mathscr{C}}^iM(p_0)) \cong \Hom_{\mathscr{C}(Q)}(U, M(p_0))\ne 0.$

\vspace{2pt}

Conversely, let $M\in R_i$ with $i\in \Z$, say $M=\tau_{\hspace{-2pt}_\mathscr{C}}^iU$ for some $U\in R_0$, be such that $\Hom_{\mathscr{C}(Q)}(M, M(p_0)) \ne 0$. Since $\tau_{\hspace{-2pt}_Q}^{-i}M(p_0)=\tau_{\hspace{-2pt}_\mathscr{C}}^{-i}M(p_0),$  by Lemma \ref{AR-translation},
we have \vspace{-2pt} $$\Hom_{\mathscr{C}(Q)}(U, M(\sigma^{-i}(p_0)))\cong \Hom_{\mathscr{C}(Q)}(U, \tau_{\hspace{-2pt}_Q}^{-i}M(p_0))\cong \Hom_{\mathscr{C}(Q)}(M, M(p_0))\ne 0.$$
On the other hand, as we have shown, $\Hom_{\mathscr{C}(Q)}(U, M(p_0))\ne 0$. If $U$ is a preprojective representation then, by Proposition \ref{wing}, $M(\sigma^{-i}(p_0))\in \mathcal{W}(M(p_0))$, and hence $i=0$. Otherwise, $U=\tau_{\hspace{-2pt}_\mathscr{C}}^jP_x$ for some $j>0$ and $x\in Q_0$. In this case, we obtain
$$\begin{array}{rcl}
\Hom_{\mathscr{C}(Q)}(P_x, M(\sigma^{-j-i}(p_0))) &\cong & \Hom_{\mathscr{C}(Q)}(P_x, \tau_{\hspace{-2pt}_Q}^{-j}M(\sigma^{-i}(p_0)))\\ \vspace{-10pt} \\
&\cong & \Hom_{\mathscr{C}(Q)}(\tau_{\hspace{-2pt}_\mathscr{C}}^jP_x, M(\sigma^{-i}(p_0)))\\ \vspace{-10pt} \\
&=& \Hom_{\mathscr{C}(Q)}(U, M(\sigma^{-i}(p_0)))\\
&\ne & 0.\end{array}$$
Similarly, $\Hom_{\mathscr{C}(Q)}(P_x, M(\sigma^{-j}(p_0)))\cong \Hom_{\mathscr{C}(Q)}(U, M(p_0))\ne 0.$ By Proposition \ref{wing}, $M(\sigma^{-j-i}(p_0))\in \mathcal{W}(M(\sigma^{-j}(p_0))$, and hence $i=0$.

More generally, let $M\in R_j$ with $j\in \Z$, say $M=\tau_{\hspace{-2pt}_\mathscr{C}}^jU$ for some $U\in R_0$, be such that $\Hom_{\mathscr{C}(Q)}(M, \tau_{\hspace{-2pt}_\mathscr{C}}^iM(p_0))\ne 0.$ This yields
$$\Hom_{\mathscr{C}(Q)}(U, \tau_{\hspace{-2pt}_\mathscr{C}}^{i-j}M(p_0))\cong \Hom_{\mathscr{C}(Q)}(M, \tau_{\hspace{-2pt}_\mathscr{C}}^iM(p_0))\ne 0.$$
As we have just shown,  $\tau_{\hspace{-2pt}_\mathscr{C}}^{i-j}M(p_0)\in \mathcal{W}(M(p_0))$, and hence $j=i$.

\vspace{1pt}

Finally, let $M\in R_i$ with $i\in \Z$ and $Y\in \mathcal{R}_R$. Then, $\Hom_{\mathscr{C}(Q)}(M, \tau_{\hspace{-2pt}_\mathscr{C}}^iM(p_0)) \ne 0.\vspace{1pt}$ If $M$ is a representation, then $\Hom_{\mathscr{C}(Q)}(M,Y) \cong \Hom_{\rep(Q)}(M,Y)\vspace{1pt}$ by Proposition \ref{special-maps}. It follows from Proposition \ref{wing} that $\Hom_{\mathscr{C}(Q)}(M,Y) \ne 0,$ if and only if, $Y\in \mathcal{W}(\tau_{\hspace{-2pt}_\mathscr{C}}^iM(p_0))$.
Otherwise, $M=\tau_{\hspace{-2pt}_\mathscr{C}}^jP_y\vspace{1pt}$ for some $j>0$ and $y\in Q_0$. In particular, $P_y\in R_{i-j}$, and hence $\Hom_{\mathscr{C}(Q)}(P_y, \tau_{\hspace{-2pt}_\mathscr{C}}^{i-j}M(p_0)) \ne 0.$ On the other hand,
$$\Hom_{\mathscr{C}(Q)}(M,Y)\cong \Hom_{\mathscr{C}(Q)}(P_y, \tau_{\hspace{-2pt}_\mathscr{C}}^{-j}Y) \cong
\Hom_{\rep(Q)}(P_y, \tau_{\hspace{-2pt}_\mathscr{C}}^{-j}Y).$$
By Proposition \ref{wing}, $\Hom_{\mathscr{C}(Q)}(M,Y)\ne 0$ if and only if $\tau_{\hspace{-2pt}_\mathscr{C}}^{-j}Y\in \mathcal{W}(\tau_{\hspace{-2pt}_\mathscr{C}}^{i-j}M(p_0))$. The latter is equivalent to $Y\in \mathcal{W}(\tau_{\hspace{-2pt}_\mathscr{C}}^iM(p_0))$. The proof of the proposition is completed.

\medskip

We shall parameterize the indecomposable objects of $\mathscr{C}(Q)$ by the arcs in $\B_\infty$, that is, to define a bijection $$\varphi: \mathscr{F}(Q) \to {\rm arc}(\B_\infty),$$ where $\mathscr{F}(Q)$ is the fundamental domain of $\mathscr{C}(Q)$. Recall that $\mathscr{F}(Q)$ consists of the objects in $\mathcal{C}$, $\mathcal{R}_R$ and $\mathcal{R}_L$. First, for each object $M$ in $\C$, there exists a unique pair of integers $(i,j)$ such that $M=L_i\cap R_j$. Sending $M$ to the connecting arc $[\ml_i, \mr_j]$ in $\B_\infty$ yields a bijection $\varphi_{_\C}$ from the objects in $\mathcal{C}$ onto the connecting arcs in $\B_\infty$. Put $\varphi(M)=\varphi_{_\C}(M),$ for each object in $\mathcal{C}.$

\medskip

Next, consider $S_L=\tau_{\hspace{-1.5pt}_\mathscr{C}}^-M(q_0),\vspace{1pt}$ a quasi-simple object in $\mathcal{R}_L$. For each $i\in \Z$, denote by $L_i^+$ the ray in $\mathcal{R}_L$ starting with $\tau_{\hspace{-1.5pt}_\mathscr{C}}^iS_L$, and by $L_i^-$ the coray ending with $\tau_{\hspace{-1.5pt}_\mathscr{C}}^iS_L.\vspace{1pt}$ For each object $X$ in $\mathcal{R}_L$, there exists a unique pair of integers $(i, j)$ with $i\le j$ such that $X=L_i^- \cap L_j^+$, and we set $\varphi_{_L}(X)= [\ml_{i-1}, \ml_{j+1}]\in {\rm arc}(\B_\infty)$. This defines a bijection $\varphi_{_L}$ from the objects in $\mathcal{R}_L$ onto the upper arcs in $\B_\infty$. Put
$\varphi(X)= \varphi_{_L}(X),$ for all $X\in \mathcal{R}_L.$ In this way, the quasi-simple objects in $\mathcal{R}_L$ are those mapped by $\varphi$ to $[\ml_i, \ml_j]$ with $|i-j|=2$.

\medskip

Finally, consider $S_R=\tau_{\hspace{-1.5pt}_\mathscr{C}}^-M(p_0),\vspace{1pt}$ a quasi-simple object in $\mathcal{R}_R$. For $i\in \Z$, denote by $R_i^+$ the ray in $\mathcal{R}_R$ starting with $\tau_{\hspace{-1.5pt}_\mathscr{C}}^iS_R;\vspace{1pt}$ and by $R_i^-$ the coray ending with $\tau_{\hspace{-1.5pt}_\mathscr{C}}^iS_R$. For each object $Y \in \mathcal{R}_R$, there exists a unique pair of integers $(i,j)$ with $i\ge j$ such that $Y=R_i^+ \cap R_j^-$, and we set $\varphi_{_R}(Y)=[\mr_{i+1}, \mr_{j-1}]\in {\rm arc}(\B_\infty)$. This yields a bijection $\varphi_{_R}$ from the objects in $\mathcal{R}_R$ onto the lower arcs $[\mr_i, \mr_j]$ in $\B_\infty$. Set $\varphi(Y)=\varphi_{_R}(Y),$ for all $Y\in \mathcal{R}_R.$ Observe that the quasi-simple objects in $\mathcal{R}_R$ are those mapped by $\varphi$ to $[\mr_i, \mr_j]$ with $|i-j|=2$.

\medskip

This yields the desired bijection $\varphi: \mathscr{F}(Q) \to {\rm arc}(\B_\infty).$ To simplify the notation, for each object $X$ in $\mathscr{F}(Q)$ and each arc $u$ in $\mathcal{B}_\infty$, we write \vspace{1pt}
$$a_{_X}=\varphi(X) \; \mbox{ and } \; M_u=\varphi^{-1}(u).$$

\medskip

{\sc Example.} In Figure \ref{fig:regular} below, the two black dots are objects in $\mathcal{R}_L$, mapped by $\varphi$ to $[\ml_2, \ml_8]$ and $[\ml_1, \ml_5]$, respectively. We see that the quasi-socle of $\varphi^{-1}[\ml_2, \ml_8]$ is $\varphi^{-1}[\ml_6, \ml_8]$ and its quasi-top is $\varphi^{-1}[\ml_2, \ml_4]$.

\vspace{3pt}

\begin{figure}[h]
\centering
\begin{tikzpicture}[xscale=3.30,yscale=1.5]

    \draw (-2.0,-1) -- (1.66,-1);

    \node at (-1.64,-1.15){}; 
    \node at (-1.31,-1.15) {$L^+_7$};
    \node at (-0.98,-1.15){};
 \node at (-0.64,-1.15){};
    \node at (-0.31,-1.15) {$L^+_4$};
    \node at (0.02,-1.15) {$L^-_3$};
    \node at (0.35,-1.15) {$L_2^-$};
    \node at (0.68,-1.15){};
 \node at (1.02,-1.15) {$S_L$};
    \node at (1.39,-1.15){}; 

\node at (-0.69, -0.58) {$\bullet$};
\node at (-0.00, -0.8) {$\bullet$};

    \draw (-2.00,-0.95) -- (-2.00,-1.0);
    \draw (-1.66,-0.95) -- (-1.66,-1.0);
    \draw (-1.33,-0.95) -- (-1.33,-1.0);
    \draw (-1.00,-0.95) -- (-1.00,-1.0);
    \draw (-0.66,-0.95) -- (-0.66,-1.0);
    \draw (-0.33,-0.95) -- (-0.33,-1.0);
    \draw (0.00,-0.95) -- (0.00,-1.0);
    \draw (0.33,-0.95) -- (0.33,-1.0);
    \draw (0.66,-0.95) -- (0.66,-1.0);
    \draw (1.00,-0.95) -- (1.00,-1.0);
    \draw (1.33,-0.95) -- (1.33,-1.0);
    \draw (1.66,-0.95) -- (1.66,-1.0);
    \draw (2.00,-0.95) -- (2.00,-1.0);

\draw (-2, 0.75) -- (-1.4, 1);
\draw (-2, 0.5) -- (-1, 1);
\draw (-2, 0.25) -- (-0.66, 1);
\draw (-2, 0) -- (-0.33, 1);
\draw (-2, -0.25) -- (0, 1);
\draw (-2, -0.5) -- (0.33, 1);
\draw (-2, -0.75) -- (0.66, 1);
\draw (-2,-1.0) -- (1.0,1);
\draw (-1.66,-1.0) -- (1.33,1);
\draw (-1.33,-1.0) -- (1.66,1);
\draw (-1.0,-1.0) -- (1.66,0.75);
\draw (-0.66,-1.0) -- (1.66,0.5);
\draw (-0.33,-1.0) -- (1.66,0.25);
\draw (0.0,-1.0) -- (1.66,0);
\draw (0.33,-1.0) -- (1.66,-0.25);
\draw (0.66,-1.0) -- (1.66,-0.5);
\draw (1,-1.0) -- (1.66,-0.75);
\draw (1.33,-1.0) -- (1.66,-1);

\draw (1.66,0.75) -- (1,1);
\draw (1.66,0.5) -- (0.66,1);
\draw (1.66,0.25) -- (0.33,1);
\draw (1.66,0) -- (0,1);
\draw (1.66,-0.25) -- (-0.33,1);
\draw (1.66,-0.5) -- (-0.66,1);
\draw (1.66,-0.75) -- (-1,1);
\draw (1.66,-1.0) -- (-1.33,1);
\draw (1.33,-1.0) -- (-1.66,1);
\draw (1,-1.0) -- (-2,1);
\draw (0.66,-1.0) -- (-2,0.75);
\draw (0.33,-1.0) -- (-2,0.5);
\draw (0,-1.0) -- (-2,0.25);
\draw (-0.33,-1.0) -- (-2,0);
\draw (-0.66,-1.0) -- (-2,-0.25);
\draw (-1,-1.0) -- (-2,-0.5);
\draw (-1.33,-1.0) -- (-2,-0.75);
\draw (-1.66,-1.0) -- (-2,-1);

\end{tikzpicture}
\caption{The regular component $\mathcal{R}_L$}.
\label{fig:regular}
\end{figure}
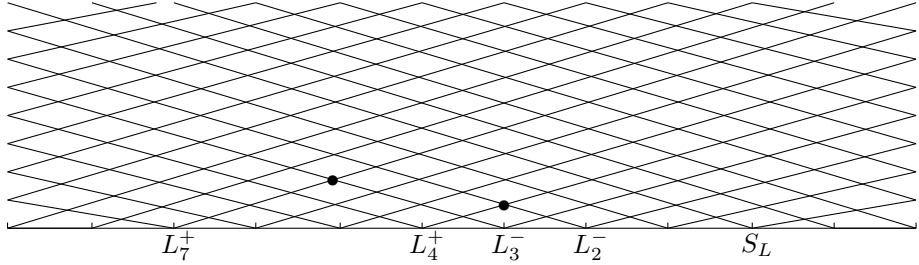


The following easy observation describes the Auslander-Reiten translation and the arrows in $\Ga_{\hspace{-1pt}_{\mathscr{C}(Q)}}$ in terms of the arcs in $\mathcal{B}_\infty$. Recall that ${\rm arc}(\B_\infty)$ is equipped with a translation $\tau$ defined in Definition \ref{arc-translation}.

\medskip

\begin{Lemma}\label{arrows}

Let $u, v$ be two distinct arcs in $\mathcal{B}_\infty$.

\begin{enumerate}[$(1)$]

\item In any case, $\tau_{\hspace{-0.5pt}_\mathscr{C}}M_u=M_{\tau u},\vspace{0.5pt}$ and $\tau_{\hspace{-0.5pt}_\mathscr{C}}^{-1}M_u=M_{\tau^{-1} u}\vspace{0.5pt}$.

\item If $u=[\ml_i, \mr_j]$, then
there exists an arrow $M_u\to M_v$ in $\Ga_{\hspace{-1pt}_{\mathscr{C}(Q)}}$ if and only if $v=[\ml_i, \mr_{j-1}]$ or $v=[\ml_{i-1}, \mr_j]$.

\vspace{1.5pt}

\item If $u=[\ml_i, \ml_j]$ with $i\le j-2$, then
there exists an arrow $M_u\to M_v$ in $\Ga_{\hspace{-1pt}_{\mathscr{C}(Q)}}$ if and only if $v=[\ml_i, \ml_{j-1}]$ with $i<j-2$ or $v=[\ml_{i-1}, \ml_j]$.

\vspace{1.5pt}

\item If $u=[\mr_i, \mr_j]$ with $i\ge j+2$, then
there exists an arrow $M_u\to M_v$ in $\Ga_{\hspace{-1pt}_{\mathscr{C}(Q)}}$ if and only if $v=[\mr_{i-1}, \mr_j]$ with $i>j+2$ or $v=[\mr_i, \mr_{j-1}]$.

\end{enumerate}

\end{Lemma}

\medskip
\medskip

The following result is essential in our investigation, and characterizes the rigidity of a pair of indecomposable objects of $\mathscr{C}(Q)$ by the non-crossing property of the corresponding arcs.

\medskip

\begin{Theo} \label{CrossingPair}

Let $u, v$ be arcs in $\mathcal{B}_\infty$. If $M_u, M_v$ are the corresponding objects in $\Ga_{\mathscr{C}(Q)}$, then $(u,v)$ is a crossing pair if and only if $\Hom_{\,\mathscr{C}(Q)}(M_u, M_v[1]) \ne 0,$ or equivalently, $\Hom_{\,\mathscr{C}(Q)}(M_v, M_u[1]) \ne 0.$

\end{Theo}

\noindent{\it Proof.} Let $u, v$ be distinct arcs in $\mathcal{B}_\infty$. Since $\mathscr{C}(Q)$ is $2$-CY, the last two statements stated in the proposition are equivalent. If one of $u, v$ is an upper arc and the other one is a lower arc, then one of $M_u, M_v$ lies in $\mathcal{R}_L$ and the other one lies in $\mathcal{R}_R$. In this case, the arcs $u, v$ do not cross. Since $\mathcal{R}_L, \, \mathcal{R}_R$ are orthogonal in $\mathscr{C}(Q)$ by Proposition \ref{idt-clust-cpt}(3), the results holds true.

Consider now the case where $u$ and $v$ are connecting arcs. Then $M_u, M_v\in \mathcal{C}.$ For any integer $i$, by Corollary \ref{tau-cross}, the arcs $u, v$ cross if and only if $\tau^i\hspace{-0.5pt}u,\, \tau^i\hspace{-0.5pt}v$ cross. On the other hand, since $\tau_{_\mathscr{C}}$ is an automorphism of $\mathscr{C}(Q)$, we have
$$\Hom_{\,\mathscr{C}(Q)}(M_{\tau^i u}, M_{\tau^i v}[1]) = \Hom_{\,\mathscr{C}(Q)}(\tau_{_\mathscr{C}}^i\hspace{-0.5pt}M_u, \tau_{_\mathscr{C}}^i\hspace{-0.5pt}M_v[1]) \cong \Hom_{\,\mathscr{C}(Q)}(M_u,M_v[1]).$$
Therefore, there is no loss of generality in assuming that $M_u$ and $\tau_{_\mathscr{C}}M_v=M_{\tau v}$ are preprojective representations in $\Ga_{\rep(Q)}$.

\vspace{1pt}

Suppose first that $u, v$ cross. By Lemma \ref{CombTech}(4), there is no loss of gene\-rality in assuming that $u=[\ml_p, \mr_q]$ and $v =[\ml_i, \mr_j]$ with $i < p$ and $j < q$. In view of Lemma \ref{arrows}(1), we obtain a path
$$M_u=M_{[\ml_p, \mr_q]} \longrightarrow  M_{[\ml_p, \mr_{q-1}]} \longrightarrow
\cdots \longrightarrow  M_{[\ml_p, \mr_{j+1}]} \longrightarrow  M_{[\ml_{p-1}, \mr_{j+1}]} \vspace{1pt} \hspace{100pt}$$
$$ \hspace{95pt}  \longrightarrow  M_{[\ml_{p-2}, \mr_{j+1}]} \longrightarrow  \cdots \longrightarrow M_{[\ml_{i+1}, \mr_{j+1}]}=M_{\tau v}\vspace{2pt}$$ in $\C$.
Since $\C$ is a standard component of $\Ga_{D^b(\rep(Q))}\vspace{1pt}$ of shape $\Z\mathbb{A}_\infty^\infty$; see (\ref{cpt-der-idt}), we deduce from Corollary \ref{A-d-inf} that $\Hom_{D^b(\rep(Q))}(M_u, M_{\tau v}) \ne 0,$ and consequently,
$$\Hom_{\hspace{0.5pt}\mathscr{C}(Q)}(M_u, M_v[1])=\Hom_{\hspace{0.5pt}\mathscr{C}(Q)}(M_u, \tau_{_\mathscr{C}}\hspace{-0.5pt}M_v) = \Hom_{\hspace{0.5pt}\mathscr{C}(Q)}(M_u, M_{\tau v}) \ne 0.$$
Suppose conversely that 
$\Hom_{\hspace{0.5pt}\mathscr{C}(Q)}(M_u, M_{\tau v}) \ne 0.$ Since $M_u, M_{\tau v}$ are assumed to be representations, by Lemma \ref{maps}(1), we have $\Hom_{D^b(\rep(Q))}(M_u, M_{\tau v}) \ne 0$ or
$$\begin{array}{rcl}
\Hom_{D^b(\rep(Q))}(M_v, M_{\tau u}) & \cong & \Hom_{D^b(\rep(Q))}(\tau_{\hspace{-1.5pt}_D} M_v, \tau_{\hspace{-1.5pt}_D} M_{\tau u})\\ \vspace{-8pt} \\
&=&\Hom_{D^b(\rep(Q))}(M_{\tau v},\tau_{\hspace{-1.5pt}_D}^2 M_u) \ne 0.
\end{array}$$
Since $\mathcal{C}$ is standard in $D^b(\rep(Q))$, we obtain a path $M_u \rightsquigarrow M_{\tau v}$ or $M_v \rightsquigarrow M_{\tau u}$, that is, a path $M_{[\ml_p, \mr_q]} \rightsquigarrow M_{[\ml_{i+1}, \mr_{j+1}]}\vspace{1pt}$ or $M_{[\ml_i, \mr_j]} \rightsquigarrow M_{[\ml_{p+1}, \mr_{q+1}]}\vspace{1pt}$ in $\mathcal{C}$. By Lemma \ref{arrows}(1), $p \le i+1$ and $q \le j+1$ in the first case, and $i \le p+1$ and $j \le q+1$ in the second case. By Lemma \ref{CombTech}(4), the arcs $u, v$ cross.

Consider next the case where $v, u$ are upper arcs, say $u =[\ml_p, \ml_q]\vspace{1pt}$ and $v =[\ml_i, \ml_j]$ with $p \le q-2$ and $i \le j-2$. Then $M_u, M_v\in \mathcal{R}_L$. Assume first that $u$ and $v$ cross. By Lemma \ref{CombTech}(5), we may assume that $i < p < j < q$. In view of Lemma \ref{arrows}(3), we see that $\mathcal{R}_L$ contains a path
$$M_u=M_{[\ml_p, \ml_q]} \longrightarrow  M_{[\ml_p, \ml_{q-1}]} \longrightarrow
\cdots \longrightarrow  M_{[\ml_p, \ml_{j+2}]} \longrightarrow  M_{[\ml_p, \ml_{j+1}]} \vspace{1pt} \hspace{100pt} $$
$$ \hspace{67pt} \longrightarrow  M_{[\ml_{p-1}, \ml_{j+1}]} \longrightarrow  M_{[\ml_{p-2}, \ml_{j+1}]} \longrightarrow  \cdots \longrightarrow M_{[\ml_{i+1}, \ml_{j+1}]}= M_{\tau v}, \vspace{2pt}$$ lying in the forward rectangle of $M_u.$ By Proposition \ref{rectangle}, $\Hom_{\hspace{0.5pt}\rep(Q)}(M_u, M_{\tau v}) \ne 0,$ and consequently, we obtain
$$\Hom_{\mathscr{C}(Q)}(M_u, M_v[1]) = \Hom_{\mathscr{C}(Q)}(M_u, M_{\tau v})\ne 0.$$

Conversely, assume that $\Hom_{\mathscr{C}(Q)}(M_u, M_v[1])=\Hom_{\mathscr{C}(Q)}(M_u, M_{\tau v})$ is non-zero. By Lemma \ref{maps}(1), we have $\Hom_{D^b(\rep(Q))}(M_u, M_{\tau v}) \ne 0$ or
$$\begin{array}{rcl}
\Hom_{D^b(\rep(Q))}(M_v, M_{\tau u}) & \cong & \Hom_{D^b(\rep(Q))}(\tau_{\hspace{-1.5pt}_D} M_v, \tau_{\hspace{-1.5pt}_D} M_{\tau u})\\ \vspace{-8pt} \\
&=&\Hom_{D^b(\rep(Q))}(M_{\tau v},\tau_{\hspace{-1.5pt}_D}^2 M_u) \ne 0.
\end{array}$$
Suppose that the first case occurs. We claim that $i < p < j < q $. Since $\mathcal{R}_L$ is standard in $D^b(\rep(Q))$ of shape $\mathbb{Z}\mathbb{A}_\infty$, by Proposition \ref{rectangle}, $M_{\tau v}$ lies in the forward rectangle of $M_u$. In particular, $\mathcal{R}_L$ contains a path $\rho: M_u=M_{[\ml_p, \ml_q]}
\rightsquigarrow M_{[\ml_{i+1}, \ml_{j+1}]}.$ By Lemma \ref{arrows}(3), we obtain $p\ge i+1>i$ and $q\ge j+1>j$.
If $\rho$ is trivial, then $p=i+1$ and $q=j+1$, and hence, $p=i+1\le j-1<j$ by our assumption. If $\rho$ is a sectional path, then we deduce from Lemma \ref{arrows}(3) that $p=i+1$ and $q>j+1$, or else, $p>i+1$ and $q=j+1$. In both cases, since $i+1\le j-1$ and $p\le q-2$ by our assumption, we obtain $p<j$. If $\rho$ is non-sectional, then $M_{[\ml_p, \ml_q]}$ is not quasi-simple, that is, $p<q-2$. Recalling that $M_{[\ml_p, \ml_q]}= L_{p+1}^-\cap L_{q-1}^+$ and $M_{[\ml_{i+1}, \ml_{j+1}]}= {L}_{i+2}^-\cap {L}_i^+$, we may choose $\rho$ to be the composite of a sectional path $M_{[\ml_p, \ml_q]} \rightsquigarrow M$ contained in the coray $ L_{p+1}^-$ and a sectional path $M\rightsquigarrow M_{[\ml_{i+1}, \ml_{i+1}]}$ contained in the ray $L_i^+$. Then, $M= {L}_{p+1}^-\cap {L}_i^+,\vspace{1pt}$ and consequently, $p+1\le j$, that is, $p<j$. This proves our claim.
Similarly, if the second case occurs, we can show that $p<i<q<j$. By Lemma \ref{CombTech}, the arcs $u, v$ cross.
The case where both $u$ and $v$ are lower arcs can be treated in a similar way.

Consider further the case where one of $u, v$ is a connecting arc and the other one is an upper arc, say $u =[\ml_p, \ml_q]$ with $p \le q-2$ and $v =[\ml_i,\mr_j]$. We have $M_u=  {L}_{p+1}^-\cap {L}_{q-1}^+\in \mathcal{R}_L$ and $M_v={L}_i\cap {L}_j\in \mathcal{C}.$ Using Corollary \ref{tau-cross} and the fact that $\tau_{_\mathscr{C}}$ is an automorphism of $\mathscr{C}(Q)$, we may assume that $M_v$ is a preprojective representation. By Corollary \ref{special-maps}, $\Hom_{\mathscr{C}(Q)}(M_v, M_u[1]) \ne 0$ if and only if $\Hom_{D^b(\rep(Q))}(M_v, M_u[1]) \ne 0.\vspace{0.5pt}$ Since
$M_u[1]=\tau_{\hspace{-0.4pt}_\mathscr{C}}M_u=M_{\tau u}=M_{[\ml_{p+1}, \ml_{q+1}]}$, the latter condition by Proposition \ref{LemmaRays}(2) is equivalent to $M_{[\ml_{p+1}, \ml_{q+1}]}\in \mathcal{W}(\tau_{\hspace{-0.4pt}_{\mathscr{C}}}^i M(q_0))$. Since
$\tau_{\hspace{-0.4pt}_{\mathscr{C}}}^i M(q_0)= \tau_{\hspace{-0.4pt}_{\mathscr{C}}}^{i+1}S_L= {L}_{i+1}^+\cap  {L}_{i+1}^-$ and $M_{[\ml_{p+1}, \ml_{q+1}]}= {L}_{p+2}^-\cap  {L}_q^+,\vspace{1pt}$ we see that $M_{[\ml_{p+1}, \ml_{q+1}]}\vspace{0.5pt}\in \mathcal{W}(\tau_{\hspace{-0.4pt}_{\mathscr{C}}}^i M(q_0))$ if and only if $i+1\ge p+2$ and $q\ge i+1$, that is, $p < i < q$. This last condition by Lemma \ref{CombTech} is equivalent to $u, v$ crossing. The case where one of $u, v$ is a connecting arc and the other one is a lower arc can be treated in a similar manner. The proof of the Theorem is completed.

\medskip

Given a strictly additive subcategory $\T$ of $\mathscr{C}(Q)$, we shall write ${\rm arc}(\T)$ for the set of arcs $a_{_T}$ with
$T \in \T \cap \mathscr{F}(Q).$ As an immediate consequence of Theorem \ref{CrossingPair} and Lemma \ref{rigidity}, we obtain the following result.

\medskip

\begin{Theo}\label{WCT-Sub}

If $\mathcal{T}$ is a strictly additive subcategory of $\mathscr{C}(Q)$, then $\T$ is weakly cluster-tilting if and only if ${\rm arc}(\T)$ is a triangulation of $\mathcal{B}_\infty$.

\end{Theo}

\medskip

Next, we shall describe the triangulations of $\mathcal{B}_{\infty}$ that correspond to the cluster-tilting subcategories of $\mathscr{C}(Q)$. Consider again the bijection $\varphi: \mathscr{F}(Q) \to {\rm arc}(\mathcal{B}_\infty).$ For each object $M\in \mathscr{F}(Q)$, we have $M=\varphi^{-1}(a_{\hspace{-0.5pt}_M})=M_{a_{\hspace{-0.5pt}_M}}$, and hence, $$\tau_{\hspace{-0.4pt}_{\mathscr{C}}}M= \tau_{\hspace{-0.4pt}_{\mathscr{C}}}M_{a_{\hspace{-0.5pt}_M}}= M_{\tau a_{\hspace{-0.5pt}_M}} = \varphi^{-1}(\tau a_{\hspace{-0.5pt}_M}),$$ that is,
$\tau a_{\hspace{-0.5pt}_M}=\varphi(\tau_{\hspace{-0.4pt}_{\mathscr{C}}}M) = a_{\hspace{-0.6pt}_{\tau_{\hspace{-0.4pt}_{\mathscr{C}}}\hspace{-1pt}M}}.$
The following easy result will be useful.

\medskip

\begin{Cor}\label{cross-lemma}

If $M, N$ are in $\mathscr{F}(Q)$, then $(a_{\hspace{-1pt}_N}, \tau a_{\hspace{-0.5pt}_M})$ is crossing if and only if
$\Hom_{\,\mathscr{C}(Q)}(M, N)\ne 0$.

\end{Cor}

\noindent{\it Proof.} Let $M, N\in \mathscr{F}(Q)$. Since $\mathscr{C}(Q)$ is 2-Calabi-Yau, we obtain
$$\Hom_{\mathscr{C}(Q)\vspace{0.5pt}}(N, \tauc M[1])=\Hom_{\mathscr{C}(Q)\vspace{0.5pt}}(N, M[2])\cong D\Hom_{\mathscr{C}(Q)\vspace{0.5pt}}(M, N).$$
By Theorem \ref{CrossingPair}, $a_{\hspace{-1pt}_N}$ and $\tau a_{\hspace{-0.5pt}_M}$ cross if and only if
$\Hom_{\mathscr{C}(Q)\vspace{0.5pt}}(N, \tauc M[1])\ne 0$, or equivalently, $\Hom_{\,\mathscr{C}(Q)}(M, N)\ne 0$. The proof of the corollary is completed.

\medskip

\begin{Lemma} \label{Ind-quo}
Let $\T$ be a weakly cluster-tilting subcategory of $\mathscr{C}(Q)$, containing a non-zero morphism $f: X\to Y$ which embeds in an exact triangle
$$\xymatrix{X \ar[r]^-f & Y \ar[r]^g & Z \ar[r]^-h & X[1]}$$
in $\mathscr{C}(Q).$ If $X, Y$ are indecomposable, then so is $Z$.

\end{Lemma}

\noindent{\it Proof.} Suppose that $X, Y$ are indecomposable and $Z = Z_1 \oplus Z_2$ with $Z_1, Z_2$ non-zero. Write $g=(g_1, g_2)^T$ and $h=(h_1, h_2)$. Assume that $Z_1[-1] \in \T$. Then $\Hom_{\mathscr{C}(Q)}(Y, Z_1)=0$, and hence, $g_1=0$. Since $g$ is a pseudo-kernel of $h$, we deduce that $h_1$ is a monomorphism. Since $\mathscr{C}(Q)$ is triangulated, $h_1$ is a section, and since $X$ is indecomposable, $h_1$ is an isomorphism. Observing that $f[1]\circ h=0$, we obtain $f=0$, a contradiction. Thus $Z_1[-1] \not \in \T$, and similarly, $Z_2[-1] \not \in \T$. Since $\T$ is weakly cluster-tilting, there exists $T_i \in \T$ such that $\Hom_{\mathscr{C}(Q)}(T_i, Z_i) \ne 0,\vspace{0.5pt}$ for $i=1, 2.$ Set $T = X \oplus Y \oplus T_1 \oplus T_2$. Since $\Hom_{\mathscr{C}(Q)\vspace{0.5pt}}(T, X[1])=0,$
applying $\Hom_{\mathscr{C}(Q)\vspace{0.5pt}}(T, -)$ to the triangle stated in the lemma yields a projective presentation of the right ${\rm End}(T)$-module $\Hom_{\mathscr{C}(Q)\vspace{1pt}}(T, Z)\vspace{0.5pt}$ as follows:
$$\xymatrix{\Hom_{\mathscr{C}(Q)\vspace{0.5pt}}(T, X) \ar[r] & \Hom_{\mathscr{C}(Q)\vspace{0.5pt}}(T, Y) \ar[r] & \Hom_{\mathscr{C}(Q)\vspace{0.5pt}}(T, Z) \ar[r] & 0.}$$ Since $\Hom_{\mathscr{C}(Q)\vspace{1pt}}(T, Z)\vspace{0.5pt}=\Hom_{\mathscr{C}(Q)\vspace{0.5pt}}(T, Z_1)\oplus \Hom_{\mathscr{C}(Q)\vspace{0.5pt}}(T, Z_2)$
is decomposable, so is $\Hom_{\mathscr{C}(Q)\vspace{0.5pt}}(T, Y)$, a contradiction to $Y$ being indecomposable in $\mathscr{C}(Q)$. The proof of the lemma is completed.

\medskip

If $\mathfrak{p}, \mathfrak{q}$ are two marked points of $\mathcal{B}_\infty$, then we define $\mathfrak{p} < \mathfrak{q}$ as follows. If $\mathfrak{p}$ is an upper marked point and $\mathfrak{q}$ is a lower marked point, then $\mathfrak{p} < \mathfrak{q}$ and $\mathfrak{q} \not < \mathfrak{q}$. If $\mathfrak{p} = \ml_i$ and  $\mathfrak{q} = \ml_j$, then $\mathfrak{p} < \mathfrak{q}$ if and only if $i < j$. If $\mathfrak{p} = \mr_i$ and  $\mathfrak{q} = \mr_j$, then $\mathfrak{p} < \mathfrak{q}$ if and only if $i < j$.

\begin{Lemma} \label{LemmaFactorization}

Let $\T$ be a weakly cluster-tilting subcategory of $\mathscr{C}(Q)$, and let $M, N, L$ be indecomposable objects in $\mathscr{C}(Q)$ with $M, N\in \T$. If $f: M\to N$ and $g: N\to L$ are non-zero morphisms in $\mathscr{C}(Q)$, then $\Hom_{\mathscr{C}(Q)}(M, L)$ is generated by $gf$ over $k$.

\end{Lemma}

\noindent{\it Proof.} We may assume that $M, N, L$ all lie in the fundamental domain $\mathscr{F}(Q)$ and  $\Hom_{\mathscr{C}(Q)}(M, L)\ne 0$. Let $f: M\to N$ and $g: N\to L$ be non-zero morphisms in $\mathscr{C}(Q)$. Since $\Hom_{\mathscr{C}(Q)}(M, L)$ is of $k$-dimension at most one by Proposition \ref{onedim}, it suffices to show that $gf\ne 0$. Assume on the contrary that $gf=0$. In particular, $f$ is not an isomorphism. Since $N$ is indecomposable, $f$ is not a section. Thus, $\mathscr{C}(Q)$ contains a non-split exact triangle
$$\xymatrix{(*) \qquad  M \ar[r]^-f & N \ar[r]^v & C \ar[r]^-w & M[1].}$$
By Lemma \ref{Ind-quo}, $C$ is indecomposable, which we may assume to be in $\mathscr{F}(Q)$. Since each of $f, v,$ $w$ is non-zero, by Proposition \ref{onedim}, each of $\Hom_{\mathscr{C}(Q)\vspace{0.5pt}}(M, N)$,
$\Hom_{\mathscr{C}(Q)\vspace{0.5pt}}(N, C)$ and $\Hom_{\mathscr{C}(Q)}(C, M[1])$ is one-dimensional. Since $f$ and $v$ are not isomorphisms, $M, N$ and $C$ are pairwise non-isomorphic. We shall need two crucial facts as follows.

\smallskip

(1) {\it Each pair of arcs $(a_{\hspace{-1pt}_N}, \tau a_{\hspace{-0.5pt}_M}), (a_{\hspace{-0.5pt}_C}, \tau a_{\hspace{-1pt}_N}), (a_{\hspace{-0.5pt}_C}, a_{\hspace{-0.5pt}_M})$ is crossing.}

\smallskip

(2) {\it Each pair of arcs $(a_{\hspace{-1pt}_N}, a_{\hspace{-0.5pt}_M}), (a_{\hspace{-0.5pt}_C}, a_{\hspace{-1pt}_N}), (a_{\hspace{-0.5pt}_C}, \tau a_{\hspace{-0.5pt}_M})$ is non-crossing.}

\smallskip

Indeed, Statement (1) follows immediately from Corollary \ref{cross-lemma} and Theorem \ref{CrossingPair}. Since $\T$ is rigid, $\Hom_{\mathscr{C}(Q)}(M, N[1])=0$, and by Theorem \ref{CrossingPair}, $(a_{\hspace{-0.5pt}_M}, a_{\hspace{-1pt}_N})$ is non-crossing. Moreover, since $\Hom_{\mathscr{C}(Q)}(C,C[1])=0$ by Corollary \ref{rigid-obj}, applying $\Hom_{\mathscr{C}(Q)\vspace{0.5pt}}(C, -)$ to the triangle $(*)$ yields an exact sequence
$$\xymatrix{\Hom_{\mathscr{C}(Q)\vspace{0.5pt}}(C, C) \ar[r] & \Hom_{\mathscr{C}(Q)\vspace{0.5pt}}(C, M[1])
\ar[r] & \Hom_{\mathscr{C}(Q)\vspace{0.5pt}}(C, N[1]) \ar[r] & 0.}$$
Since $\Hom_{\mathscr{C}(Q)}(C, M[1])$ is one-dimensional, we have $\Hom_{\mathscr{C}(Q)}(C,N[1])=0,\vspace{0.5pt}$ and by Theorem \ref{CrossingPair}, $(a_{\hspace{-0.5pt}_C}, a_{\hspace{-1pt}_N})$ is non-crossing.
Further, since $\Hom_{\mathscr{C}(Q)}(M, M[1])=0,$ applying $\Hom_{\mathscr{C}(Q)\vspace{0.5pt}}(M, -)$ to the triangle $(*)$ yields an exact sequence
$$\xymatrix{\Hom_{\mathscr{C}(Q)\vspace{0.5pt}}(M, M) \ar[r] & \Hom_{\mathscr{C}(Q)\vspace{0.5pt}}(M, N)
\ar[r] & \Hom_{\mathscr{C}(Q)\vspace{0.5pt}}(M, C) \ar[r] & 0.}$$
Since  $\Hom_{\mathscr{C}(Q)\vspace{0.5pt}}(M, N)$ is one-dimensional, we have
$\Hom_{\mathscr{C}(Q)\vspace{0.5pt}}(M, C)=0$, and thus, $(a_{\hspace{-0.5pt}_C}, \tau a_{\hspace{-0.5pt}_M})$ is non-crossing by Corollary \ref{cross-lemma}. This establishes Statement (2).

\vspace{1pt}

Denote by $\mf{p}, \mf{q}$ with $\mf{p}<\mf{q}$ the endpoints of $a_{\hspace{-0.5pt}_M},$ and by $\tau \mf{p}, \tau \mf{q}$ the endpoints of $\tau a_{\hspace{-0.5pt}_M}$. Observe that $\mf{p}, \mf{q}, \tau \mf{p}, \tau \mf{q}$ are pairwise distinct. Let $\it\Omega$ be a simply connected region in $\mathcal{B}_\infty$ enclosed by a closed simple curve which is the composite of two edge arcs $u'\in [\mf{p}, \tau \mf{p}]$ and $v'\in [\mf{q}, \tau \mf{q}]$, an arc curve $u\in [\mf{p}, \tau \mf{q}]$ and a simple curve $v\in [\tau \mf{p}, \mf{q}]$.

We shall say that an arc $w$ in $\mathcal{B}_\infty$ {\it crosses} $\it\Omega$ if every simple curve in the isotopy class of $w$ intersects the interior region of $\it\Omega$. Since $(a_{\hspace{-1pt}_N}, \tau a_{\hspace{-0.5pt}_M})$ is crossing,
$a_{\hspace{-1pt}_N}$ crosses $\it\Omega$. Since $u'$ and $v'$ are edge curves, either $a_{\hspace{-1pt}_N}$ crosses both $u$ and $v$, or else, $a_{\hspace{-1pt}_N}$ has as endpoint one of the marked points $\mf{p}, \mf{q}, \tau \mf{p}, \tau \mf{q}.$ In the first situation, $a_{\hspace{-1pt}_N}$ crosses $a_{\hspace{-0.5pt}_M}$, a contradiction to Statement (2).

If $a_{\hspace{-1pt}_N}$ shares an endpoint with $\tau a_{\hspace{-0.5pt}_M}$, then this contradicts the fact that $(a_{\hspace{-1pt}_N}, \tau a_{\hspace{-0.5pt}_M})$ is crossing. Therefore, $a_{\hspace{-0.5pt}_M}, a_{\hspace{-0.5pt}_N}$ share an endpoint, which we may assume to be $\mf{p}$. Write $a_{\hspace{-0.5pt}_N} = [\mf{p}, \mr]$, where $\mr$ is a marked point different from any of $\mf{p}, \mf{q}, \tau \mf{p}, \tau \mf{q}$.

Since $(a_{\hspace{-0.5pt}_C}, \tau a_{\hspace{-0.5pt}_N})$ is crossing and $(a_{\hspace{-0.5pt}_C}, a_{\hspace{-0.5pt}_N})$ is non-crossing, a similar argument shows that $a_{\hspace{-0.5pt}_C}, a_{\hspace{-0.5pt}_N}$ share an endpoint. Since $(a_{\hspace{-0.5pt}_C}, a_{\hspace{-0.5pt}_M})$ is crossing, this common endpoint is different from $\mf{p}$, and hence, it is $\mr$.

Finally, since $(a_{\hspace{-0.5pt}_C},a_{\hspace{-0.5pt}_M})$ is crossing and $(a_{\hspace{-0.5pt}_C},\tau a_{\hspace{-0.5pt}_M})$ is non-crossing, $a_{\hspace{-0.5pt}_C}, \tau a_{\hspace{-0.5pt}_M}$ share an endpoint. Since $(a_{\hspace{-0.5pt}_C}, \tau a_{\hspace{-1pt}_N})$ is crossing, this endpoint is different from $\tau \mf{p}$, and hence, it is $\tau \mf{q}$. This yields $a_{\hspace{-0.5pt}_C} = [\mr, \tau \mf{q}]$.

\smallskip

Since $\Hom_{\mathscr{C}(Q)}(M, L)$ is nonzero, by Corollary \ref{cross-lemma}, $(a_{\hspace{-0.5pt}_L}, \tau a_{\hspace{-0.5pt}_M})$ is crossing. Similarly, since $g: N\to L$ is nonzero,  $(a_{\hspace{-0.5pt}_L}, \tau a_{\hspace{-0.5pt}_N})$ is crossing.
Let $\it\Omega'$ be a simply connected region in $\mathcal{B}_\infty$ enclosed by a closed simple curve which is the composite of three arcs $u_1\in [\tau \mf{p}, \tau \mr]$, $u_2\in [\tau \mr, \tau \mf{q}]$, and $u_3\in [\tau \mf{q}, \tau \mf{p}]$. Observe that $u_1$ is isotopic to $\tau a_N$ and $u_3$ is isotopic to $\tau a_M$. Since we know that $a_{\hspace{-0.5pt}_L}$ crosses both $\tau a_{\hspace{-0.5pt}_M}, \tau a_{\hspace{-0.5pt}_N}$, the arc $a_{\hspace{-0.5pt}_L}$ cannot cross the other arc, which is $u_2$. Consider now $\it\Omega''$ a simply connected region in $\mathcal{B}_\infty$ enclosed by a closed simple curve which is the composite of three arcs $v_1 = u_2\in [\tau \mf{q}, \tau \mr]$, $v_2\in [\tau \mr, \mr]$, and $v_3\in [\mr, \tau \mf{q}]$. Observe that $v_2$ is an edge arc. We claim that $a_{\hspace{-0.5pt}_L}$ does not cross $\it\Omega''$. Suppose the contrary. Since $v_2$ is an edge arc, either $a_{\hspace{-0.5pt}_L}$ crosses both $v_1, v_3$, or else $a_{\hspace{-0.5pt}_L}$ only crosses one of $v_1, v_3$ but then has endpoint $\tau \mr$ or $\mr$. Since we already know that $a_{\hspace{-0.5pt}_L}$ does not cross $u_2=v_1$, the first case does not occur. Hence, $a_{\hspace{-0.5pt}_L}$ crosses $v_3$ and has endpoint $\tau \mr$, which contradicts the fact that $(a_{\hspace{-0.5pt}_L}, \tau a_{\hspace{-0.5pt}_N})$ is crossing. This proves our claim that $a_{\hspace{-0.5pt}_L}$ does not cross $\it\Omega''$. Since $a_{\hspace{-0.5pt}_C}$ is isotopic to $v_3$, this yields that $a_{\hspace{-0.5pt}_L}$ does not cross $a_{\hspace{-0.5pt}_C}$.

By Theorem \ref{CrossingPair}, $\Hom_{\mathscr{C}(Q)}(C[-1],L)=0$. Applying $\Hom{\mathscr{C}(Q)}(-, L)$ to the exact triangle $(*)$, we obtain the following exact sequence
$$\xymatrix{\Hom_{\mathscr{C}(Q)\vspace{0.5pt}}(C, L) \ar[r] & \Hom_{\mathscr{C}(Q)\vspace{0.5pt}}(N, L)
\ar[r] & \Hom_{\mathscr{C}(Q)\vspace{0.5pt}}(M, L) \ar[r] & 0.}$$
Since both $\Hom_{\mathscr{C}(Q)}(N,L)$ and $\Hom_{\mathscr{C}(Q)}(M,L)$ are one-dimensional by Proposition \ref{onedim}, $\Hom_{\mathscr{C}(Q)}(f,L): \Hom_{\mathscr{C}(Q)}(N,L) \to \Hom_{\mathscr{C}(Q)}(M,L)$ is an isomorphism. In particular, $gf\ne 0$.
This completes the proof of the lemma.

\medskip


Now, we are ready to determine the cluster-tilting subcategories of $\mathscr{C}(Q)$ in terms of the triangulations of $\B_\infty$.

\medskip

\begin{Theo} \label{functioriallyfinite}

Let $Q$ be a quiver with no infinite path of type $\mathbb{A}_\infty^\infty$. The following three statements are equivalent for a strictly additive subcategory $\T$ of $\mathscr{C}(Q)$.

\vspace{-2pt}

\begin{enumerate}[$(1)$]

\item The subcategory $\T$ is cluster-tilting.

\item ${\rm arc}(\T)$ is a compact triangulation of $\B_\infty$.

\item ${\rm arc}(\T)$ contains some connecting arcs and every marked point in $\B_\infty$ is either ${\rm arc}(\T)$-bounded or a ${\rm arc}(\T)$-fountain base.

\end{enumerate}

\vspace{-2pt}

\noindent In this case, moreover, ${\rm arc}(\T)$ has at most two fountains, and if it has two fountains, then one is a left fountain and the other one is a right fountain.

\end{Theo}

%
%
%
%
\noindent{\it Proof.} In view of Theorem \ref{TheoClusterTilting} and Lemmas  \ref{Full-fountain} and \ref{LemmaFountain}, it suffices to show the equivalence of the first two statements.
Assume first that ${\rm arc}(\mathcal{T})$ is a compact triangulation of $\B_\infty$. In particular, $\,\T$ is weakly cluster-tilting. We shall need to prove that $\T$ is functorially finite. We first show that every indecomposable object $M \in \mathscr{C}(Q)$ admits a right $\T$-approximation. Since $\mathscr{C}(Q)$ is 2-Calabi-Yau, for each indecomposable object $N\in \mathscr{C}(Q)$, we have
$$\Hom_{\mathscr{C}(Q)}(N,M) \cong D\Hom_{\mathscr{C}(Q)}(\tau_{_\mathscr{C}}^-M,N[1]),$$ which is of $k$-dimension at most one by Proposition \ref{onedim}.

By the assumption, the set ${\rm arc}(\T)_{\tau^-a_{\hspace{-0.5pt}_M}}\vspace{0.5pt}$ of arcs of ${\rm arc}(\T)$ which cross $\tau^-a_{\hspace{-0.5pt}_M}$ is compact.  Let $\it\Sigma$ be a finite subset of ${\rm arc}(\T)_{\tau^-a_{\hspace{-0.5pt}_M}}$ satisfying the condition stated in Definition \ref{compact-set}. Observe that $\tau_{_\mathscr{C}}^-M=\tau_{_\mathscr{C}}^-M_{\hspace{-0.5pt}a_{\hspace{-0.5pt}_M}}=M_{\tau^-a_{\hspace{-0.5pt}_M}}$.
If $v\in \it\Sigma$, since $v$ crosses $\tau^-a_{\hspace{-0.5pt}_M}$, we deduce from Theorem \ref{CrossingPair} that
$$\Hom_{\mathscr{C}(Q)}(M_v, M)\cong D\Hom_{\mathscr{C}(Q)}(\tau_{_\mathscr{C}}^-M, M_v[1])=D\Hom_{\mathscr{C}(Q)}(M_{\tau^-a_{\hspace{-0.5pt}_M}}, M_v[1])\ne 0.$$ In particular, we may choose a nonzero morphism $f_v: M_v \to M$ for each $v\in \it\Sigma$.
Set $f=\oplus_{v\in \it\Sigma}\,f_v: \oplus_{v \in {\it\Sigma}}\, M_v \to M,$ which we claim is a right $\T$-approximation of $M$. Indeed, assume that $\Hom_{\mathscr{C}(Q)}(N,M) \ne 0$, for some indecomposable object $N \in \T$. Then $\Hom_{\mathscr{C}(Q)}(\tau_{_\mathscr{C}}^-M,N[1]) \ne 0$, and by Theorem \ref{CrossingPair}, $a_{\hspace{-1pt}_N}$ crosses $\tau^-a_{\hspace{-0.5pt}_M}$, that is, $a_{\hspace{-1pt}_N} \in {\rm arc}(\T)_{\tau^-a_{\hspace{-0.5pt}_M}}$. By assumption, there exists an arc $w\in \it\Sigma$ such that $a_{\hspace{-1pt}_N}$ crosses $\tau^-w.$ This implies that
$$\Hom_{{\mathscr{C}(Q)}}(N, M_w)\cong D\Hom_{\mathscr{C}(Q)}(\tau_{_\mathscr{C}}^-M_w, N[1])  \ne 0 .$$
Let $g_w: N\to M_w$ be a non-zero morphism. By Lemma \ref{LemmaFactorization}, every morphism $g: N \to M$ is a multiple of  $f_w \hspace{0.4pt} g_w$. In particular, $g$ factors through $f$. This establishes our claim. Using the dual of Lemma \ref{LemmaFactorization} and the compactness of ${\rm arc}(\T)_{\tau a_{\hspace{-0.5pt}_M}},\vspace{0.5pt}$ we may show that $M$ admits a left $\T$-approximation. 
This establishes the sufficiency.

Conversely, assume that $\T$ is functorially finite in $\mathscr{C}(Q)$. Fix $u \in {\rm arc}(\mathcal{B}_\infty).$ By assumption, $M_{\tau u}$ admits a minimal right $\T$-approximation $f: \oplus_{w \in \it\Sigma^-} M_w \to M_{\tau u}$, where $\it\Sigma^-$ is a finite subset of ${\rm arc}(\T)$. For any $w\in \it\Sigma^-,$ since $f$ is right minimal, restricting $f$ to $M_w$ yields a non-zero morphism $f_w: M_w\to M_{\tau u}$. Observing that
$M_u=\tau_{_\mathscr{C}}^-M_{\tau u}$, we obtain
$$\Hom_{\mathscr{C}(Q)}(M_u, M_w[1])=\Hom_{\mathscr{C}(Q)}(\tau_{_\mathscr{C}}^-M_{\tau u}, M_w[1])\cong D\Hom_{\mathscr{C}(Q)}(M_w, M_{\tau u}) \ne 0.$$
By Theorem \ref{CrossingPair}, $w$ crosses $u$. This shows that $\it\Sigma^-\subseteq {\rm arc}(\T)_u$. On the other hand, for $v\in {\rm arc}(\T)_u$, since $M_{\tau u}=\tau_{_\mathscr{C}} M_u=M_u[1]$, we deduce from Theorem \ref{CrossingPair} that
$$\Hom_{\mathscr{C}(Q)}(M_v, M_{\tau u})=\Hom_{\mathscr{C}(Q)}(M_v, M_u[1])\cong D\Hom_{\mathscr{C}(Q)}(M_u, M_v[1])\ne 0.$$
Thus, there exists a nonzero morphism $g: M_v\to M_{\tau u},$ which factors through $f: \oplus_{w \in \it\Sigma^-} M_w \to M_{\tau u}$. In particular, there exists $v_1 \in \it \Sigma^-$ such that
$$\Hom_{\mathscr{C}(Q)}(M_{\tau^- v_1}, M_v[1])=\Hom_{\mathscr{C}(Q)}(\tau^-_{_\mathscr{C}}M_{v_1}, M_v[1])\cong D\Hom_{\mathscr{C}(Q)}(M_v, M_{v_1}) \ne 0.$$ That is, $v$ crosses $\tau^-v_1$. Similarly, considering a minimal left $\T$-approximation for $M_{\tau^-u}$, we may show that there exists a finite subset $\it\Sigma^+$ of ${\rm arc}(\T)_u$ such that every $v\in {\rm arc}(\T)_u$ crosses $\tau v_2$ for some arc $v_2\in \it\Sigma^+$. In particular, ${\rm arc}(\T)_u$ is compact. Hence, ${\rm arc}(\T)$ is compact. The proof of the theorem is completed.

\medskip
\medskip

Here is an example of a triangulation of $\B_\infty$, having two fountains, corresponding to a cluster-tilting subcategory of $\mathscr{C}(Q)$.

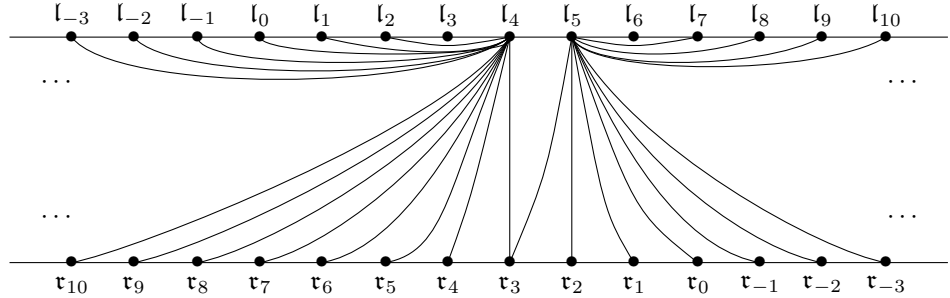
\begin{figure}[h]
  \centering
  \begin{tikzpicture}[xscale=2.50,yscale=1.5]

    \path (-2.4,0.6) node{$\cdots$};
    \path (2.1,0.6) node{$\cdots$};

   \path (-2.4,-0.6) node{$\cdots$};
    \path (2.1,-0.6) node{$\cdots$};

    \draw (-2.66,1) -- (2.33,1);
    \draw (-2.66,-1) -- (2.33,-1);

     \node at (-2.33,1.2) {$\ml_{-3}$};
    \node at (-2.0,1.2) {$\ml_{-2}$};
    \node at (-1.66,1.2) {$\ml_{-1}$};
    \node at (-1.33,1.2) {$\ml_0$};
    \node at (-1,1.2) {$\ml_1$};
 \node at (-0.66,1.2) {$\ml_{2}$};
    \node at (-0.33,1.2) {$\ml_{3}$};
    \node at (0,1.2) {$\ml_{4}$};
    \node at (0.33,1.2) {$\ml_5$};
    \node at (0.66,1.2) {$\ml_6$};
 \node at (1,1.2) {$\ml_{7}$};
    \node at (1.33,1.2) {$\ml_{8}$};
    \node at (1.66,1.2) {$\ml_{9}$};
    \node at (2,1.2) {$\ml_{10}$};

    \node at (-2.33,-1.2) {$\mr_{10}$};
    \node at (-2.0,-1.2) {$\mr_{9}$};
    \node at (-1.66,-1.2) {$\mr_{8}$};
    \node at (-1.33,-1.2) {$\mr_7$};
    \node at (-1,-1.2) {$\mr_6$};
 \node at (-0.66,-1.2) {$\mr_{5}$};
    \node at (-0.33,-1.2) {$\mr_{4}$};
    \node at (0,-1.2) {$\mr_{3}$};
    \node at (0.33,-1.2) {$\mr_2$};
    \node at (0.66,-1.2) {$\mr_1$};
 \node at (1,-1.2) {$\mr_{0}$};
    \node at (1.33,-1.2) {$\mr_{-1}$};
    \node at (1.66,-1.2) {$\mr_{-2}$};
    \node at (2,-1.2) {$\mr_{-3}$};

    \node at (-2.33,1.00) {$\bullet$};
    \node at  (-2.00,1.00) {$\bullet$};
    \node at  (-1.66,1.0){$\bullet$};
    \node at  (-1.33,1.0){$\bullet$};
    \node at  (-1.00,1.0){$\bullet$};
    \node at  (-0.66,1.0){$\bullet$};
    \node at  (-0.33,1.0){$\bullet$};
    \node at  (0.00,1.0){$\bullet$};
    \node at  (0.33,1.0){$\bullet$};
    \node at  (0.66,1.0){$\bullet$};
    \node at  (1.00,1.0){$\bullet$};
    \node at  (1.33,1.0){$\bullet$};
    \node at  (1.66,1.0){$\bullet$};
    \node at  (2.00,1.0){$\bullet$};

    \node at (-2.33,-1.0){$\bullet$};
    \node at (-2.00,-1.0){$\bullet$};
    \node at (-1.66,-1.0){$\bullet$};
   \node at (-1.33,-1.0){$\bullet$};
   \node at (-1.00,-1.0){$\bullet$};
    \node at (-0.66,-1.0){$\bullet$};
    \node at (-0.33,-1.0){$\bullet$};
    \node at (0.00,-1.0){$\bullet$};
    \node at (0.33,-1.0){$\bullet$};
    \node at(0.66,-1.0){$\bullet$};
    \node at (1.00,-1.0){$\bullet$};
   \node at (1.33,-1.0){$\bullet$};
    \node at(1.66,-1.0){$\bullet$};
    \node at (2.00,-1.0){$\bullet$};

\draw (-2.33, -1) .. controls (-2, -0.9) and (-0.33, 0.2) .. (0, 1);
\draw (-2, -1) .. controls (-1.66, -0.9) and (-0.33, 0.1) .. (0, 1);
\draw (-1.66, -1) .. controls (-1.33, -0.9) and (-0.33, 0) .. (0, 1);
\draw (-1.33, -1) .. controls (-1, -0.9) and (-0.33, -0.1) ..  (0, 1);
\draw (-1, -1) .. controls (-0.66, -0.9) and (-0.33, -0.2) ..  (0, 1);
\draw (-0.66, -1) .. controls (-0.33, -0.9) and (-0.33, -0.3) .. (0, 1);
\draw (-0.33, -1) .. controls (-0.2, -0.4) .. (0, 1);
\draw (0, -1) -- (0, 1);
\draw (0, -1) .. controls (0.2, -0.2).. (0.33, 1);
\draw (0.33, -1) -- (0.33, 1);

\draw (0.66, -1) .. controls (0.5, -0.55) .. (0.33, 1);
\draw (1.0, -1) .. controls (0.66, -0.55) .. (0.33, 1);
\draw (1.33, -1) .. controls (1, -0.9) and (0.66, -0.3) .. (0.33, 1);
\draw (1.66, -1) .. controls (1.33, -0.9) and (0.66, -0.2) .. (0.33, 1);
\draw (2, -1) .. controls (1.66, -0.9) and (0.66, 0) .. (0.33, 1);

\draw (-0.66, 1) .. controls (-0.33, 0.9) .. (0,1);
\draw (-1, 1) .. controls (-0.5, 0.8) and (-0.1, 0.85).. (0,1);
\draw (-1.33, 1) .. controls (-1.30, 0.8) and (0.1, 0.8) .. (0,1);
\draw (-1.66, 1) .. controls (-1.65, 0.7) and (-0.2, 0.7) .. (0,1);
\draw (-2, 1) .. controls (-1.8, 0.6) and (-0.3, 0.6) .. (0,1);
\draw (-2.33, 1) .. controls (-2.1, 0.5) and (-0.4, 0.5) .. (0,1);

\draw (0.33, 1) .. controls (0.66, 0.9) .. (1,1);
\draw (0.33, 1) .. controls (0.5, 0.8) and (1, 0.8) .. (1.33,1);
\draw (0.33, 1) .. controls (0.4, 0.7) and (1.33, 0.7) .. (1.66,1);
\draw (0.33, 1) .. controls (0.4, 0.65) and (1.66, 0.65) .. (2.0,1);

  \end{tikzpicture}
\label{fig:triangulation2}
\end{figure}

%
%
%
%

\medskip


\medskip

Let $\T$ be a cluster-tilting subcategory of $\mathscr{C}(Q)$. We would like to have a combinatorial description of the irreducible morphisms of $\T$ in terms of the arcs of ${\rm arc}(\T)$. This would give a combinatorial description of the quiver $Q_{_\T}$ of $\T$. Given two arcs $u, v\in {\rm arc}(\T)$, we shall write $u \vdash v$ if $u, v$ share an endpoint $\mathfrak{p}$ and $v$ is obtained by rotating $u$ in the counter-clockwise direction around $\mathfrak{p}$.

\begin{Prop} \label{PropIrredMorphisms}

Let $\T$ be a cluster-tilting subcategory of $\mathscr{C}(Q)$. If $M, N$ are non-isomorphic objects in $\mathscr{F}(Q) \cap \mathcal{T}$, then $\Hom_{\,\mathscr{C}(Q)}(M,N) \ne 0$ if and only if $a_{\hspace{-0.5pt}_M} \vdash a_{\hspace{-1pt}_N}$.

\end{Prop}

\noindent{\it Proof.} Let $M, N\in \mathscr{F}(Q) \cap \mathcal{T}$ be non-isomorphic. Suppose that $\Hom_{\mathscr{C}(Q)}(M,N) \ne 0$. Then $\Hom_{\mathscr{C}(Q)}(M,\tau_{_{\mathscr{C}}}^- N[1]) \ne 0$. Since $\T$ is cluster-tilting, $\Hom_{\mathscr{C}(Q)}(M,N[1]) = 0$. Hence, by Theorem \ref{CrossingPair}, $a_{\hspace{-0.5pt}_M}, \tau^-a_{\hspace{-1pt}_N}$ cross but $a_{\hspace{-0.5pt}_M}, a_{\hspace{-1pt}_N}$ do not cross. This clearly means that $a_{\hspace{-0.5pt}_M}, a_{\hspace{-1pt}_N}$ share an endpoint $\mathfrak{p}$. Suppose first that $a_{\hspace{-1pt}_N} = [\ml_p, \ml_q]$ with $p < q-1$. The conditions that $a_{\hspace{-0.5pt}_M}, \tau^-a_{\hspace{-1pt}_N}$ cross but $a_{\hspace{-0.5pt}_M}, a_{\hspace{-1pt}_N}$ do not cross clearly imply that $\mathfrak{p} = \ml_p$.  If $a_{\hspace{-0.5pt}_M}$ is not a connecting arc, then $a_{\hspace{-0.5pt}_M} = [\ml_p, \ml_r]$ with $r > q+1$ or $a_{\hspace{-0.5pt}_M} = [\ml_r, \ml_p]$ with $r < p-1$. In both cases, we see that $a_{\hspace{-0.5pt}_M} \vdash a_{\hspace{-1pt}_N}$. Clearly, if $a_{\hspace{-0.5pt}_M}$ is a connecting arc, then $a_{\hspace{-0.5pt}_M} \vdash a_{\hspace{-1pt}_N}$ as well. A similar argument handles the case where $a_{\hspace{-1pt}_N}$ lies on the lower boundary component. So assume $a_{\hspace{-1pt}_N} = [\ml_p, \mr_q]$. Then the conditions that $a_{\hspace{-0.5pt}_M}, \tau^-a_{\hspace{-1pt}_N}$ cross but $a_{\hspace{-0.5pt}_M}, a_{\hspace{-1pt}_N}$ do not cross give the following possibilities for $a_{\hspace{-0.5pt}_M}$. If $\mathfrak{p} = \ml_p$, then $a_{\hspace{-0.5pt}_M} = [\ml_p, \mr_i]$ with $i > q$ or $a_{\hspace{-0.5pt}_M} = [\ml_t, \ml_p]$ with $t < p-1$. If $\mathfrak{p} = \mr_q$, then $a_{\hspace{-0.5pt}_M} = [\ml_j, \mr_q]$ with $j > p$  or  $a_{\hspace{-0.5pt}_M} = [\mr_q,\mr_s]$ with $s < q-1$. In all four cases, $a_{\hspace{-0.5pt}_M} \vdash a_{\hspace{-1pt}_N}$. This proves the necessity. Conversely, assume that $a_{\hspace{-0.5pt}_M} \vdash a_{\hspace{-1pt}_N}$. Then $a_{\hspace{-0.5pt}_M}, a_{\hspace{-1pt}_N}$ share some endpoint. With no loss of generality, we may assume that $\ml_p$ is a common endpoint of $a_{\hspace{-0.5pt}_M}, a_{\hspace{-1pt}_N}$. Then a small neighborhood of $\ml_p$ looks as follows:

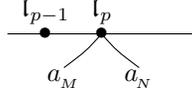
\begin{figure}[h]
  \centering
  \begin{tikzpicture}[xscale=2.50,yscale=1.5]

    \draw (-0.5,1) -- (0.5,1);

    \node at (0,1.2) {$\ml_{p}$};
 \node at (-0.3,1.2) {$\ml_{p-1}$};
 \node at (0,1) {$\bullet$};
 \node at (-0.3,1) {$\bullet$};
 \node at (-0.2,0.6) {$a_{\hspace{-0.5pt}_M}$};
 \node at (0.2,0.6) {$a_{\hspace{-1pt}_N}$};

\draw (0,1) .. controls (-0.08,0.85) .. (-0.2, 0.7);
\draw (0,1) .. controls (0.08,0.85) .. (0.2, 0.7);
\end{tikzpicture}
\label{fig:Neighborhood}
\end{figure}

\vspace{-5pt}

\noindent In view of this figure, we see that  $a_{\hspace{-0.5pt}_M}$ crosses $\tau^- a_{\hspace{-1pt}_N}$. By Theorem \ref{CrossingPair}, we obtain \vspace{-3pt}

$$\Hom_{\mathscr{C}(Q)}(M,N) \cong \Hom_{\mathscr{C}(\T)}(M,\tau_{_{\mathscr{C}}}^- N[1]) \ne 0.\vspace{2pt}$$
The proof of the proposition is completed.

\medskip

Let $\T$ be a cluster-tilting subcategory of $\mathscr{C}(Q)$. Given two arcs $u, v \in {\rm arc}(\T)$, we shall say that $v$ {\it covers} $u$ with respect to $\vdash$ provided that $u \vdash v$ and there exists no $w \in {\rm arc}(\T)$ such that $u \vdash w$ and $w \vdash v$. As an easy consequence of Proposition \ref{PropIrredMorphisms}, Lemmas \ref{LemmaFactorization} and \ref{onedim}, we obtain the following result.

\medskip

\begin{Prop} \label{ConstructQuiver}

Let $\T$ be a cluster-tilting subcategory of $\mathscr{C}(Q)$. If $M, N$ are objects in $\mathscr{F}(Q) \cap \mathcal{T}$, then $Q_{_\T}$ has an arrow $M \to N$ if and only if $a_{\hspace{-1pt}_N}$ covers $a_{\hspace{-0.5pt}_M}$ with respect to $\vdash$.

\end{Prop}

\medskip

If $\mathfrak{p}, \mathfrak{q}$ are two marked points on the same boundary line of $\B_\infty$, then the line segment  between $\mathfrak{p}$ and $\mathfrak{q}$ is called a \emph{boundary segment}. The following is an easy consequence of Proposition \ref{ConstructQuiver}.

\medskip

\begin{Cor} \label{LemmaOrderArcs}

Let $\T$ be a cluster-tilting subcategory of $\mathscr{C}(Q)$. If $M,N \in \mathscr{F}(Q) \cap \T$, then $M, N$ lie in the same connected component of $Q_\T$ if and only if there exists a simple closed curve $S$ in $\B_\infty$ which is the composite of some arcs of ${\rm arc}(\T)$ and possibly some boundary segments such that the region enclosed by $S$ contains $a_M,a_N$ and at most finitely many arcs.
\end{Cor}

\medskip

We conclude this section with the following statement.

\medskip

\begin{Prop}\label{cpt-cts}

Let $\T$ be a cluster-tilting subcategory of $\mathscr{C}(Q)$. If ${\rm arc}(\T)$ has $m$ full fountains and $n$ non-full fountains, then $Q_{_\T}$ has $2m+ n+1$ connected components.
In particular, $Q_{\T}$ is connected if and only if ${\rm arc}(\T)$ has no fountain, or equivalently, every marked point in $\B_\infty$ is an endpoint of at most finitely many arcs of ${\rm arc}(\T)$.

\end{Prop}

\noindent{\it Proof.} If ${\rm arc}(\T)$ has no fountain then, by Corollary \ref{LemmaOrderArcs}, any pair of objects $M,N$ of $\mathscr{F}(Q) \cap \T$ lie in the same connected component of $Q_\T$. Therefore, we may assume that ${\rm arc}(\T)$ has at least one fountain.

Suppose first that  ${\rm arc}(\T)$ has either a full fountain (which is then the unique fountain) or two non-full fountains (that is, one left fountain and one right fountain). We shall construct non-empty subsets $\Sa_1, \Sa_2, \Sa_3$ of ${\rm arc}(\T)$ such that ${\rm arc}(\T)$ is a disjoint union of $\Sa_1, \Sa_2, \Sa_3$. Moreover, for $M,N \in \mathscr{F}(Q) \cap \T$, the arcs $a_M, a_N$ lie in the same connected component of $Q_\T$ if and only if $a_M, a_N \in \Sa_i$ for some $1\le i\le 3$.

For the first case, assume that there exists a full fountain. We may assume that the fountain base is $\ml_p$ for some integer $p$. Let $\Sa_1$ be the set of non-upper arcs of ${\rm arc}(\T)$, and $\Sa_2$ be the set of arcs of ${\rm arc}(\T)$ of form $[\ml_i, \ml_j]$ with $i,j \le p$, and $\Sa_3$ be the set of arcs of ${\rm arc}(\T)$ of form $[\ml_i, \ml_j]$ with $i,j \ge p$. By Corollary \ref{LemmaOrderArcs}, $\Sa_1, \Sa_2, \Sa_3$ satisfy the desired properties.

For the second case, assume that there exists a left fountain and a right fountain whose fountain bases lie on the same boundary line of $\B_\infty$. We may assume that these fountain bases are $\ml_p, \ml_q$ with $p < q$. In this case, let $\Sa_1$ be the set of arcs of ${\rm arc}(\T)$ which are connecting arcs, lower arcs, or upper arcs of form $[\ml_i, \ml_j]$ with $p \le i,j \le q$, and $\Sa_2$ be the set of arcs of ${\rm arc}(\T)$ of form $[\ml_i, \ml_j]$ with $i, j \le p$, and $\Sa_3$ be the set of arcs of ${\rm arc}(\T)$ of form $[\ml_i, \ml_j]$ with $i,j \ge q$. By Corollary \ref{LemmaOrderArcs}, $\Sa_1, \Sa_2, \Sa_3$ satisfy the desired properties.

For the third case, assume that there exists a left fountain and a right fountain whose fountain bases do not lie on the same boundary line of $\B_\infty$. We may assume that these fountain bases are $\ml_p, \mr_q$ with $\ml_p$ a left fountain base and $\mr_q$ a right fountain base. In this case, let $\Sa_1$ be the arcs of ${\rm arc}(\T)$ which are connecting arcs, lower arcs of form $[\mr_i, \mr_j]$ with $i,j \ge q$ or upper arcs of form $[\ml_i, \ml_j]$ with $i, j \ge p$. Let $\Sa_2$ be the arcs of ${\rm arc}(\T)$ of form $[\ml_i, \ml_j]$ with $i,j \le p$, and $\Sa_3$ be the arcs of ${\rm arc}(\T)$ of form $[\mr_i, \mr_j]$ with $i,j \le q$. By Corollary \ref{LemmaOrderArcs}, $\Sa_1, \Sa_2, \Sa_3$ satisfy the desired properties.

It remains to consider the case where ${\rm arc}(\T)$ has a unique fountain, which is not a full fountain. We may assume that $\ml_p$ is the fountain base, which is a left ${\rm arc}(\T)$-fountain base. In this case, let $\Sa_1$ be the set of arcs of ${\rm arc}(\T)$ of form $[\ml_i, \ml_j]$ with $i,j \le p$, and $\Sa_2$ be the set of other arcs of ${\rm arc}(\T)$. By Corollary \ref{LemmaOrderArcs}, $\Sa_1, \Sa_2$ satisfy the desired properties. The proof of the proposition is completed.

\end{document}